\documentclass[11pt]{article}
\usepackage{amssymb,amsmath,color,graphicx}
\usepackage{trfsigns}

\def\eins{\mbox{1\hskip-0.24em l}}
\def\T{^{\sf T}}
\def\mT{^{\sf -T}}
\def\ones{\mbox{\normalfont{1\hskip-0.24em l}}}

\newcommand{\R}{ {\mathbb R} }
\newcommand{\N}{ {\mathbb N} }

\newcommand{\MM}{ {\mathbb M} }
\newcommand{\diag}{\,\mbox{diag}}
\newcommand{\cc}{{\bf c}}
\newcommand{\CC}{{\cal C}}
\newcommand{\LL}{{\cal L}}
\newcommand{\PP}{{\cal P}}
\newcommand{\QQ}{{\cal Q}}
\newcommand{\ZZ}{{\cal Z}}
\DeclareMathOperator*{\argmin}{arg\,min}
\newcommand{\qed}{\qquad\mbox{$\square$}}
\newcommand{\gdw}{\ \iff\ }

\newtheorem{remark}{Remark}[section]
\newtheorem{theorem}{Theorem}[section]
\newtheorem{lemma}{Lemma}[section]
\newtheorem{example}{Example}[section]

\begin{document}
\title{Implicit Peer Triplets in Gradient-Based Solution Algorithms for ODE Constrained Optimal Control}
\author{Jens Lang \\
{\small \it Technical University Darmstadt,
Department of Mathematics} \\
{\small \it Dolivostra{\ss}e 15, 64293 Darmstadt, Germany}\\
{\small lang@mathematik.tu-darmstadt.de} \\ \\
Bernhard A. Schmitt \\
{\small \it Philipps-Universit\"at Marburg,
Department of Mathematics,}\\
{\small \it Hans-Meerwein-Stra{\ss}e 6, 35043 Marburg, Germany} \\
{\small schmitt@mathematik.uni-marburg.de}}
\maketitle

\begin{abstract}
It is common practice to apply gradient-based optimization algorithms to
numerically solve large-scale ODE constrained optimal control problems.
Gradients of the objective function are most efficiently computed by
approximate adjoint variables. High accuracy with moderate computing time
can be achieved by such time integration methods that satisfy a
sufficiently large number of adjoint order conditions and supply gradients
with higher orders of consistency. In this paper, we upgrade our former implicit
two-step Peer triplets constructed in [Algorithms, 15:310, 2022] to
meet those new requirements. Since Peer methods use several stages of the
same high stage order, a decisive advantage is their lack of order
reduction as for semi-discretized PDE problems with boundary control.
Additional order conditions for the control and certain positivity
requirements now intensify the demands on the Peer triplet.
We discuss the construction of $4$-stage methods with order pairs $(3,3)$
and $(4,3)$ in detail and provide three Peer triplets of practical interest.
We prove convergence of order $s-1$, at least, for $s$-stage methods if state, 
adjoint and control satisfy the corresponding order conditions.
Numerical tests show the expected order of convergence for the new Peer triplets.
\end{abstract}

\noindent{\em Key words.} Implicit Peer two-step methods, nonlinear optimal control,
gradient-based optimization, first-discretize-then-optimize, discrete adjoints

\newpage

\section{Introduction}
The numerical solution of optimal control problems governed by time-dependent
differential equations is still a challenging task in designing and
analyzing higher-order time integrators. An essential solution
strategy is the so called \textit{first-discretize-then-optimize} approach, where
the continuous control problem is first discretized into a nonlinear programming
problem which is then solved by state-of-the-art gradient-based optimization algorithm.
Nowadays this direct approach is the most commonly used method due to its easy
applicability and robustness.
Consistent gradients of the objective function are derived from control, state and adjoint variables given by first-order necessary conditions of Karush-Kuhn-Tucker type, and are used in
an iterative minimization algorithm to calculate the wanted optimal control.
In this solution strategy, the unique gain of higher-order time integrators is twofold: increasing the efficiency in computing time by the use of larger and fewer time steps and, even more important for large-scale problems, thus reducing at the same time the memory requirement caused by the necessity to store all variables for the computation of gradients.
\par
There are one-step as well as multistep time integrators in common use to solve
ODE constrained optimal control problems. Symplectic Runge-Kutta methods \cite{BonnansLaurentVarin2006,LiuFrank2021,SanzSerna2016} and backward differentiation formulas \cite{AlbiHertyPareschi2019,BeigelMommerWirschingBock2014} are prominent classes, but also partitioned and implicit-explicit Runge–Kutta methods \cite{HertyPareschiSteffensen2013,MatsudaMiyatake2021} and explicit stabilized Runge-Kutta-Chebyshev methods \cite{AlmuslimaniVilmart2021} have been proposed.
However, fully implicit one-step methods often request the solution of large systems
of coupled stages and might suffer from serious order reduction due to their lower stage
order. This is especially the case when they are applied to semi-discretized PDEs with
general time-dependent boundary conditions arising from boundary control problems
\cite{OstermannRoche1992,LangSchmitt2023a}. Matlab codes for \cite{LangSchmitt2023a} are available at \cite{LangSchmitt2023b}. In general, further consistency conditions have to be satisfied \cite{Hager2000,LangVerwer2013} in order the achieve a higher
classical order for the discrete adjoint variables. Multistep methods avoid order reduction and
have a simple structure, but higher order comes with restricted stability properties and adjoint initialization steps are usually inconsistent approximations which may lead to severe losses of accuracy.
Moreover, the appropriate approximation of initial values and its structural consequences for the
adjoint variables are further unsolved inherent difficulties that have limited
the application of higher-order multistep methods for optimal control problems in a first-discretize-then-optimize solution strategy.
\par
Recently, we have proposed a new class of implicit two-step Peer triplets that aggregate the attractive properties of one- and multistep methods through the use of several stages of one and the same high order and at the same time avoid their deficiencies with the aid of a two-step form
\cite{LangSchmitt2022,LangSchmitt2022b,LangSchmitt2023a}. The incorporation of
different but matching start and end steps increases the flexibility of Peer methods for the solution of optimal control problems, especially also allowing higher-order approximations of the control, adjoint variables and the gradient of the objective function.
\par
The class of $s$-stage implicit two-step Peer methods was introduced in \cite{SchmittWeiner2004} in linearly implicit form to solve stiff ODEs of the form $y'(t)=f(y(t))$, $y_0=y(0)$, $y\in\R^m$.
Later, the methods were simplified as implicit two-step schemes
\begin{align}\label{peer-std}
Y_n =&\; (Q\otimes I_m)Y_{n-1} + h (R\otimes I_m) F(Y_n),\; n=0,1,\ldots,
\end{align}
where the constant time step $h>0$ will be considered here.
An application of \eqref{peer-std} to large scale problems with Krylov solvers was discussed in  \cite{BeckWeinerPodhaiskySchmitt2012}.
The stage solutions $Y_n=(Y_{ni})_{i=1}^s$ are
approximations of $(y(t_n+c_ih))_{i=0}^s$ with equal accuracy and stability
properties, which motivates the attribute \textit{peer} and is the key for avoiding order reduction.
The off-step nodes $c_1,\ldots,c_s$ are associated with the interval $[0,1]$ but some may lie outside.
In \eqref{peer-std}, $R\in\R^{s\times s}$ is lower triangular
and invertible, $Q\in\R^{s\times s}$, $I_m\in\R^{m\times m}$ is the identity matrix, and $F(Y_n)=(f(Y_{ni}))_{i=1}^s$. Note that the
stages $Y_{ni}$ can be successively computed for $i=1,\ldots,s$ due to the
triangular structure of $R$. Several variants of Peer methods have been
developed and successfully applied to a broad class of differential
equations, e.g. \cite{GerischLangPodhaiskyWeiner2009,JebensKnothWeiner2009,MassaNoventaLoriniBassiGhidoni2018,PodhaiskyWeinerSchmitt2005,SchmittWeinerBeck2013,SchneiderLangWeiner2021}.
\par
The application of Peer methods to optimal control problems requires a
couple of modifications. A first direct attempt in \cite{SchroederLangWeiner2014} was unsatisfactory, mainly due to the restricted, first-order approximation of the
adjoint variables. In \cite{LangSchmitt2022}, we found that the general redundant
formulation of the above standard Peer method,
\begin{align}\label{peer-trf}
(A\otimes I_m)Y_n =&\; (B\otimes I_m)Y_{n-1} + h (K\otimes I_m) F(Y_n),\; n=0,1,\ldots,
\end{align}
with invertible lower triangular matrix $A\in\R^{s\times s}$ and diagonal matrix $K\in\R^{s\times s}$ is admittedly equivalent in terms of the state variables to \eqref{peer-std}, but surprisingly not for the adjoint variables with the same coefficients $Q,R$.
The additional degrees of freedom given by $K$ together with a careful design of  start and
end methods with different coefficient matrices laid the foundation for improved Peer triplets with higher-order convergence. Triplets with good stability properties could be found
\cite{LangSchmitt2022b}. Formulation \eqref{peer-trf} will also be the starting point in this paper. In contrast to our former approach in \cite{LangSchmitt2022,LangSchmitt2022b}, where
the control $u(t)$ has been eliminated and a boundary value problem has been solved, we now compute
the optimal control in an iterative procedure, making use of gradients of the objective function.
\par
An important advantage of the approximations of the state
$Y_{ni}\approx y(t_n+c_ih)$ and the adjoint multipliers $P_{ni}\approx p(t_n+c_ih)$ in
the discrete time points $t_n+c_ih$
(see the next Chapter for the details of the notation) with equal
high accuracy is the opportunity to simply use the discrete control variables
\begin{align}\label{peer-opt}
U_{ni} = \Phi(Y_{ni},P_{ni})\approx u(t_n+c_ih),\; i=1,\ldots,s,\,n=0,1,\ldots,
\end{align}
in a gradient-based optimization algorithm. Although, the function
$\Phi$ is only implicitly given, the higher order of $Y_{ni}$ and
$P_{ni}$ is directly transferable to the control vector $U_{ni}$.
Interpolation in time is easily realizable, the data structure keeps simple.
However, additional order conditions for the control derived from
\eqref{peer-opt} and positivity requirements for column sums in the
matrix triplet $(K_0,K,K_N)$, where $K_0$ and $K_N$ are the matrices
for the start and end method, intensify the demands on the Peer methods in
the triplet. The arising bottlenecks in the design caused by stronger
entanglement of all matrices must be resolved by a more sophisticated analysis.
\par
The paper is organized as follows. In Chapter \ref{SecOpt}, we formulate
the optimal control problem and define its discretization. The gradient of the
cost function is derived in Chapter \ref{SecGradM}. Order conditions and
their algebraic consequences are discussed in Chapter \ref{SecOrdCond}. Two classes
of four-stage Peer triplets are studied in Chapter \ref{SP4o43} and Chapter \ref{SP4o33}
and three triplets of practical interest are constructed. In Chapter \ref{SGlobErr},
we give a detailed convergence proof for unconstrained controls.
Numerical examples collected in Chapter \ref{SecNum} illustrate the theoretical findings.
The paper concludes with a summary in Chapter \ref{SecSum}.
\section{The optimal control problem and its discretization}\label{SecOpt}
We are interested in the numerical solution of
the following ODE-constrained nonlinear optimal control problem:
\begin{align}
\mbox{minimize } \CC\big(y(T)\big) \label{OCprob_objfunc} &\\
\mbox{subject to } y'(t) =& \,f\big(y(t),u(t)\big),\quad
u(t)\in U_{ad},\;t\in(0,T], \label{OCprob_ODE}\\
y(0) =& \,y_0, \label{OCprob_ODEinit}
\end{align}
where the state $y(t)\in\R^m$, the control $u(t)\in\R^d$,
$f: \R^m\times\R^d\mapsto\R^m$, the objective function $\CC: \R^m\mapsto\R$, and the set of admissible controls $U_{ad}\subset\R^d$ is closed and convex. Introducing for any $u\in U_{ad}$ the normal cone mapping
\begin{align*}
N_U(u) =&\, \{ w\in\R^d: w^T(v-u)\le 0 \mbox{ for all } v\in U_{ad}\},
\end{align*}
the first-order optimality conditions read \cite{Hager2000,Troutman1996}
\begin{align}
y'(t) =& \,f\big(y(t),u(t)\big),\quad t\in(0,T],\quad y(0)=y_0, \label{KKT_state}\\
p'(t) =& \,-\nabla_y f\big(y(t),u(t)\big)\T p(t),\quad t\in[0,T),
\quad p(T)=\nabla_y {\CC}\big(y(T)\big)\T, \nonumber\\
& \,-\nabla_u f\big(y(t),u(t)\big)\T p(t) \in N_U\big(u(t)\big),\quad t\in[0,T]. \label{KKT_ctr}
\end{align}
Under appropriate regularity conditions, there exists a local solution
$(y^\star,u^\star)$ of the optimal control problem
\eqref{OCprob_objfunc}-\eqref{OCprob_ODEinit} and a Lagrange multiplier
$p^\star$ such that the first-order optimality conditions
\eqref{KKT_state}-\eqref{KKT_ctr} are necessarily satisfied at $(y^\star,p^\star,u^\star)$.
If, in addition, the Hamiltonian
$H(y,p,u):=p\T f(y,u)$ satisfies a coercivity assumption, then these
conditions are also sufficient \cite{Hager2000}. The control
uniqueness property introduced in \cite{Hager2000} yields the existence
of a locally unique minimizer $u=u(\hat{y},\hat{p})$ of the Hamiltonian over all
$u\in U_{ad}$, if $(\hat{y},\hat{p})$ is sufficiently close to $(y^\star,p^\star)$.
\par
Many other optimal control problems can be transformed to the Mayer form
$\CC(y(T))$ which only uses terminal solutions. For example, terms given in
the Lagrange form
\begin{align}\label{objfunc-lagrange}
\CC_L(y,u) :=&\;\int_{0}^{T}l(y(t),u(t))\,dt
\end{align}
can be equivalently reduced to the Mayer form by adding a
new differential equation $y_{m+1}'(t)=l(y(t),u(t))$ and initial
values $y_{m+1}(0)=0$ to the constraints. Then \eqref{objfunc-lagrange}
simply reduces to $y_{m+1}(T)$.
\par
On a time grid $\{t_0,\ldots,t_N\}\subset [0,T]$ with fixed step size length $h=t_{n+1}-t_n$ Peer methods use $s$ stage approximations
$Y_{ni}\approx y(t_{ni})$ and $U_{ni}\approx u(t_{ni})$ per time step at points $t_{ni}=t_n+c_ih$, $i=1,\ldots,s,$ associated with fixed nodes $c_1,\ldots,c_s$.
All $s$ stages share the same properties like a common high stage order equal to the global order, preventing order reduction.
Applying the two-step Peer method for $n\ge 1$ and an exceptional starting step for $n=0$ to the problem (\ref{OCprob_ODE})--(\ref{OCprob_ODEinit}),
we get the discrete constraint nonlinear optimal
control problem
\begin{align}
\mbox{minimize } \CC\big(y_h(T)\big) \nonumber &\\[1mm]
\mbox{subject to } A_0 Y_0=&\,a\otimes y_0+hK_0F(Y_0,U_0),\label{OCprob_peer_init}\\[1mm]
A_n Y_n=&\,B_nY_{n-1}+hK_nF(Y_n,U_n),\ n=1,\ldots,N,\label{OCprob_peer}
\end{align}
with long vectors $Y_n=(Y_{ni})_{i=1}^s\in\R^{sm}$,
$U_n=(U_{ni})_{i=1}^s\in\R^{sd}$, and $F(Y_n,U_n)=\big(f(Y_{ni},U_{ni})\big)_{i=1}^s$.
Further, $y_h(T)=(w\T\otimes I)Y_N\approx y(T)$, $a,w\in\R^s$,
$A_n,B_n,K_n\in\R^{s\times s}$, and
$I\in\R^{m\times m}$ being the identity matrix.
As a change to the introduction,
we will use the same symbol for a coefficient matrix like $A$ and its Kronecker product
$A\otimes I$ as a mapping from
the space $\R^{sm}$ to itself. Throughout the paper, $e_i$ denotes the $i$-th
cardinal basis vector and $\eins:=(1,\ldots,1)\T\in\R^s$, sometimes with an additional 
index indicating a different space dimension.
\par
On each subinterval $[t_n,t_{n+1}]$, Peer methods may be
defined by three coefficient matrices $A_n,B_n,K_n$, where $A_n$ is
assumed to be non-singular.
For practical reasons, this general version will not be used. We choose a
fixed Peer method $(A_n,B_n,K_n)\equiv (A,B,K)$, $n\!=\!1,\ldots,N-1$, in the inner
grid points with lower triangular non-singular $A$, which allows a consecutive computation of
the solution vectors $Y_{ni}$, $i\!=\!1,\ldots,s$, in \eqref{OCprob_peer_init},
\eqref{OCprob_peer}. Exceptional coefficients
$(A_0,K_0)$ and $(A_N,B_N,K_N)$ in the first and last forward steps are taken
to allow for a better approximation in the initial step and of the end conditions.
\par
The first order optimality conditions now read
\begin{align}
A_0 Y_0=&\,a\otimes y_0+hK_0F(Y_0,U_0),\label{KKT_state_peer_init}\\[1mm]
A_n Y_n=&\,B_nY_{n-1}+hK_nF(Y_n,U_n),\ n=1,\ldots,N,\label{KKT_state_peer}\\[1mm]
A_N\T P_N=&\,w\otimes p_h(T)+h\nabla_YF(Y_N,U_N)\T K_N\T P_N,\label{KKT_adj_peer_init}\\[1mm]
A_n\T P_n=&\,B_{n+1}\T P_{n+1}+h\nabla_YF(Y_n,U_n)\T K_n\T P_n,\ 0\le n\le N-1,\label{KKT_adj_peer}\\[1mm]
&\,-h\nabla_UF(Y_n,U_n)\T K_n\T P_n \in N_{U^s}\big(U_n\big),\ 0 \le n\le N.\label{KKT_ctr_peer}
\end{align}
Here, $p_h(T)=\nabla_y \CC\big(y_h(T)\big)\T$ and the Jacobians of $F$ are block diagonal matrices $\nabla_YF(Y_n,U_n)=\diag_i\big(\nabla_{Y_{ni}}f(Y_{ni},U_{ni})\big)$ and
$\nabla_UF(Y_n,U_n)=\diag_i\big(\nabla_{U_{ni}}f(Y_{ni},U_{ni})\big)$.
The generalized normal cone mapping $N_{U^s}\big(U_n\big)$ is defined by
\begin{align*}
N_{U^s}(u) =&\, \left\{ w\in\R^{sd}: w^T(v-u)\le 0 \mbox{ for all } v\in U_{ad}^s\subset\R^{sd}\right\}.
\end{align*}
If  $K_n:=(\kappa_{ij}^{[n]})$ is diagonal and $\kappa_{ii}^{[n]}=0$, then \eqref{KKT_ctr_peer} is automatically 
satisfied for stage $i\in\{1,\ldots,s\}$.
Assuming otherwise $k_{ni}:=e_i^TK_n^T\ones\not=0$, and defining
\begin{align*}
Q_{ni} := \frac{1}{k_{ni}}\,\sum_j \kappa_{ji}^{[n]}P_{nj},
\ 0\le n\le N,\ i=1,\ldots,s,
\end{align*}
(\ref{KKT_ctr_peer}) can be equivalently reformulated as
\begin{align}\label{use_qni}
-k_{ni}\,\nabla_{U_{ni}}f(Y_{ni},U_{ni})\T Q_{ni} \in N_{U}\big(U_{ni}\big),\
0 \le n\le N,\ i=1,\ldots,s.
\end{align}
Note that $Q_{ni}=P_{ni}$, if the matrix $K_n$ is diagonal.
A severe new restriction on the Peer triplet comes from the need to preserve the correct sign in \eqref{use_qni} requiring that $k_{ni}>0$.
Then, we can divide by it and the control
uniqueness property guarantees the existence
of a local minimizer $U_{ni}$ of the Hamiltonian $H(Y_{ni},Q_{ni},U)$ over all
$U\in U_{ad}$ since $Q_{ni}$
can be seen as an approximation to the multiplier $P_{ni}\approx p(t_n+c_ih)$.
Such positivity conditions also arise in the context of classical Runge-Kutta methods
or W-methods, see e.g. \cite[Theorem 2.1]{Hager2000} and \cite[Chapter 5.2]{LangVerwer2013}.
\par
The need to sacrifice the triangular resp. diagonal form of the matrix coefficients $A_n,K_n$ in the boundary steps comes from the fact that the starting steps \eqref{KKT_state_peer_init}, and backwards \eqref{KKT_adj_peer_init} are single-step methods with $s$ outputs.
With a triangular form of $A_0,A_N$ and $K_0,K_N$ their first stages (backward for $n=N$) would represent simple implicit Euler steps with a local order limited to 2,
see Section~5 in \cite{LangSchmitt2022} for a discussion.
\section{The gradient of the cost function}\label{SecGradM}
We first introduce the vector of control values for the entire interval $[0,T]$
\begin{align*}
U =&\,(U_{01}\T,\ldots,U_{0s}\T,U_{11}\T,\ldots,U_{Ns}\T)\T\in\R^{sd(N+1)}
\end{align*}
and let $\CC(U):=\CC(y_h(U))$ be the cost function associated with these controls. The
first order system (\ref{KKT_state_peer_init})--(\ref{KKT_ctr_peer}) provides a
convenient way to compute the gradient of $\CC(U)$ with respect to $U$. Following the
approach in \cite{HagerRostamian1987}, we find
\begin{align}\label{gradient_C}
\nabla_{U_{ni}}\CC(U) =&\,h\nabla_{U_{ni}}f(Y_{ni},U_{ni})\T (e_i\T K_n\T\otimes I) P_n,\
0 \le n\le N,\ i=1,\ldots,s,
\end{align}
where all approximations $(Y_n,P_n)$ are computable by a forward-backward marching scheme.
The state variables $Y_n$ are obtained from the discrete state equations (\ref{KKT_state_peer_init})--(\ref{KKT_state_peer}) for $n=0,\ldots,N$, using the given values of the control vector $U$.
Then, using the updated values $Y_n$, one computes
$p_h(T)=\nabla_y \CC\big(y_h(T)\big)\T$ with $y_h(T)=(w\T\otimes I)Y_N$ before marching the steps (\ref{KKT_adj_peer_init})--(\ref{KKT_adj_peer}) backwards for $n=N,\ldots,0$, solving the discrete costate equations
for all $P_n$.
\par
The gradients from \eqref{gradient_C} can now be employed in gradient-based optimization algorithms which have
been developed extensively since the 1950s.
Many good algorithms are now available to
solve nonlinear optimization problems in an iterative procedure
\begin{align}\label{nopt_iter}
U^{(k+1)} :=&\, U^{(k)}+\triangle U^{(k)},\quad k=0,1,\ldots
\end{align}
starting from an initial estimate $U^{(0)}$ for the control vector. Evaluating the objective
function, its gradient and, in some cases, its Hessian, an efficient update $\triangle U^{(k)}$ of the control can be computed.
Based on the principle (\ref{nopt_iter}),
several good algorithms have been implemented in commercial software packages like
{\sc Matlab}, {\sc Mathematica}, and others. We will use the {\sc Matlab} routine
\texttt{fmincon} in our numerical experiments. It offers several optimization algorithms
including \texttt{interior-point} \cite{ByrdGilbertNocedal2000} and \texttt{trust-region-reflective} \cite{ColemanLi1994} for large-scale sparse
problems with continuous objective function and first derivatives.
\par
Since the optimal control $u(t)$ minimizes the Hamiltonian $H(y,p,u)=p\T f(y,u)$, we may compute an improved approximation of the control by the following minimum principle:
\begin{align}\label{ctrpp}
U^\ddagger_{ni} = \argmin_{U\in U_{ad}} H(Y_{ni},P_{ni},U),\quad
0 \le n\le N,\ i=1,\ldots,s,
\end{align}
if $Y_{ni}$ or $P_{ni}$ are approximations of higher-order. We note that the function
$\Phi$ in \eqref{peer-opt} provides the solution in \eqref{ctrpp}, when $H$ is replaced by
its discrete approximation defined by the Peer triplet.
\section{Order conditions for the Peer triplet in the unconstrained case}\label{SecOrdCond}
We recall the conditions for local order $r\le s$ for the forward schemes
and order $q\le s$ for the adjoint schemes, see \cite{LangSchmitt2022,LangSchmitt2022b}.
These conditions use the Vandermonde matrices $V_q:=(\eins,\cc,\cc^2,\ldots,\cc^{q-1})\in\R^{s\times q}$
with the column vector of nodes
$\cc=(c_i)_{i=1}^s\in\R^s$, and the non-singular Pascal matrix $\PP_q=\big({j-1\choose i-1}\big)=\exp(\tilde E_q)\in\R^{q\times q}$, where $\tilde E=\big(i\delta_{i+1,j}\big)\in\R^{q\times q}$ is nilpotent.
There are 5 conditions for the forward scheme and its adjoint method:
\begin{align}
A_0V_r =&\, ae_1\T +K_0V_r\tilde E_r,\ n=0,\label{OBvStrt}\\[1mm]
A_nV_r =&\, B_nV_r\PP_r^{-1}+K_nV_r\tilde E_r,\ 1\le n\le N, \label{OBvStd}\\[1mm]
w\T V_r =&\,\eins\T,\label{OBvEnd}\\[1mm]
A_n\T V_q=&\,B_{n+1}\T V_q\PP_q-K_n\T V_q\tilde E_q,\ 0\le n\le N-1,\label{OBaStd}\\[1mm]
A_N\T V_q=&\,w\eins\T -K_N\T V_q\tilde E_q,\ n=N. \label{OBaEnd}
\end{align}
We remind that the coefficient matrices from interior grid intervals belong to a standard scheme $(A_n,B_n,K_n)\equiv(A,B,K),\,1\le n\le N-1$.
The whole triplet consists of 8 coefficient matrices $(A_0,K_0),\,(A,B,K),\,(A_N,B_N,K_N)$.
\par
Next, we focus on the new optimality condition (\ref{KKT_ctr_peer})
in the unconstrained case with $N_U=\{0\}$.
It reads stage-wise
\begin{align}
\label{KKT_ctr_peer_all}
\nabla_u f(Y_{nj},U_{nj})\T\sum_{i=1}^sP_{ni}\kappa_{ij}^{[n]}=0,\ j=1,\ldots,s,\ 0\le n\le N.
\end{align}
Order conditions are obtained by Taylor expansions, where approximations are
replaced by exact solutions $(y^\star(t_n+c_ih),p^\star(t_n+c_ih),u^\star(t_n+c_ih))$ and the
(continuous) optimality condition
\begin{align}\label{KKT_ctr_ode}
\nabla_u f\big(y^\star(t),u^\star(t)\big)\T p^\star(t) =&\; 0,\ t\in [0,T]
\end{align}
is used.
Defining the partial sums $\exp_q(z):=\sum_{j=0,\ldots,q-1}z^j/j!$ with $q$ terms,
Taylor's theorem for the expansion of a smooth function $v(t)$, $v\in C^q[0,T]$, at $t_{nj}:=t_n+c_jh$ may be written as
\begin{align*}
v(t_{n}+c_ih) =&\;\exp_q((c_i-c_j)z)v|_{t=t_{nj}}+O(h^q),\; z:=h
\frac{d}{dt},
\end{align*}
with some slight abuse of notation.
Then, the corresponding expansion of the residuals in (\ref{KKT_ctr_peer_all})
for order $q+1$ gives
\begin{align}\label{NBres}
 \nabla_u f\big(y^\star(t_{nj}),u^\star(t_{nj})\big)\T
 \sum_{i=1}^s\kappa_{ij}^{[n]}\exp_q\big((c_i-c_j)z\big)p^\star(t_{nj})
 \stackrel!=O(z^q),\ j=1,\ldots,s.
\end{align}
\begin{lemma}\label{Lknspf}
Let the solution $p^\star$ be smooth, $p^\star\!\in\! C^q[0,T]$, and $C\!:=\!\diag(c_1,\ldots,c_s)$ the diagonal matrix containing the nodes and assume that
\begin{align}\label{add_cond}
(\cc^{l-1})\T K_n - \ones\T K_nC^{l-1} =&\;0,\ l=2,\ldots,q,
\end{align}
for all $n=0,\ldots,N$.
Then, for these $n$ and $j=1,\ldots,s$ it holds with $k_{nj}$ from \eqref{use_qni}
\begin{align}\label{FeKP}
 &\sum_{i=1}^s\kappa_{ij}^{[n]}p^\star(t_{ni})- k_{nj} p^\star(t_{nj})=O(h^q),
\\\label{DeftauU}
 \tau_{nj}^U:=&\,\nabla_u f(y^\star(t_{nj}),u^\star(t_{nj}))\T\sum_{i=1}^sp^\star(t_{ni})
 \kappa_{ij}^{[n]}=O(h^q).
\end{align}
\end{lemma}
{\it Proof:}
The assumption \eqref{add_cond} is obtained from \eqref{NBres} by removing mixed powers in the conditions $\sum_{i=1}^s\kappa_{ij}^{[n]}(c_i-c_j)^{l-1}=0$ with the corresponding equations for lower degrees.
In fact, these condition prove the stronger version \eqref{FeKP} which will be needed below.
Of course, \eqref{DeftauU} is a simple consequence due to \eqref{KKT_ctr_ode}.
\qed
\par
For the standard method with a diagonal matrix $K_n\equiv K$, $1\le n<N$,
the condition \eqref{OBaStd} is sufficient for adjoint local order $q$ since \eqref{add_cond} is satisfied trivially.
However, for the more general matrices $K_0,K_N$ required in the first
and last forward steps, it has been
shown in \cite[Chapter 2.2.4]{LangSchmitt2022b} that
additional conditions have to be satisfied due to an unfamiliar form of one-leg-type applied to the linear adjoint equation
$p'=-J(t)p$ with $J(t)=\nabla_yf(y(t),u(t))\T$.
These conditions are now covered by \eqref{add_cond} for $l\!=\!2$.
In our present context with unknown control, the additional
constraint equation \eqref{KKT_ctr_peer_all} sharpens these requirements and leads to
much stronger restrictions on the design of the whole Peer triplet.
As discussed in connection with \eqref{use_qni}, we also require positive column sums in the boundary steps and non-negative ones in the standard scheme,
\begin{align}\label{pos-cond}
\eins\T K_0 > 0\T,\ \eins\T K_N > 0\T,\ \eins\T K\ge 0\T.
\end{align}
In the special case of a Peer triplet of FSAL type as constructed in Section
\ref{SecFSAL}, we allow $e_1\T K_0=0\T$ since the corresponding control component $U_{01}$ can be eliminated.
We note that even in the unconstrained case, $N_U=\{0\}$, positivity \eqref{pos-cond} is required in order to preserve the positive definiteness of the Hesse matrix.
\subsection{Combined conditions for the order pair $(r,q)$}
With $r,q\le s$ the full set of conditions may
lead to practical problems for deriving formal solutions with the aid of
algebraic manipulation software due to huge algebraic expressions.
Possible alternatives like numerical search procedures will suffer from
the large dimension of the search space consisting of the entries of 8
coefficient matrices.
Fortunately, many of these parameters may be eliminated temporarily by
solving certain condensed necessary conditions first.
Afterwards, the full set \eqref{OBvStrt}--\eqref{OBaEnd} may be more easily
solved in decoupled form.
The combined conditions will be formulated with the aid of the linear operator
\begin{align}\label{LLqr}
 X\mapsto \LL_{q,r}(X):=\tilde E_q\T X+X\tilde E_r,\quad X\in\R^{q\times r}.
\end{align}
We note that the map $\LL_{q,r}$ is singular since $\tilde E_r$ is nilpotent
and $\tilde E_r e_1=0$ for any $r\in\N$.
Hence, the first entry of its image always vanishes
\begin{align}\label{LL11}
\big(\LL_{q,r}\big)_{11}=0.
\end{align}
The combined conditions are presented in the order in which they would be
applied in practice, with the standard scheme $(A,B,K)$ in the first place.
In all these conditions the matrix
\begin{align*}
 \QQ_{q,r}:=V_q\T BV_r\PP_r^{-1}
\end{align*}
plays a central role.
\begin{lemma}\label{LKombOB}
Let the matrices of the Peer triplet $(A_0,K_0)$, $(A,B,K)$, $(A_N,B_N,K_N)$ satisfy the order conditions \eqref{OBvStrt}--\eqref{OBaEnd} with $q\le r\le s$.
Then, also the following equations hold,
\begin{align}\label{cocos}
 \LL_{q,r}(V_q\T KV_r)=&\PP_q\T V_q\T BV_r-V_q\T BV_r\PP_r^{-1}
  =\PP_q\T \QQ_{q,r}\PP_r-\QQ_{q,r},\\\label{coco0}
 \LL_{q,r}(V_q\T K_0V_r)=&\PP_q\T V_q\T BV_r-V_q\T a e_1\T
  =\PP_q\T \QQ_{q,r}\PP_r-V_q\T a e_1\T,\\\label{cocoN}
 \LL_{q,r}(V_q\T K_NV_r)=&\eins_q\eins_r\T-V_q\T BV_r\PP_r^{-1}
  =\eins_q\eins_r\T-\QQ_{q,r}.
\end{align}
\end{lemma}
{\it Proof:}
Considering the cases $1\le n<N$ first, equation \eqref{OBvStd}
is multiplied by $V_q\T$ from the left and the transposed condition
\eqref{OBaStd} by $V_r$ from the right.
This gives
\begin{align*}
 V_q\T AV_r=&V_q\T BV_r\PP_r^{-1}+(V_q\T K V_r)\tilde E_r, \\
  V_q\T AV_r =&\PP_q\T V_q\T BV_r -\tilde E_q\T(V_q\T KV_r) .
\end{align*}
Subtracting both equations eliminates $A$ and yields \eqref{cocos}.
Equation \eqref{coco0} follows in the same way from \eqref{OBvStrt}
and \eqref{OBaStd} for $n=0$.
The end method has to satisfy three conditions.
Multiplying \eqref{OBvStd} for $n=N$ again from the left by $V_q\T$
and the transposed end condition \eqref{OBaEnd} by $V_r$ from the
right and combining both eliminates $A_N$ and gives
\begin{align}\label{cocoNN}
 \LL_{q,r}(V_q\T K_NV_r)=\eins_q\eins_r\T-V_q\T B_NV_r\PP_r^{-1}
 =\eins_q\eins_r\T-\QQ_{q,r}
\end{align}
since $V_r\T w=\eins$ by \eqref{OBvEnd}.
The third condition for $B_N$ is \eqref{OBaStd} with $n=N-1$.
It reduces to $B_N\T V_q\PP_q=A\T V_q+K\T V_q\tilde E_q=B\T V_q\PP_q$,
which means $V_q\T B_N=V_q\T B$.
Hence, the matrix $B_N$ in \eqref{cocoNN} may simply be replaced by $B$.
\qed
\par
Since the operator $\LL_{q,r}$ is singular, solutions for
\eqref{cocos}--\eqref{cocoN} exist for special right-hand sides only.
For instance, in \eqref{cocoN}, \eqref{coco0} the property \eqref{LL11} requires that
\begin{align}\label{eBe}
 1=e_1\T \QQ_{q,r}e_1=\eins\T B\eins\mbox{ and } \eins_q\T a=\eins\T B\eins=1.
\end{align}
Also the map $X\mapsto \PP_q\T X\PP_r-X$ in \eqref{cocos} is singular, but here,  \eqref{LL11} imposes no restrictions since we always have $e_1\T\big(\PP_q\T X\PP_r-X)e_1=0$.
However, many further restrictions are due to special structural properties of the matrices $V_q\T K_nV_r$ which are the arguments of $\LL_{q,r}$. The following algebraic background highlights the hidden
Hankel structure of certain matrices.
\par
A matrix $X=(x_{ij})\in\R^{q\times r}$ is said to have Hankel form, if its elements are constant along anti-diagonals, i.e. if $x_{ij}=\xi_{i+j-1}$ for $1\le i\le q$, $1\le j\le r$.
Some simple, probably well-known, properties are the following.
\begin{lemma}\label{LHankel}
a) If $K\in\R^{s\times s}$ is a diagonal matrix, then $V_q\T KV_r$ has Hankel form.\\[1mm]
b) Congruence transformations with Pascal matrices and the operator $\LL_{q,r}$ preserve Hankel form:
if $X\in\R^{q\times r}$ has Hankel form, then also $\PP_q\T X\PP_r$ and $\LL_{q,r}(X)$.\\[1mm]
c) The operator $\LL_{q,r}$ is homogeneous for multiplication with Pascal matrices:
\[ \PP_q\T\LL_{q,r}(X)=\LL_{q,r}(\PP_q\T X),\ \LL_{q,r}(X)\PP_r=\LL_{q,r}(X\PP_r),
 \quad X\in\R^{q\times r}.\]
\end{lemma}
{\it Proof:} a) For $K=\diag(\kappa_{ii})$, we get $(\cc^{i-1})\T K\cc^{j-1}=\sum_{\ell=1}^s\kappa_{\ell\ell}c_\ell^{i+j-2}$.
\\ b) For $X=(\xi_{i+j-1})$ the explicit expression
\[ e_i\T \PP_q\T X\PP_r e_j=\sum_{k=1}^{i+j-1}{i+j-2\choose k-1}\xi_k\]
should be well-known and is easily shown with the aid of the Vandermonde identity.
And from $\tilde E_q=(i\delta_{i+1,j})$ it follows that
\[ e_i\T(\tilde E_q\T X+X\tilde E_r)e_j=(i+j-2)\xi_{i+j-1},\]
where the factor $i+j-2=0$ leads to \eqref{LL11} for $i=j=1$.
\\ Assertion c) holds, since $\tilde E_q$ and $\PP_q=\exp(\tilde E_q)$ commute.
\qed
\subsubsection{Consequences of the one-leg-conditions}
For convenience, the trivial case $l=1$ is included in the one-leg conditions \eqref{add_cond} for adjoint order $q$.
Recalling the matrix $C:=\diag(c_i)$, these conditions may be written as
\begin{align}\label{olq}
 \left.\begin{array}{c}
 (\cc^{l-1})\T K_n=\eins\T K_nC^{l-1}\\
 l=1,\ldots,q
 \end{array}\right\}
 \gdw
  V_q\T K_n=\begin{pmatrix}
  \eins\T K_n\\\eins\T K_nC\\ \vdots\\\eins\T K_nC^{q-1} \end{pmatrix}
  \in\R^{q\times s}.
\end{align}
In particular, the matrix $V_q\T K_n$ depends on the row vector $\eins\T K_n$ of the column sums of $K_n$ only.
The consequences of \eqref{olq} on the matrices in Lemma~\ref{LKombOB} are far-reaching.
It will be seen that \eqref{olq} reduces the combined order conditions \eqref{coco0}, \eqref{cocoN} for the boundary methods to over-determined linear systems of the column sums $\eins\T K_0,\,\eins\T K_N$ alone.
This will lead to bottlenecks in the design of Peer triplets for $q+r>s+2$.
A first simple restriction is shown now.
\begin{lemma}\label{LKol}
If a matrix $K_n\in\R^{s\times s}$ satisfies \eqref{olq}, then the matrices $V_q\T K_nV_k$ and $\LL_{q,k}(V_q\T K_nV_k)$ have Hankel form for any $k\in\N$.
\end{lemma}
{\it Proof:} For $l\le q$ and $k\in\N$, Hankel form follows from \eqref{olq} by
\begin{align*}
 (\cc^{l-1})\T K_n\cc^{k-1}=\eins\T K_nC^{l-1}\cc^{k-1}=\eins\T K_n\cc^{l+k-2}.
\end{align*}
In fact, $V_q\T K_n V_k=V_q\T D_K V_k$ with the diagonal matrix $D_K$ containing the column sums of $K_n$.
The Hankel property of $\LL_{q,k}(V_q\T K_nV_k)$ then is a consequence of Lemma~\ref{LHankel}.
\qed
\par\noindent
\begin{remark}
By this Lemma, the conditions for the end method lead to severe restrictions on the standard method through equation \eqref{cocoN}.
Since its left hand side $\LL_{q,r}(V_r\T K_NV_r)$ is in Hankel form, also its right hand side $\eins_q\eins_r\T-\QQ_{q,r}$, has to be so restricting the shape of $B$ further.
In particular it means that
\begin{align}\label{HankelQ}
 \QQ_{q,r}=V_r\T BV_q\PP_q^{-1}\mbox{ has Hankel form}.
\end{align}
Then, by Lemma~\ref{LHankel}, also $\PP_q\T\QQ_{q,r}\PP_r$ and $\LL_{q,r}(V_q\T K_0V_r)$ for the starting method have Hankel form, which leaves only $V_q\T a=e_1$ for the slack variable in equation \eqref{coco0}.
\end{remark}
\begin{remark}
The Hankel form of all three matrices $V_q\T K_nV_r$ in Lemma~\ref{LKombOB} also prohibits the presence of certain kernel elements of $\LL_{q,r}$ in possible solutions.
In this Lemma there are obvious similarities between the equations for $K,K_0,K_N$ and for the method \texttt{AP4o43p} below, it was observed that
\[  V_3\T(K_0+K_N-K)V_4=\begin{pmatrix}
  1&\frac12&\frac13&\frac14\\[1mm]
  \frac12&\frac13&\frac14&\frac15\\[1mm]
  \frac13&\frac14&\frac15&\ast
 \end{pmatrix}
\]
is a Hilbert matrix with the exception of its last element.
\end{remark}
Due to the one-leg conditions \eqref{olq} the matrix equations \eqref{coco0}, \eqref{cocoN} collapse to heavily over-determined linear systems for the column sums $\eins\T K_0,\eins\T K_N$, implying additional restrictions for their right-hand sides, in particular for the matrix $\QQ_{q,r}$.
Looking at an equation
\[ \LL_{q,r}(V_q\T K_nV_r)=\Theta:=(\vartheta_{ij})\in \R^{q\times r}\]
of the form of \eqref{coco0} or \eqref{cocoN}, $n\in\{0,N\}$,
and considering element $(l,k)$ of it, with \eqref{olq} we get
\begin{align}\notag
 & e_l\T (\tilde E_q\T V_q\T K_nV_r+V_q\T K_nV_r\tilde E_r)e_k\\[1mm]\notag
  &=(l-1)(\cc^{l-2})\T K_n\cc^{k-1}+(k-1)(\cc^{l-1})\T K_n\cc^{k-2}\\\label{LKGam}
  &=(\eins\T K_n)(l+k-2)\cc^{l+k-3}
  \stackrel!=\vartheta_{lk},\quad 1\le l\le q,\ 1\le k\le r.
\end{align}
These are $q\cdot r$ conditions for the only $s$ degrees of freedom in $\eins\T K_n$, where the index $j:=l+k-2$ lies in the range $\{0,\ldots,q+r-2\}$.
Obviously, this system is only solvable if $\vartheta_{11}=0$, see \eqref{LL11}, and if $\vartheta_{lk}$ does only depend on $l+k$, which means that $\Theta$ has Hankel form.
For $r+s-2>s$, even more restriction follow.
These additional requirements will be collected farther down for different choices of $q$ and $r$.
\begin{lemma}\label{LSpSum}
Let \eqref{olq} hold and assume that the equations \eqref{coco0} and \eqref{cocoN} have solutions $K_0,K_N$.\\
a) If $q+r\ge s+2$, then the column sums $\eins\T K_0,\eins\T K_N$ are uniquely determined by the nodes and the matrix $B$ of the standard scheme alone.
These sums are solutions of the non-singular linear systems
\begin{align}\label{einsK0a}
 (\eins\T K_0)V_sD_s=&\,\big(((\cc+\eins)^{l-1})\T B\cc^{j-l+1}\big)_{j=1}^s,
 \\\label{einsKNa}
 (\eins\T K_N)V_sD_s=&\,\eins\T-\big((\cc^{l-1})\T B(\cc-\eins)^{j-l+1}\big)_{j=1}^s,
\end{align}
where $D_s=\diag(1,2,\ldots,s)$.
In particular, solvability requires that the expressions on the right-hand sides of these equations do not depend on the index $l$, $1\le l\le\min\{q,j+1\}$.\\
b) For $r=s$, the column sums may be given by
\begin{align}\label{SpSumK0}
 \eins\T K_0=&\,(\cc+\eins)\T BV_sD_s^{-1}V_s^{-1},\\
\label{SpSumKN}
 \eins\T K_N=&\,(\eins\T-\cc\T BV_s\PP_s^{-1})D_s^{-1}V_s^{-1}.
\end{align}
\end{lemma}
{\it Proof:}
a) Omitting the trivial first equation $0=\vartheta_{11}$ in \eqref{LKGam}, for $q+r-3\ge s-1$ the next $s$ cases with $1\le l+k-2\le s$ may be combined to the system
\begin{align*}
 \eins\T K_n\big(1,2\cc,\ldots,s\cc^{s-1})=\eins\T K_n V_sD_s
 \stackrel!=\big(\vartheta_{l,j-l+1}\big)_{j=1}^s,
\end{align*}
$n\in\{0,N\}$, where the index $l$ selects one of several possible choices.
Since the matrix $V_sD_s$ is non-singular, solutions are unique, if they exist.
Equations \eqref{einsK0a} and \eqref{einsKNa} are obtained with $\Theta=\PP_q\T V_q\T BV_r-e_1e_1\T$ and $\Theta=\eins_q\eins_r\T -\QQ_{q,r}$, respectively, where \eqref{eBe} implies $\vartheta_{11}=0$.\\
b) For $r\!=\!s$, the choice $l\!=\!2$ gives a complete row vector of length $s$ in both equations.
\qed
\par\noindent
\begin{remark}
There are practical consequences for the design of Peer triplets.
In the beginning, one may have hoped that positivity $\eins\T K_0>0\T,\eins\T K_N>0\T$ could be obtained in the final design of the end methods with the aid of the many remaining degrees of freedom in the matrices $K_0,K_N$.
However, the Lemma prohibits that.
Instead, for $q+r\ge s+2$ all row sums are determined by the standard method and the positivity restrictions may now be included in search procedures for the standard method alone having fewer degrees of freedom than the whole triplet.
\end{remark}
For methods with forward order $r=s$ and adjoint order $q=s-1$, it holds $q+r-2=2s-3>s$ for $s\ge4$.
In this case there are still more requirements for the solvability of \eqref{LKGam} which are discussed for the special case $(r,q)=(4,3)$ in Section~\ref{SExt43} below.
\section{Four-stage triplets for the order pair $(r,q)=(3,3)$}\label{SP4o33}
The restriction to positive column sums of all matrices $K_n$ reduces the degrees of freedom in the Peer triplets.
In order to regain flexibility in their design, we start with the order pair $(r,q)=(3,3)$ with smaller order $r=3$ than in \cite{LangSchmitt2022b}.
Since Lemma~\ref{LSpSum} represents a bottleneck in the design process, we have to consider the boundary methods first.
Now, the large number of parameters for the boundary methods and huge algebraic expressions bring the formal elimination of the 3 original order conditions for the end method to its limits.
One may circumvent these difficulties by solving the combined condition \eqref{cocoN} formally (with free parameters) for $K_N$ first and then the other conditions in a step-by-step fashion as follows:
\begin{align*}
 \mbox{solve }\left\{\begin{array}{rll}
 \LL_{3,3}(V_3\T K_N V_3)=&\eins_3\eins_3\T -\QQ_{3,3}&\mbox{ for }K_N\\[1mm]
 A_N\T V_3=&w\eins_3\T -K_N\T V_3\tilde E_3&\mbox{ for }A_N,\\[1mm]
 B_N\T V_3=&B\T V_3&\mbox{ for }B_N.
 \end{array}\right.
\end{align*}
\par
The column sums $\eins\T K_0,\,\eins\T K_N$ of the boundary methods are still uniquely determined by the standard method through \eqref{LKGam} since $q+r-2=4=s$.
We may write such a system with $n\in\{0,N\}$ in the form
\begin{align}\label{zsys33}
\R^{1\times 4}\ni\eins\T K_n(\eins,2\cc,3\cc^2,4\cc^4)
 =\begin{pmatrix}
 \vartheta_{12}&\vartheta_{13}&\ast&\ast\\
 \vartheta_{21}&\vartheta_{22}&\vartheta_{23}&\ast\\
 \ast&\vartheta_{31}&\vartheta_{32}&\vartheta_{33}
\end{pmatrix}.
\end{align}
This over-determined system with $s=4$ columns is written with some abuse of notation, 
where potentially conflicting entries on the right-hand side are stacked in the same column.
Asterisks indicate slack variables where no conditions have to be satisfied.
Solutions will exist only if entries within each column have the same value.
This corresponds to the Hankel form of $\Theta\in\R^{3\times 3}$, which means Hankel form of $\QQ_{3,3}$ for $n=N$.
Then, Hankel form of $\PP_3\T V_3\T BV_3=\PP_3\T \QQ_{3,3}\PP_3$ for $n=0$ follows from Lemma~\ref{LHankel}.
Since the matrix $(\eins,2\cc,3\cc^2,4\cc^4)=V_4D_4$ is non-singular, the row sums $\eins\T K_n$ are uniquely determined by \eqref{zsys33}.
A possible choice of elements $\vartheta_{ij}$ is used in the following equation which enforces the positivity conditions \eqref{pos-cond} with the representations
\begin{align}\label{Kpos33}
 \eins\T K_n=(\vartheta_{21},\vartheta_{22},\vartheta_{23},\vartheta_{33})D_4^{-1}V_4^{-1}
 \stackrel!{\ge}\kappa_\ast\eins\T,\ n\in\{0,N\},
\end{align}
with $\Theta=\PP_3\T\QQ_{3,3}\PP_4$ for $n=0$ and $\Theta=\eins_3\eins_3\T-\QQ_{3,3}$ for $n=N$.
\par
Considering the good performance of a triplet based on the backward difference formula BDF4 in \cite{LangSchmitt2022b} and \cite{LangSchmitt2023a}, it is of interest to look for a version with boundary methods satisfying the additional one-leg-conditions \eqref{add_cond} or \eqref{olq} for $q=3$.
According to \eqref{HankelQ}, Hankel form of $\QQ_{3,r}$ is required now.
However, the matrix $\QQ_{4,4}$ for the BDF standard method has the Hankel property in its first 4 anti-diagonals only with $e_1\T \QQ_{4,4}=(1,\frac18,\frac1{96},0)$.
Hence, only $\QQ_{3,3}$ has Hankel form since it contains only one element from the fifth anti-diagonal.
Now, for a method with $q=r=3$, the column sums of the boundary methods are explicitly determined by the standard method through \eqref{zsys33}.
Here, for BDF the end method satisfies \eqref{Kpos33} with $\kappa_*=55/576$, but not the starting method since $\eins\T K_0e_2=-21/64<0$.
Hence, no positive triplet based on BDF4 exists satisfying the one-leg-condition with $q\ge3$.
\subsection{Method \texttt{AP4o33pa}}
So far, we have only discussed the normal order conditions for the Peer methods and their adjoints and the resulting restrictions.
For methods to be efficient further requirements have to be considered.
First, the conditions \eqref{OBvStd} and \eqref{OBaStd} for the pair $(r,q)$ relate to the (local) orders of consistency.
In order to also establish convergence of (global) order $O(h^r)$ and $O(h^q)$ in \cite{LangSchmitt2022b}, the following two conditions for super-convergence of the forward and adjoint scheme have been added,
\begin{align}\label{SupKnv}
 \eins\T\big(A\cc^r-B(\cc-\eins)^r-rK\cc^{r-1}\big)=&\,0,
 \\\label{SupKna}
 \eins\T\big(A\T \cc^q-B\T(\cc+\eins)^q+qKc^{q-1}\big)=&\,0,
\end{align}
which cancel the leading term in the global error.
In practice, the super-convergence effect may be observed for a sufficiently fast damping of secondary modes of the stability matrix only and requires that
\begin{align}\label{lbd2}
 |\lambda_2(A^{-1}B)|\le\gamma<1,\quad \gamma\cong 0.8,
\end{align}
where $\lambda_2$ denotes the absolutely second largest eigenvalue of the stability matrix $\bar B:=A^{-1}B$ of the standard scheme.
We note, that the stability matrix $\tilde B\T:=(BA^{-1})\T$ of the adjoint time steps has the same eigenvalues as $\bar B$.
In several numerical tests in \cite{LangSchmitt2022b}, the given value $\gamma=0.8$ was sufficiently small to produce super-convergence reliably and it does so in our tests at the end.
Superconvergence  \eqref{SupKnv}, \eqref{SupKna} only cancels the leading error terms of the methods.
In order to cover other essential error contributions also the norms
\begin{align*}
 err_r:=&\frac1{r!}\|\cc^r-A^{-1}B(\cc-\eins)^r-rA^{-1}K\cc^{r-1}\|_\infty,\\
 err^\dagger_q:=&\frac1{q!}\|\cc^q-A\mT B\T(\cc+\eins)^q+q A\mT Kc^{q-1}\|_\infty,
\end{align*}
are monitored  as the essential error constants, see \cite{LangSchmitt2022b}.
Furthermore, the norm $\|A^{-1}B\|_\infty$ of the stability matrix is of interest since it may be a measure for the propagation of rounding errors.
\par
Application of the Peer methods to stiff problems requires good stiff stability properties.
$A(\alpha)$-stability is defined here by the requirement that the spectral radius $\varrho\big((A-zK)^{-1}B\big)<1$ of the stability matrix of the standard scheme is below one for $z$ in the open sector of the complex plain with aperture $2\alpha$ centered at the negative real axis.
The adjoint stability matrix $(A-zK)\T B\T$ and $\big(A-zK)^{-1}B$ possess the same eigenvalues.
Details on the computation of $\alpha$ can be found in \cite[\S5.2]{LangSchmitt2022}. The angles for the different methods are contained in Table~\ref{TPT}.
\par
Very mild eigenvalue restrictions for the boundary methods are also taken from \cite{LangSchmitt2022b}.
In order to guarantee the solvability of the stage systems for the first and last steps we require
\begin{align}\label{LsbRS}
 \mu_0:=\min_j\Re\lambda_j(K_0^{-1}A_0)>0,\quad
 \mu_N:=\min_j\Re\lambda_j(K_N^{-1}A_N)>0.
\end{align}
\par
A new requirement is the non-negativity condition \eqref{pos-cond} imposed by the use of a gradient-based method to update the control vector in \eqref{nopt_iter}.
Lemma~\ref{LSpSum} has shown that the column sums of the boundary methods are already fully determined by the standard method and their positivity $\eins\T K_0>0\T,\eins\T K_N>0\T$ can be included in the search for it.
In practice, the performance of the gradient method \eqref{nopt_iter} may suffer badly if the column sums have differing magnitudes.
Hence, the search was narrowed to methods with moderate positive values of the {\bf c}olumn {\bf s}um {\bf q}uotient
\begin{align}\label{csq}
 csq:=\max\Big\{\frac{\max_i|\eins\T K_0e_i|}{\min_i\eins\T K_0e_i},\frac{\max_i|\eins\T K_Ne_i|}{\min_i\eins\T K_Ne_i}\Big\}\stackrel!>0,
  \eins\T K\ge 0\T,
\end{align}
and \eqref{pos-cond} where $\eins\T K_0,\eins\T K_N$ are determined by the standard method, see \eqref{SpSumK0}, \eqref{SpSumKN}.
Although there are rather tight restrictions on the column sums of $K_0,K_N$, there still exists a null-space in the conditions for these matrices and it was necessary to restrict the norms $\|K_0\|,\|K_N\|$ in addition to \eqref{LsbRS} in the final search for the boundary methods.
\par
For easier reference, the full set of conditions is collected in Table~\ref{TOBed}.
We remind that the conditions in line (d) there  ensure the existence of the boundary methods and allow for the detached construction of the standard method alone.
\par
Although searches for positive standard methods with $r=3=s-1$ came very close to A-stability, no truly A-stable methods could be found.
This might be due to the restriction \eqref{lbd2} on the sub-dominant eigenvalue and we suspect that a (formally) A-stable method might have a multiple eigenvalue 1 if it exists.
In fact, with a rather unsafe restriction $|\lambda_2(A^{-1}B)|\le0.9$, a method was found with stability angle $89.976^o$ extremely close to A-stability.
Slightly relaxing the requirement on the angle to $\alpha=89.90^o$, the following
{\bf a}lmost A-stable method \texttt{AP4o33pa} was constructed.
Its node vector is given by
\begin{align*}
  \cc\T=\left(\frac{46}{5253},\frac{29}{51},\frac{1723}{2193},\frac{17131}{12189}\right)
  \doteq(0.00876,0.5686,0.7857,1.4054)
\end{align*}
with monotonic nodes.
It is super-convergent with \eqref{SupKnv},\eqref{SupKna} for $r=q=3$ with a good damping factor $|\lambda_2(A^{-1}B)|<0.66$.
The error constants are almost equal with $err_3=0.298/3!\approx 0.050$ and $err^\dagger_3=0.279/3!\approx 0.046$.
Further data of this method are collected in Table~\ref{TPT}.
All boundary steps are zero-stable but only for the end method a block structure could be obtained with block sizes \texttt{blksz}=$3+1$.
These data are presented in Table~\ref{TPTB}. The complete set of coefficients is given in Appendix A.1.
\begin{table}
 \centerline{\begin{tabular}{|c|l|c|c|c|c|c|c|c|}\hline
  triplet &$(r,q)$&nodes&$\alpha$&$\|A^{-1}B\|_\infty$&$|\lambda_2|$&$err_r$&$err_q^\dagger$&$csq$\\\hline
  AP4o33pa&$(3,3)$&$[0,1.41]$&$89.90^o$&$8.2$&$0.66$&0.050&0.046&$33.4$\\
  AP4o33pfs&$(3,3)$&$[0,1]$&$77.53^o$&$16.0$&$0.46$&0.031&0.030&$1.72$\\
  AP4o43p&$(4,3)$&$[0.1,0.9]$&$59.78^o$&8.5&0.58&0.0038&0.024&11.0\\
  \hline
 \end{tabular}}
 \caption{Properties of the 4-stage standard methods of Peer triplets.}\label{TPT}
\end{table}
\begin{table}
 \centerline{\begin{tabular}{|l|c|c|c|c|c|c|c|}\hline
 &\multicolumn{3}{c|}{Starting method}&\multicolumn{4}{c|}{End method}\\\cline{2-8}
 triplet &blksz&$\mu_0$&$\varrho(BA_0^{-1})$& blksz&$\mu_N$&$\varrho(A_N^{-1}B_N)$&$\varrho(B_NA^{-1})$\\\hline
 AP4o33pa&4&2.03&1& 1+3&2.21&1&1\\
 AP4o33pfs&1+3&4.92&1& 1+3&1.61&1&1\\
 AP4o43p&4&4.13 &1 &4 &4.36 & 1& 1.09\\
\hline
 \end{tabular}}
 \caption{Properties of the boundary methods of Peer triplets.}\label{TPTB}
\end{table}
\subsection{FSAL method \texttt{AP4o33pfs}}\label{SecFSAL}
The number of stage equations to be solved numerically may be reduced with the aid of the FSAL property ({\em first stage as last}) frequently used in the design of one-step methods, where the last stage of the previous time step equals the first stage of the new step.
For Peer methods, this property has been discussed in \cite{SchmittWeinerBeck2013}.
In our formulation \eqref{KKT_state_peer} it means that
\begin{align}\label{PFSAL}
 c_1=0,\,c_s=1,\quad e_1\T K_n=0\T,\ e_1\T A_n=a_{11}^{[n]}e_1\T,
  \ e_1\T B_n=a_{11}^{[n]}e_s\T,
\end{align}
implying $Y_{n,1}=Y_{n-1,s}\cong y(t_n)$, $n\ge1$.
A convenient benefit of Peer methods is that, due to their high stage order, the interpolation of all $s$ stages provides an accurate polynomial approximation of the solution being also continuous if the FSAL property holds.
\par
Since the order conditions for methods of type \texttt{AP4o33*} leave a 10-parameter family of standard methods, the additional restrictions \eqref{PFSAL} can easily be satisfied for $(A,B,K)$.
However, some properties of the boundary methods imply further restrictions on the standard method through \eqref{zsys33}.
Although the matrices $K_0,K_N$ in the boundary steps are not restricted to diagonal form, the one-leg conditions \eqref{olq} with $q=s-1$ request that they are rank-1-changes of diagonal matrices only.
Then, the condition $e_1\T K_n=0\T$ leaves off-diagonal elements in their first columns only.
However, for matrices of such a form, the constraint \eqref{KKT_ctr_peer_all} reads
\begin{align*}
 \nabla_u f(Y_{nj},U_{nj})\T P_{nj} \kappa_{jj}^{[n]}=&\,0,\ j=2,\ldots,s,\\
 \nabla_u f(Y_{n1},U_{n1})\T\sum_{i=2}^sP_{ni} \kappa_{i1}^{[n]} =&\,0.
\end{align*}
These are $s$ conditions on $P_{n2},\ldots,P_{ns}$, which may not always be satisfiable since $\nabla_u f\T$ is evaluated at different places.
Hence, the condition $e_1\T K_n=0\T$ requires that $K_n$ is diagonal with zero as the first diagonal element.
However, the property $K_ne_1=0$ also leads to $k_{n1}=\eins\T K_ne_1=0$ and introduces via \eqref{zsys33} one additional restriction on the matrix $\QQ_{3,3}$ from the standard method both for $K_0$ and $K_N$.
In condition \eqref{Kpos33}, the first component, being zero, can be deleted since also the control $U_{n1}$ is no longer present in \eqref{KKT_ctr_peer_all}.
\par
Unfortunately, only in the starting step a diagonal matrix $K_0\succeq0$ with $\kappa_{11}^{[0]}=0$ is possible leading to an exact start $Y_{01}=y_0$.
However, no non-negative triplet seems to exist with a final FSAL step.
Hence, $K_N$ is chosen lower triangular, having a dense first column and the first row $e_1\T K_N=\frac13e_1\T$, leading to a small jump $Y_{n1}-Y_{n-1,s}=O(h^3)$ at $t_N$ only.
\par
Computer searches found the method \texttt{AP4o33pfs} with stability angle $\alpha=77.53^o$, having node vector $c\T=(0,\frac9{86},\frac{321}{602},1)$, error constants $\eta_3=0.187/3!\approx 0.031$, $\eta_3^\dagger=0.180/3!\approx 0.030$ and a small damping factor $\gamma=0.46$.
More data are given in Table~\ref{TPT} and Table~\ref{TPTB}.
Of course, the computation of the quotient $csq$ in \eqref{csq} and the real part $\mu_0$ in \eqref{LsbRS} was restricted to the nontrivial lower $3\times 3$ block of $K_0$.
The coefficients of \texttt{AP4o33pfs} are given in Appendix A.2.
\begin{table}[t!]
\begin{center}
\begin{tabular}{|cl|c|c|}\hline
 &Steps&forward: $r\le s=4$&adjoint: $q=s-1=3$\\\hline
 (a)&Start, $n=0$&\eqref{OBvStrt}&\eqref{OBaStd}, $n=0$, \eqref{add_cond}, \eqref{pos-cond}\\\hline
 (b)&Standard, $1\le n<N$&\eqref{OBvStd} &\eqref{OBaStd}\\
 (c)&Superconvergence  &\eqref{SupKnv} &\eqref{SupKna}\\\cline{2-4}
 (d)&Compatibility &\multicolumn{2}{c|}{\eqref{eBe}, \eqref{HankelQ}}\\
 &$\ldots$ for $r=s=4$:&\multicolumn{2}{c|}{\eqref{Korth43}, \eqref{Bnot43a}}\\\hline
 (e)&Last step&\eqref{OBvStd}, $n=N$&\eqref{OBaStd}, $n=N-1$\\
 (f)&End point& \eqref{OBvEnd}&\eqref{OBaEnd}, \eqref{add_cond}, \eqref{pos-cond}\\\hline
\end{tabular}
\caption{Combined order conditions for the Peer triplets \texttt{AP4o43p} and \texttt{AP4o33pa}}\label{TOBed}
\end{center}
\end{table}
\section{Four-stage triplets for the order pair $(r,q)=(s,q)=(4,3)$}\label{SP4o43}
Obeying the additional order and super-convergence conditions for $r=s=4$ means that the Peer method will have global order 4 for pure initial value problems without control.
Unfortunately, in Section~\ref{SGlobErr} we will not be able to prove the same improvement for the state solutions of the control problem if the adjoint equations remain at order $q=3$.
However, the numerical experience in our previous papers \cite{LangSchmitt2022b,LangSchmitt2023a} indicates that the coupling between both errors is often rather weak and the higher order may show up in numerical tests.
The A-stable methods from \cite{LangSchmitt2022b} already used this combination of $(r,q)$ for problems where $u$ may be explicitly eliminated, leading to a boundary value problem for $(y,p)$.
However, these methods are not suited within a gradient-based optimization method discussed in Section~\ref{SecGradM} since $K$ possesses negative diagonal elements.
Still, two methods \verb+AP4o43bdf+ and \verb+AP4o43dif+ satisfy condition \eqref{pos-cond} and may be used in a gradient-based optimization.
In our numerical tests, \texttt{AP4o43bdf} will be compared to the new methods derived now and which satisfy one set of one-leg conditions more, see \eqref{add_cond}.
These stronger requirements on all methods of the triplet lead to a severe bottleneck in the augmented order conditions: any appropriate standard method $(A,B,K)$, $K=\diag(\kappa_{ii}),$ has a blind third stage with $\kappa_{33}=0$.
This observation is a consequence of equation \eqref{cocoN} and Lemma~\ref{LKol}.
The full set of conditions will be collected at the end of this section.
\subsection{Consequences of the Hankel form of $\QQ_{3,4}$}
In this subsection, it is shown that the restriction $\kappa_{33}=0$ is a consequence of only the forward order conditions \eqref{OBvStd} with $r=s=4$ and the $q\times s=3\times 4$-Hankel form \eqref{HankelQ} of the matrix
\begin{align}\label{Q34}
 \QQ_{q,s}=V_q\T BV_s\PP_s^{-1}=V_q\T (AV_s-KV_s\tilde E_s).
\end{align}
With the shift matrices $S_q=\big(\delta_{i,j-1}\big)\in\R^{q\times q}$ and the projection $\check I_q:=I_q-e_qe_q\T$, Hankel form of the matrix \eqref{Q34} is equivalent with
\begin{align}\label{Hank}
 0=\check I_q(S_q\QQ_{q,s}-\QQ_{q,s}S_s\T)\check I_s.
\end{align}
Also, column shifts in the Vandermonde matrix have a simple consequence, $V_qS_q\T=CV_q\check I_q$.
\begin{theorem}\label{TAKH}
a) Hankel form of the matrix \eqref{Q34} is equivalent with the condition
\begin{align}\label{AKHm}
 V_{q-1}\T(AC-CA-K)V_{s-1}=0.
\end{align}
b) For $q=3$, $r=s=4$, equations \eqref{Q34}, \eqref{Hank} imply
$$\kappa_{33}=e_3\T Ke_3=0.$$
\end{theorem}
{\it Proof:}
Since $\tilde E_r=D_rS_r=S_r(D_r-I_r)$ with $D_r=\diag_{i=1}^r(i)$, we have $\tilde E_r S_r\T=D_r\check I_r$ and $CV_r\tilde E_r=V_r(D_r-I_r)$.
Now, the Hankel condition \eqref{Hank} for the matrix \eqref{Q34} reads
\begin{align*}
 0=&\check I_q\big(S_q V_q\T (AV_s-KV_s\tilde E_s)-V_q\T (AV_s-KV_s\tilde E_s)S_s\T\big)\check I_s\\
 =&\check I_q( V_q\T C(AV_s-KV_s\tilde E_s)-V_q\T ACV_s+V_q\T KV_s D)\check I_s\\
 =&\check I_q V_q\T( CAV_s- ACV_s-K(CV_s\tilde E_s-V_s D))\check I_s\\
 =&\check I_q V_q\T(CA- AC+K)V_s\check I_s.
\end{align*}
Since the matrices $\check I_q,\check I_s$ simply eliminate the last row or column, this equation is equivalent with the assertion \eqref{AKHm}.\\
b) For $q=3,\,s=4$, condition \eqref{AKHm} consists of 6 equations.
The commutator $[A,C]=AC-CA$ is strictly lower triangular since the diagonal of $A$ cancels out.
Ignoring the diagonal, the map $A\mapsto V_{q-1}\T(AC-CA-K)V_{s-1}$ still has a rank deficiency if it is considered as a function of the 6 subdiagonal elements of $A$ only.
In fact, there exists a nontrivial kernel of its adjoint having rank-1 structure.
Consider
\[ V_3\begin{pmatrix}c_1c_2\\ -c_1-c_2\\1 \end{pmatrix}
 =\begin{pmatrix}0\\0\\(c_3-c_1)(c_3-c_2)\\(c_4-c_1)(c_4-c_2)\end{pmatrix}=:x_R,\quad
 V_2\begin{pmatrix}-c_4\\1\end{pmatrix}
 =\begin{pmatrix}\ast\\\ast\\c_4-c_3\\0 \end{pmatrix}=:x_L.
\]
Since $[A,C]$ is strictly lower triangular, the vector $[A,C]x_R$ has 3 leading zeros and its inner product with $x_L$ vanishes.
Hence,
\begin{align*}
 (-c_4,1)V_{2}\T(AC-CA-K)V_{3}
 \begin{pmatrix}c_1c_2\\ -c_1-c_2\\1 \end{pmatrix}
 =&\, -x_L\T Kx_R \\
 =&\, -(c_4-c_3)(c_3-c_2)(c_3-c_1)\kappa_{33}
\end{align*}
and \eqref{AKHm} implies $\kappa_{33}=0$ for non-confluent nodes. \qed
\par
Since $K$ is diagonal, the matrix $V_{q-1}\T KV_{s-1}$ has Hankel form again and  \eqref{AKHm} is an over-determined system for the column sums $\eins\T K=(\kappa_{11},\ldots,\kappa_{ss})$.
Solutions may only exist if elements of $V_{q-1}\T(AC-CA)V_{s-1}$ within each anti-diagonal have the same value.
For $q=3,\,s=4$, this leads to the following restrictions on $A$ alone:
\begin{align}\label{KompA}
 (\cc^2)\T A\cc^j-2\cc\T A\cc^{j+1}+\eins\T A\cc^{j+2}=0,\ j=0,1.
\end{align}
Since a $2\times 3$ matrix possesses 4 antidiagonals, under assumption \eqref{KompA} the system \eqref{AKHm} reduces to
\begin{align*}
 \eins\T KV_4=&\,\eins\T(AC-CA)(\eins,\cc,\cc^2,0)+\cc\T(AC-CA)\cc^2\\
 \gdw \eins\T K=&\,\eins\T(AC-CA)+\beta e_4\T V_4^{-1},\\
 \mbox{where } \beta=&\,\cc\T(AC-CA)\cc^2-\eins\T(AC-CA)\cc^3
  =\eins\T[C,[A,C]]\cc^2.
\end{align*}
\begin{remark}
Since $\kappa_{33}=0$, the third stage of the standard method uses no additional function evaluation of $f(Y_{n3},U_{n3})$ and it seems that it does not provide any additional information.
In fact, this stage can be eliminated but the resulting method will be a 3-stage 3-step Peer method.
\end{remark}
\begin{remark}\label{Rblst}
The blind third stage has consequences both for the analysis and the implementation of the standard Peer method.
In the equations \eqref{KKT_state_peer_init}--\eqref{KKT_ctr_peer} of the Peer steps, it is seen that there is no coupling between the controls $U_n$ of different time steps.
Now, the contribution to the Lagrange function from the third stage of the standard method in time step $n$ with multiplier $P_{n3}$ is given by
\begin{align*}
 \ldots+P_{n3}\T\left(\sum_{j=1}^3a_{3j}Y_{nj}-\sum_{j=1}^4b_{nj}Y_{n-1,j}\right)+\ldots
\end{align*}
missing the unknown $U_{n3}$.
Since $U_{n3}$ does not appear anywhere else, it is non-existent and the unknown $U_{n3}$ should be discarded as well as the corresponding stage equation $0\cdot\nabla_u f(Y_{n3},U_{n3})\T P_{n3}=0$ from \eqref{KKT_ctr_peer_allk}.
This measure will also be used in the analysis of Section~\ref{SGlobErr}.
\end{remark}
\subsection{Further requirements for the existence of Peer triplets}\label{SExt43}
For $q+r-2>s$, the restriction of $\QQ_{q,r}$ to Hankel form is not the whole picture yet.
If $q+r-2=s+1$, which case occurs for $r=s=4$ and $q=s-1=3$, the system \eqref{LKGam} does not possess full rank having more than $s$ columns.
The vector $(\psi\T,1)\T$ with $\psi:=-V_s^{-1}\cc^s$ spans the kernel of the extended Vandermonde matrix $V_{s+1}=(\eins,\cc,\ldots,\cc^s)$ since it contains the coefficients of the node polynomial $\hat\psi(t)=(t-c_1)\cdots(t-c_s)=t^s+\sum_{j=1}^s\psi_jt^{j-1}$.
Rewriting this property in the form of \eqref{LKGam},
\begin{align*}
 \eins\T K_n\big(1,2\cc,\ldots,(s+1)\cc^s\big)\cdot D_{s+1}^{-1}\begin{pmatrix}
  \psi\\1
 \end{pmatrix}=0,
\end{align*}
it follows that solutions only exist if also
\begin{align}\label{slvb43}
 \big(\vartheta_{l,j-l+1}\big)_{j=1}^{s+1}\cdot D_{s+1}^{-1}{\psi\choose 1}=0.
\end{align}
Since $r=s$, a convenient choice for the indices $l$ here is
\begin{align}\label{Auswahl}
\big(\vartheta_{21},\ldots,\vartheta_{2s},\vartheta_{3s}\big)
=\left\{\begin{array}{ll}
 \big(\eins\T-\cc\T BV_s\PP_s^{-1},1-(\cc^2)\T B(\cc-\eins)^{s-1}\big)&\mbox{for }K_N,\\[1mm]
 \big((\cc+\eins)\T BV_s,((\cc+\eins)^2)\T B\cc^{s-1}\big)&\mbox{for }K_0.
\end{array}\right.
\end{align}
We summarize all conditions in the following lemma.
\begin{lemma}
Let $r=4$ and $q=3$.
Then, necessary conditions for the existence of four-stage Peer triplets satisfying \eqref{OBvStrt}--\eqref{OBaEnd} and \eqref{olq} are Hankel form of $\QQ_{3,4}=\PP_3\T BV_4\PP_4^{-1}$ as well as
\begin{align}\label{Gnot43a}
 \cc\T BV_4P_4^{-1}D_4^{-1}\psi+\frac15(\cc^2)\T B(\cc-\eins)^3=&\,\int_0^1\hat\psi(t)dt
 \\\label{Gnot43b}
 (\cc+\eins)\T BV_4 D_4^{-1}\psi+\frac15((\cc+\eins)^2)\T B\cc^3=&\,0.
\end{align}
\end{lemma}
{\it Proof:}
As a first step, we consider the constant term $\eins_{s+1}\T$ in the equation for $K_N$ in \eqref{Auswahl}.
In \eqref{slvb43} it gives rise to the contribution
\begin{align*}
 \eins_s\T D_s^{-1}\psi+\frac1{s+1}=\sum_{j=0}\frac1j\psi_j+\frac1{s+1}=\int_0^1\hat\psi(t)dt.
\end{align*}
Now, the conditions \eqref{Gnot43a}, \eqref{Gnot43b} correspond to the vanishing of the inner products \eqref{slvb43} with the two vectors from \eqref{Auswahl}.
\qed
\par\noindent
\begin{remark}
In practice it was found that the two conditions \eqref{Gnot43a}, \eqref{Gnot43b} seem to be equivalent with the two simpler equations
\begin{align}\label{Korth43}
 \int_0^1\hat\psi(t)dt=&\,0,\\\label{Bnot43a}
 \cc\T BV_4P_4^{-1}D_4^{-1}\psi+\frac15(\cc^2)\T B(\cc-\eins)^3=&\,0,
\end{align}
in conjunction with the many other order conditions.
\end{remark}
In general, Peer methods are invariant under a common shift of the nodes, which means in practice that this shift may be fixed after the construction of some method, e.g. by choosing $c_s=1$.
However, orthogonality \eqref{Korth43} strongly depends on the absolute positions of the nodes.
Still, for $\int_0^1\hat\psi(t)dt$ this dependence is only linear for the node differences and \eqref{Korth43} may be easily solved for one of those, e.g. for $d_4=c_4-c_2$.
\subsection{Method \texttt{AP4o43p}}
Without the non-negativity condition several different regions in the parameter space of Peer triplets did exist for $(r,q)=(4,3)$ in \cite{LangSchmitt2022b}.
Some of the standard methods found there have non-monotonic nodes and negative diagonal elements in $K$.
Now, non-negativity seems to leave only one such region with ordered nodes $c_i$ and we present only one such method with a nearly maximal stability angle.
Method \texttt{AP4o43p}
has a stability angle of $\alpha=59.78^o$ with node vector
\begin{align}\notag
 \cc\T=\big(\frac{4657}{46172},\frac{43}{97},\frac{3991}{6596},\frac{21111803999}{23798723875}\big)
 \doteq\big(0.1009,0.4432,0.6050,0.8871\big).
\end{align}
The node $c_4$ has a rather long representation since it was used to solve condition \eqref{Korth43}.
The damping factor $\gamma=0.58$ from \eqref{lbd2} is well below one, the error constants are
$err_4=0.092/4!\approx 0.0038$ and $err^\dagger_3=0.144/3!\approx 0.024$ and the quotient \eqref{csq} is $csq=11.0$.
Further data are collected in Table~\ref{TPT}.
In order to obtain acceptable properties for the stability and definiteness \eqref{LsbRS} of the boundary methods, the matrices $A_0,A_N$ have full block size 4, denoted by \texttt{blksz}=4 in Table~\ref{TPTB}. The coefficients of \texttt{AP4o43p} are given in Appendix A.3.
\section{The global error}\label{SGlobErr}
Convergence of the Peer triplets for $h\to0$ will be discussed for the unconstrained case $N_U=\{0\}$ only.
Here, the additional constraint \eqref{KKT_ctr_peer}, \eqref{KKT_ctr_peer_all} complicates the situation compared to \cite{LangSchmitt2022,LangSchmitt2022b} and we will extend these proofs here.
Node vectors of the exact solution are denoted by bold face, e.g., ${\bf y}_n=\big(y^\star(t_{ni})\big)_{i=1}^s$, ${\bf y}:=\big({\bf y}_n\big)_{n=0}^N$, and the global errors by checks, e.g., $\check Y_{n}=Y_n-{\bf y}_n$, $\check Y:=\big({\check Y}_n\big)_{n=0}^N$.
For ease of writing, we also introduce the combined vectors ${\bf z}:=({\bf y}\T,{\bf p}\T,{\bf u}\T)\T$ and $\check Z=(\check Y\T,\check P\T,\check U\T)\T$, and use abbreviations like $F(Z_n):=F(Y_n,U_n)$ and $\nabla_z f=(\nabla_y f,0,\nabla_u f)$.
\par
In the error discussion, a notational difficulty arises since the right-hand sides of the adjoint equations already contain a first derivative $(\nabla_y f)\T p$.
In order to avoid ambiguities with second derivatives, we introduce an additional notation $\langle.,.\rangle$ for the standard inner product in $\R^m$ which is exclusively dedicated to the product of the Lagrange multiplier $p$ and the components of the function $f=(f_i)_{i=1}^m$ or its derivatives as in
\begin{align*}
 \langle p,f\rangle=\sum_{i=1}^m(e_i\T p)f_i,\quad
 \langle p,\nabla_z f\rangle=\sum_{i=1}^m(e_i\T p)\nabla_z f_i.
\end{align*}
The notation is particularly used for second derivatives, where the matrix
$\nabla_{uu}\langle p,f\rangle=\sum_{i=1}^m(e_i\T p)\nabla_{uu}f_i$ is symmetric and a linear combination of Hessian matrices of the components of $f$.
The notation carries over to compound expressions of a whole time step, e.g.
\begin{align*}
 \left\langle K_n\T P_n,\nabla_{UU}F(Y_n,U_n)\right\rangle=\diag_{j=1,\ldots,s}\Big(\left\langle\sum_{i=1}^s P_{ni}\kappa^{[n]}_{ij},\nabla_{uu}f(Y_{nj},U_{nj})\right\rangle\Big).
\end{align*}
\par
With these notations, we may introduce the local errors for the two-step equations
\begin{align}\label{LokFey0}
\tau_0^Y:=&\,{\bf y}_0-\eins\otimes y_0-hA_0^{-1}K_0 F({\bf z}_0),
\\\label{lokFey}
\tau_n^Y:=&\,{\bf y}_n-A_n^{-1}\big(B_n{\bf y}_{n-1}+hK_nF({\bf z}_n)\big),\quad 1\le n\le N,
\\\label{lokFep}
\tau_n^P:=&\,{\bf p}_n-A_n\mT\big(B_{n+1}\T{\bf p}_{n+1}+h\nabla_Y F({\bf z}_n)\T K_n\T{\bf p}_n\big),\ 0\le n<N,
\\\label{lokFepN}
\tau_N^P:=&\,{\bf p}_N-\eins\otimes p(T)-hA_N\mT \nabla_Y F({\bf z}_N)\T K_N\T{\bf p}_N.
\end{align}
The order conditions from Table~\ref{TOBed} have the following consequences for these errors (see \cite{LangSchmitt2022b} and Lemma~\ref{Lknspf}),
\begin{align}\label{KonsFe}
 \tau_n^Y=O(h^r\|(y^\star)^{(r)}\|_\infty),\quad \tau_n^P=O(h^q\|(p^\star)^{(q)}\|_\infty),\ 0\le n\le N.
\end{align}
Furthermore, line (c) of Table~\ref{TOBed} gives for the interior time steps with $1\le n\le N-1$ that
\begin{align}\label{KonsSup}
  \eins\T A\tau_n^Y=O(h^{r+1}\|(y^\star)^{(r+1)}\|_\infty),\quad \eins\T A\T \tau_n^P=O(h^{q+1}\|(p^\star)^{(q+1)}\|_\infty).
\end{align}
The global error in these time steps may be analyzed by multiplying \eqref{KKT_state_peer_init}, \eqref{KKT_state_peer} by inverses $A_n^{-1}$ and \eqref{KKT_adj_peer_init}, \eqref{KKT_adj_peer} by $A_n\mT$ and subtracting the corresponding equations \eqref{lokFey} or \eqref{lokFep}, respectively.
This yields the equations
\begin{align}\label{errYn}
\check Y_n-\bar B_n\check Y_{n-1}=&\,R_n^Y(\check Z)-\tau_n^Y,\ n=0,\ldots,N,\\\notag
 R_n^Y(\check Z):=&\,h\bar K_n\big(F({\bf z}_n+\check Z_n)-F({\bf z}_n)\big),
 \\\label{errPn}
 \check P_n-\tilde B_{n+1}\T\check P_{n+1}=&\,R_n^P(\check Z_n)-\tau_n^P,\ n=0,\ldots,N-1,
 \\\notag
 R_n^P(\check Z):=&\,hA_n\mT\big(\langle K_n\T({\bf p}_n+\check P_n),\nabla_Y 
 F({\bf z}_n+\check Z_n)\rangle-\langle\tilde K_n\T{\bf p}_n,\nabla_Y F({\bf z}_n)\rangle\big),
\end{align}
with $\bar K_n=A_n^{-1}K_n$.
Recall that $\bar B_n=A_n^{-1}B_n$ and $\tilde B_n\T =A_n\mT B_n\T$.
The starting step fits in \eqref{errYn} by setting $\bar B_0=0$.
Restricting the objective function $\CC$ to polynomials of degree two at most the equation for the error in step \eqref{KKT_adj_peer_init} becomes
\begin{align}\label{errPN}
 \check P_N-\big((\eins w\T)\otimes\nabla_{yy}\CC\big)\check Y_N=&\,R_N^P(\check Z)-\tau_N^P,\\\notag
 R_N^P(\check Z)=&\,hA_N^{-1}\big(\langle K_N\T({\bf p}_N+\check P_N),\nabla_Y F({\bf z}_N+\check Z_N)\rangle\\\notag
 &-\langle\tilde K_n\T{\bf p}_N,\nabla_Y F({\bf z}_N)\rangle\big).
\end{align}
\par
The convergence proof will employ a fixed-point argument and its basic principle is to collect all terms of size $O(1)$ (i.e. independent of $h$) implicitly as part of a fixed linear system and terms of size $h$ as Lipschitz functions.
Accordingly, in the equations \eqref{errYn}--\eqref{errPN} the implicit two-step terms were written on the left-hand sides and $O(h)$-terms on the right.
Unfortunately, this separation is not so obvious for the new constraint \eqref{KKT_ctr_peer_all} which may be written as
\begin{align}\label{KKT_ctr_peer_allk}
\left\langle K_n\T P_{n},\nabla_U F(Y_{n},U_{n})\right\rangle\T=0
\end{align}
in our new notation.
We remind that by Remark~\ref{Rblst} the unknowns $U_{n3}$ and the corresponding equations for blind stages ($\kappa_{33}^{[n]}=0$, $1\le n<N$) should be discarded.
\par
For \eqref{KKT_ctr_peer_allk} we resort to Taylor expansion at the exact solution ${\bf z}$ and use the Jacobian at this point in the implicit part of the error equation.
Considering the definition \eqref{DeftauU} of the local error $\tau_n^U$ and that the exact solution $(y^\star(t),p^\star(t),u^\star(t))$ satisfies $\nabla_u f(y^\star,u^\star)\T p^\star=\langle p^\star,\nabla_y f(y^\star,u^\star)\rangle\T\equiv0$, see \eqref{KKT_ctr_ode}, we have
\begin{align*}
 -(\tau_n^U)\T=&\;\langle K_n\T P_n,\nabla_U F(Z_n)\rangle-\langle K_n\T {\bf p}_n ,\nabla_U F({\bf z}_n)\rangle\\[1mm]
 =&\;\langle K_n\T\check P_n,\nabla_U F({\bf z}_n)\rangle
 +\langle K_n\T{\bf p}_n,\nabla_U F(Z_n)-\nabla_U F({\bf z}_n)\rangle\\[1mm]
 &+\langle K_n\T\check P_n,\nabla_U F(Z_n)-\nabla_U F({\bf z}_n)\rangle\\[1mm]
 =&\;\langle K_n\T\check P_n,\nabla_U F({\bf z}_n)\rangle\\[1mm]
 &\;+\langle K_n\T{\bf p}_n,\nabla_{UY}F({\bf z}_n)\check Y_n\rangle
 +\langle K_n\T{\bf p}_n,\nabla_{UU}F({\bf z}_n)\check U_n\rangle\\[1mm]
 &\;+\langle K_n\T{\bf p}_n,\nabla_U F(Z_n)-\nabla_U F({\bf z}_n)\rangle
 - \langle K_n\T{\bf p}_n,\nabla_{UY}F({\bf z}_n)\check Y_n\rangle\\[1mm]
 &\; -\langle K_n\T{\bf p}_n,\nabla_{UU}F({\bf z}_n)\check U_n\rangle
 +\langle K_n\T\check P_n,\nabla_U F(Z_n)-\nabla_U F({\bf z}_n)\rangle.
\end{align*}
Hence, in each time step there is an additional equation
\begin{align}\label{ZGlUn}\notag
 &\Big(\langle K_n\T{\bf p}_n,\nabla_{UY}F({\bf z}_n)\rangle\T,
 \nabla_U F({\bf z}_n)\T K_n,
 \langle K_n\T{\bf p}_n,\nabla_{UU}F({\bf z}_n)\rangle\T\Big)\check Z_n
\\
 &=\;-\tau_n^U+R_n^U(\check Z_n).
\end{align}
The function on its right-hand side is given by
\begin{align}\notag
 R_n^U(\check Z)=
  &\; -\langle K_n\T\check P_n,\nabla_U F({\bf z}_n+\check Z_n)-\nabla_U F({\bf z}_n)\rangle\T
  \\[1mm]\label{ResU}
 &\;-\langle K_n\T{\bf p}_n,\nabla_U F({\bf z}_n+\check Z_n)
 -\nabla_U F({\bf z}_n)-\nabla_{UY}F({\bf z}_n)\check Y_n-\nabla_{UU}F({\bf z}_n)\check U_n\rangle\T.
\end{align}
The important point in this equation is that the matrix of the left-hand side of \eqref{ZGlUn} is independent of $\check Z$ and that the right-hand side \eqref{ResU} has a Lipschitz constant $O(\epsilon)$ in an $\epsilon$-neighborhood of the origin as will be shown below.
\par
Combining now all equations \eqref{errYn}, \eqref{errPn}, \eqref{errPN} and \eqref{ZGlUn}, the complete system for the error $\check Z$ has the form
\begin{align}\label{MFPGl}
\MM_0\check Z
=&\begin{pmatrix}
 -\tau^Y+R^Y(\check Z)\\
 -\tau^P+R^P(\check Z)\\
 -\tau^U+R^U(\check Z)
\end{pmatrix},
\\\label{DefM0}
\MM_0:=&
\begin{pmatrix}
 M_{11}\otimes I_m&0&0\\
 M_{21}\otimes\nabla_{yy}\CC_N&M_{22}\otimes I_m&0\\
  \langle{\bf K}\T{\bf p},\nabla_{UY}F({\bf z})\rangle\T
 &\nabla_U F({\bf z})\T {\bf K}
 &\Omega
\end{pmatrix},
\end{align}
where the last block in $\MM_0$, being a symmetric block diagonal matrix, was abbreviated as
\begin{align}\label{DefOmega}
\Omega:=\langle{\bf K}\T{\bf p},\nabla_{UU}F({\bf z})\rangle=\diag_n(\Omega_n),
\ \Omega_n:=\langle K_n\T{\bf p}_n,\nabla_{UU}F({\bf z}_n)\rangle,
\end{align}
since it requires detailed investigation.
In order to avoid confusion with the coefficient $K$ from the standard scheme, we denoted the block diagonal matrix of all $K_n$ by ${\bf K}:=\diag(K_0,K,\ldots,K,K_N)$.
We note that $R^Y(0)=R^P(0)=0$, $R^U(0)=0$.
\par
For ease of reading, we recall a few details for $\MM_0$ from \cite{LangSchmitt2022}.
The index range $0\ldots N$ of the grid is also used for this matrix and its different blocks.
From \eqref{errYn} and \eqref{errPn}, we see that the submatrices $M_{11}$ and $M_{22}$ are block bi-diagonal matrices with identities $I_s$ in the main diagonal.
$M_{11}$ is lower bi-diagonal with subdiagonal blocks $(M_{11})_{n,n-1}=-\bar B_n:=-A_n^{-1}B_n,\,1\le n\le N,$ and $M_{22}$ is upper bi-diagonal with super-diagonals $(M_{22})_{n,n+1}=-\tilde B_{n+1}\T:=-(BA^{-1})\T,\,0\le n<N$.
Due to requirement \eqref{lbd2} there exist norms such that $\|\bar B_n\|=\|\tilde B_{n+1}\T\|=1$ hold.
Hence, all non-trivial blocks of the inverses $M_{11}^{-1}$ and $M_{22}^{-1}$ have norm one (with the possible exception of one single block at each boundary).
The third block $M_{21}$ is very sparse with a rank-one entry $\eins w\T$ in its last diagonal $s\times s$-block only due to \eqref{KKT_adj_peer_init}.
\par
Temporarily assuming non-singularity of $\Omega$ (which will be considered later on) the system \eqref{MFPGl} may be transformed to fixed-point form
\begin{align}\label{GSFPF}
 \check Z=\Phi(\check Z)
 =\begin{pmatrix}\Phi^Y\\\Phi^P\\\Phi^U \end{pmatrix}
 :=-\MM_0^{-1}\tau +\MM_0^{-1}R(\check Z)
\end{align}
with the vectors
\begin{align*}
\tau:=\begin{pmatrix}\tau^Y\\ \tau^P\\\tau^U\end{pmatrix},\quad
R(\check Z):=\begin{pmatrix}R^Y(\check Z)\\ R^P(\check Z)\\R^U(\check Z)\end{pmatrix}.
\end{align*}
Now, the inverse of the fixed matrix $\MM_0$ may be given in factored form as
\begin{align}\label{FakMinvU}
 \MM_0^{-1}=&\begin{pmatrix}
 I&\\
 0&I&\\
 -\Omega^{-1}\langle{\bf K}\T{\bf p},\nabla_{UY}F({\bf z})\rangle\T
 &-\Omega^{-1}\nabla_U F({\bf z})\T {\bf K}&\Omega^{-1}
 \end{pmatrix}\\\label{FakMinvYP}
 &\cdot
 \begin{pmatrix}
 M_{11}^{-1}\otimes I_m\\
 -M_{22}^{-1}M_{21}M_{11}^{-1}\otimes \nabla_{yy}\CC_N&M_{22}^{-1}\otimes I_m\\
 0&0&I
 \end{pmatrix}.
\end{align}
Since there is no coupling between the stages of different time steps in $R(\check Z)$, all Jacobians $\nabla_Y R,\nabla_P R,\nabla_U R$ for each part $R^Y,R^P,R^U$ are block diagonal matrices with blocks of size $(sm)\times(sm)$ or $(sd)\times(sd)$,
since the two-step structure of the Peer methods is completely represented in the block bi-diagonal matrices $M_{11},M_{22}$.
Since every single nontrivial $s\times s$-block of the block triangular inverses $M_{11}^{-1},M_{22}^{-1}$ has norm one according to the discussion following equation \eqref{DefOmega}, the inverse \eqref{FakMinvYP} exists with a norm of size $O(N)=O(h^{-1})$.
Together with smoothness requirements on the function $f$ Lipschitz, estimates for the $\check Y,\check P$ equations in \eqref{GSFPF} will follow below by standard arguments.
However, the last part with $R^U$, lacking the factor $h$, which is essentially covered by the left factor \eqref{FakMinvU} of $\MM_0^{-1}$, needs additional inspection.
\par
Obviously, boundedness of this factor \eqref{FakMinvU} requires that the block diagonal matrix $\Omega$ from \eqref{DefOmega} has a uniformly bounded inverse.
Since $p(t)$ is assumed to be smooth, Lemma~\ref{Lknspf} shows that $K_n\T{\bf p}_n=D_{K_n}{\bf p}_n+O(h^q)$, where $D_{K_n}\succ0$ is the diagonal matrix with the row sums of $K_n\T$.
Hence, $D_{K_n}^{-1}\Omega_n$ is a small perturbation of a block diagonal matrix with blocks being the Hesse matrices $\nabla_{uu} H({\bf y}_{ni},{\bf p}_{ni},{\bf u}_{ni})=\langle {\bf p}_{ni},\nabla_{uu}f({\bf y}_{ni},{\bf u}_{ni})\rangle,\,1\le i\le s,$ of the Hamiltonian $H(y,p,u)=p\T f(y,u)$ at the solution.
Now, the {\em control-uniqueness property} from \cite{Hager2000} assumes that for any $t\in[0,T]$ the Hamiltonian $H\big(y(t),p(t),\tilde u\big)$ has a unique minimum with respect to $\tilde u$ in small neighborhoods of $u(t)$.
An appropriate condition for this property is the definiteness of the Hessian
\begin{align}\label{DefHess}
 \nabla_{uu} H\big(y(t),p(t),u(t)\big)\succ \eta I_d,\ t\in[0,T],\ \eta>0,
\end{align}
implying bounded invertibility of this Hessian which is not essentially affected by small perturbations of $p(t)$.
In the main theorem below, we will use a slightly weaker assumption \eqref{Vorausfuu}, but we will show now that \eqref{DefHess} is satisfied for an interesting class of control problems.
\par\noindent
\begin{example}\label{Ectrprb}
A common type of optimal control problems of tracking type has right-hand sides $f(y,u)$ which depend linearly on $u\in\R^d$.
Only the objective function is quadratic in the form \eqref{objfunc-lagrange},
\begin{align}\label{ZFInt}
 \frac12\int_0^T(y\T\Upsilon y+u\T{W} u)dt
\end{align}
with positive definite matrices $\Upsilon,{W}\succ0$.
The transformation to standard form \eqref{OCprob_objfunc} uses the additional differential equation
\begin{align}\label{ZDglm}
 y_{m+1}'=\frac12\big(y\T\Upsilon y+u\T{W} u),\ y_{m+1}(0)=0,
\end{align}
and \eqref{ZFInt} becomes $\CC(\bar y(T))=y_{m+1}(T)$ with extended $\bar y\T=(y\T,y_{m+1})$.
Also, $\bar f$, $\bar p$ are extended versions.
Since the right-hand side of \eqref{ZDglm} does not depend on $y_{m+1}$, the adjoint equation for the last Lagrange multiplier simply reads $p_{m+1}'=0$ with end condition $p_{m+1}(T)=1$ yielding $p^\star_{m+1}(t)\equiv 1$.
Hence, $\langle p^\star,\nabla_{uu}H(z^\star)\rangle=$ $\langle\bar p^\star,\nabla_{uu}\bar f(y^\star,u^\star)\rangle={W}\succ 0$ is definite.
\par
We note that here even the $m+1$-th component of the discrete solution $\bar P_n$ is exact.
Denoting the stage vector of its $m+1$-th component by $\phi_n\in\R^s$, one sees that the first column of \eqref{OBaEnd} simply states $A_N\T V_se_1=A_N\T\eins=w$ for $r\ge s-1$ and that the end condition \eqref{KKT_adj_peer_init} reads $A_N\T\phi_N=w$.
Hence $\phi_N=\eins$ is the unique solution due to the non-singularity of $A_N$.
Then, the adjoint recursion \eqref{KKT_adj_peer} reduces to $A_n\T \phi_n=B_{n+1}\T\phi_{n+1}$, which leaves the vector $\eins$ unchanged by \eqref{OBaStd}.
Hence, $\phi_n\equiv\eins$ for $n=0,\ldots,N$.
Since the original right-hand side $f$ is linear in $u$, condition \eqref{KKT_ctr_peer} depends on $u$ only in the $m+1$-th component of $\bar F$ and its derivative with respect to $U$ is
 $D_{K_n}\otimes{W}\succ0$
by \eqref{pos-cond}.
We remind that for \texttt{AP4o43p} equations from the third stage with $\kappa_{33}=0$ have been removed from the system, see Remark~\ref{Rblst}.
Hence, the corresponding diagonal block in Newton's method is non-singular.
\end{example}
\par
For the error estimates below norms are required on three different levels.
On the highest, the grid level, the maximum norm is used for convenience.
On the step level it is essential to use appropriate weighted norms for $\check Y_n,\check P_n$ such that $\|\bar B\|=\|\tilde B\T\|=1$ holds.
On the lowest, the problem level, any norm may be appropriate.
If, for instance, \eqref{DefHess} is given, the Euclidean norm may be considered for $\check U_{ni}\in\R^d$.
However, in the following theorem, we will use the slightly more general assumption \eqref{Vorausfuu} in an appropriate norm.
In order to prove contraction in the equations for $\check Y,\check P$ in \eqref{GSFPF}, some degree of smoothness is required for the right-hand sides of both differential equations \eqref{KKT_state}, \eqref{KKT_ctr}.
We assume that there exist constants $\Lambda_j,\,j=1,2,3,$ such that
\begin{align}\label{consLbd}
 \|\nabla_z f(z)\|\le\Lambda_1,\quad
 (1+\|z\|)\|\nabla_{zz} f(z)\|\le\Lambda_2,\quad
 \|\langle p,\nabla_{uzz}f(z)\rangle\|\le\Lambda_3,
\end{align}
in some tubular neighborhood of the solution $z^\star=((y^\star)\T,(p^\star)\T,(u^\star)\T)\T$.
The constant $\Gamma:=\|\nabla_{yy}\CC\|$ vanishes for linear objective functions.
One further detail concerns the triplet \texttt{AP4o33pa} having a node larger than one, which means that the last off-step node $t_{Ns}$ in the grid exceeds $T$.
Hence, the smoothness assumptions on the solution are required in a slightly larger interval $[0,T^*]\supseteq[0,T]$.
\begin{theorem}\label{TExstGlobal}
Considering the unconstrained case with $N_U=\{0\}$, let the objective function $\CC$ in \eqref{OCprob_objfunc} be a polynomial of degree less or equal two and assume
\begin{align}\label{Vorausfuu}
 \|\langle \tilde p,\nabla_{uu}f\big(y^\star(t),u^\star(t)\big)\rangle^{-1}\|\le \omega
\end{align}
for all $\tilde p$ with $\|\tilde p-p(t)\|\le\epsilon<1$, $t\in[0,T^*]$, $T^*>T$.
Assume also that a unique solution $\big(y^\star(t),p^\star(t),u^\star(t)\big)$ of \eqref{KKT_state}--\eqref{KKT_ctr} exists with $y^\star,p^\star,u^\star\in C^q[0,T^*]$.
Also let the right-hand side $f$ of \eqref{OCprob_ODE} satisfy \eqref{consLbd} such that
\begin{align}\label{Glatth1}
(\zeta_1^P\Lambda_1+\zeta_2^P\Lambda_2+\zeta_c^P\Gamma\Lambda_1)T\le\frac23,\\
\label{Glatth2}
(\zeta_1^U\Lambda_1+\zeta_2^U\Lambda_2)\Lambda_1\omega\le\frac13,
\end{align}
where the constants $\zeta$ are determined by the actual triplet only and do not depend on the grid size.

\par
Let the Peer triplet satisfy the order conditions from lines (a,b,d,e,f) of Table~\ref{TOBed} with $2\le q\le r\le s$, and the eigenvalue conditions \eqref{lbd2}, \eqref{LsbRS}.
\par
Then, for $h\le h_0$ the fixed point problem \eqref{GSFPF} has a unique solution $\check Z$ in a sufficiently small tubular neighborhood of the exact solution $\big(y^\star(t),p^\star(t),u^\star(t)\big)$ of \eqref{KKT_state}--\eqref{KKT_ctr}.
The solution $\check Z$ satisfies
\begin{align}\label{FeSchranke}
 \|\check Z\|=\max\{\|Y-{\bf y}\|,\|P-{\bf p}\|,\|U-{\bf u}\|\}=O(h^{q-1})
\end{align}
where $(Y,P,U)$ is the solution of the discrete boundary value problem \eqref{KKT_state_peer_init}--\eqref{KKT_ctr_peer}.
\end{theorem}
{\it Proof:}
a)
Considering $\tilde Z,\hat Z$ in a neighborhood ${\cal N}_\ZZ:=\{\ZZ:\,\|\ZZ\|\le\varepsilon\}$ of the origin, a Lipschitz condition for $R_n^Y$ in \eqref{errYn} follows from
\begin{align}\notag
 \|R_n^Y(\tilde Z)-R_n^Y(\hat Z)\|=&\,h\|\bar K_n\big(F({\bf z}_n+\tilde Z_n)-F({\bf z}_n+\hat Z_n)\big)\|\\\label{LDifRY}
 \le&\, h\|\bar K_n\|\Lambda_1\|\tilde Z_n-\hat Z_n\|.
\end{align}
Now, $R^Y$ is multiplied by a lower block triangular matrix possessing the blocks $(M_{11}^{-1})_{nk}=\bar B_n\cdots\bar B_{k+1}$, $k<n$.
With the possible exception of $\bar B_N$, all other factors $\bar B$ have norm one in a suitable norm which exists according to \eqref{lbd2}.
So a Lipschitz condition in the equation for $\check Y$ holds as
\begin{align}\label{LipBY}
 \|M_{11}^{-1}\big(R^Y(\tilde Y)-R^Y(\hat Y)\big)\|\le& L_Y\|\tilde Z-\hat Z\|,\ L_Y= \zeta^Y\Lambda_1T,
\end{align}
since $\|M_{11}^{-1}\|\le N\|B_N\|+1\le \|B_N\|Th^{-1}$.
In \eqref{LipBY}, $\zeta^Y=\|\bar B_N\|\max_n\|\bar K_n\|$ is a constant depending on the actual triplet only, while $\Lambda_1T$ is problem-dependent.
\par
Using the formal identity $2(ab-\alpha\beta)=(a-\alpha)(b+\beta)+(a+\alpha)(b-\beta)$, the Lipschitz difference for $R^P$ may be rewritten as
\begin{align*}
\|R_n^P(\tilde Z)-R_n^P(\hat Z)\|=&\,
 h\|A_n\mT\big(\langle K_n\T({\bf p}_n+\tilde P_n),\nabla_Y F({\bf z}_n+\tilde Z_n)\rangle\\
 &\,-\langle K_n\T({\bf p}_n+\hat P_n),\nabla_Y F({\bf z}_n+\hat Z_n)\rangle\big)\|\\
 =&\,\frac{h}2\|A_n\mT\big(\langle K_n\T(\tilde P_n-\hat P_n),\nabla_YF({\bf z}_n+\tilde Z_n)+\nabla_YF({\bf z}_n+\hat Z_n)\rangle\\
 &\,+\langle K_n\T(2{\bf p}_n+\tilde P_n+\hat P_n),\nabla_YF({\bf z}_n+\tilde Z_n)-\nabla_YF({\bf z}_n+\hat Z_n)\rangle\big)\|\\
 \le&\,h(\tilde\zeta_1^P\Lambda_1+\tilde\zeta_2^P\Lambda_2)\|\tilde Z_n-\hat Z_n\|,
\end{align*}
with constants $\tilde\zeta_j^P$.
Looking again at \eqref{GSFPF}, \eqref{FakMinvYP}, the matrix $M_{22}^{-1}$ is block upper triangular with off-diagonal blocks essentially being powers of $\tilde B\T$ where $\|\tilde B\T\|=1$.
Again, only the norms in the boundary steps may be different and we get $\|M_{22}^{-1}\|\le\|\tilde B_1\T\|(N+\|\tilde B_N\T\|)\le 2\|\tilde B_1\T\|Th^{-1}$ for small $h$.
The matrix $M_{22}^{-1}M_{21}M_{11}^{-1}$ in the subdiagonal of \eqref{FakMinvYP} is a simple rank-one matrix \cite{LangSchmitt2022}.
Since the two factors of the only nontrivial last block of $(M_{21})_{NN}=\eins w\T=\eins\eins\T A_N$ are right resp. left eigenvectors of the blocks in $M_{22}^{-1}$ resp. $M_{11}^{-1}$, each $s\times s$-Block of  $M_{22}^{-1}M_{21}M_{11}^{-1}$ is equal to $\eins w\T$ and the norm of this matrix is $O(N)$ again.
Hence, in the Lipschitz difference of $\Phi^P$, we also inherit some contribution $O(\Gamma L_Y\|\tilde Z-\hat Z\|)$ from \eqref{LDifRY} yielding
\begin{align}\label{LipBP}
 \|\Phi^P(\tilde Z)-\Phi^{P}(\hat Z)\|\le& \,L_P\|\tilde Z_n-\hat Z_n\|
 \\\notag
 L_P=&\,(\zeta_1^P\Lambda_1+\zeta_2^P\Lambda_2+\zeta_c^P\Gamma\Lambda_1)T.
\end{align}
In order to simplify assumptions slightly, we will assume $\zeta_1^P\ge\zeta^Y$ yielding $L_P\ge L_Y$.
\par\noindent
b)
The map $R_n^U$ given explicitly in \eqref{ResU} consists of two parts and we will consider Lipschitz differences $R_n^U(\tilde Z)-R_n^U(\hat Z)$ for both parts separately by using Taylor's Theorem with integral remainder.
For the contribution in the first line of \eqref{ResU}, we get with 
$(\nabla_Y,0,\nabla_U)=\nabla_Z$ that
\begin{align*}
& -\langle K_n\T(\tilde P_n-\hat P_n),\nabla_U F({\bf z}_n)\rangle\T\\
 & +\langle K_n\T\tilde P_n,\nabla_U F({\bf z}_n+\tilde Z_n)\rangle\T
  -\langle K_n\T\hat P_n,\nabla_U F({\bf z}_n+\hat Z_n)\rangle\T
 \\[1mm]
 =&\;\langle K_n\T(\tilde P_n-\hat P_n),\nabla_U F({\bf z}_n+\tilde Z_n)-\nabla_U F({\bf z}_n)\rangle\T\\
 &\;+\langle K_n\T\hat P_n,\int_0^1\nabla_{UZ} F({\bf z}_n+\xi\tilde Z_n+(1-\xi)\hat Z_n)\, d\xi\,(\tilde Z_n-\hat Z_n)\rangle\T\\
 =&\;\langle K_n\T(\tilde P_n-\hat P_n),\int_0^1\nabla_{UZ} F({\bf z}_n+\xi\tilde Z_n)\,d\xi\,\tilde Z_n\rangle\T\\
 &\;+\langle K_n\T\hat P_n,\int_0^1\nabla_{UZ} F({\bf z}_n+\xi\tilde Z_n+(1-\xi)\hat Z_n)d\xi(\tilde Z_n-\hat Z_n)\rangle\T.
\end{align*}
This part is bounded by
\begin{align}
\notag
 &\|K_n\T\|\Lambda_2\|\tilde P_n-\hat P_n\|(\|\tilde Y\|+\|\tilde U\|)+L_2\|\hat P_n\|\|\tilde Z_n-\hat Z_n\|\\\label{LipB1}
 &\le\zeta^U\Lambda_2(\|\tilde Z_n\|+\|\hat Z_n\|)\|\tilde Z_n-\hat Z_n\|.
\end{align}
In the difference for the remaining term from \eqref{ResU}, the constant part $\langle K_n\T{\bf p}_n,\nabla_UF({\bf z}_n)\rangle\T$ cancels out and the others contribute
\begin{align*}
&
-\langle K_n\T{\bf p}_n,\nabla_U F({\bf z}_n+\tilde Z_n)
-\nabla_U F({\bf z}_n+\hat Z_n)\rangle\T
+\langle K_n\T{\bf p}_n,\nabla_{UZ}F({\bf z}_n)(\tilde Z_n-\hat Z_n)\rangle\T\\
=&\;
-\langle K_n\T{\bf p}_n,\int_0^1\nabla_{UZ} F({\bf z}_n+\xi\tilde Z_n+(1-\xi)\hat Z_n)\,d\xi\,(\tilde Z_n-\hat Z_n)\rangle\T\\[1mm]
&\;+\langle K_n\T{\bf p}_n,\nabla_{UZ}F({\bf z}_n)(\tilde Z_n-\hat Z_n)\rangle\T
 \\[1mm]
=&\;
-\langle K_n\T{\bf p}_n,\int_0^1\left(\nabla_{UZ} F({\bf z}_n+\xi\tilde Z_n+(1-\xi)\hat Z_n)\rangle\T
 -\nabla_{UZ}F({\bf z}_n)\right)\,d\xi\,(\tilde Z_n-\hat Z_n)\rangle\T.
\end{align*}
Obviously, Lipschitz continuity of the second derivatives would suffice to obtain the desired result.
Using the bound from assumption \eqref{consLbd}, it is seen that this contribution to the Lipschitz difference of $R_n^U$ is bounded by $\tilde\zeta_3^U\Lambda_3(\|\tilde Z_n\|+\|\hat Z_n\|)\|\tilde Z_n-\hat Z_n\|$.
Hence, we have shown that
\begin{align}\label{LipRU}
 \|R^U(\tilde Z)-R^U(\hat Z)\|\le(\tilde\zeta_2^U\Lambda_2+\tilde\zeta_3^U\Lambda_3)(\|\tilde Z\|+\|\hat Z\|)\|\tilde Z-\hat Z\|.
\end{align}
c) We conclude the proof that \eqref{GSFPF} is a contractive fixed point problem by showing that the first factor \eqref{FakMinvU} of $\MM_0^{-1}$ is bounded uniformly in $h\le h_0$.
The matrix $\Omega$ in its last block is again a block diagonal matrix with blocks $\Omega_n=\langle K_n\T{\bf p}_n,\nabla_{UU}F({\bf z}_n)\rangle$.
Since ${\bf p}_n$ contains the node values of the smooth solution $p^\star(t)$, the estimate \eqref{FeKP} in Lemma~\ref{Lknspf} shows that
\[ \|K_n\T{\bf p}_n-D_{K_n}{\bf p}_n\|=O(h^q) <\epsilon \]
for $h\le h_0$ where $D_{K_n}:=\diag_i(\eins\T K_ne_i)$ differs from $K_n\T$ in the boundary steps only, i.e., for $n=0,N$.
Hence, assumption \eqref{Vorausfuu} shows that $\max_n\|\Omega_n^{-1}\|\le\omega$ and that the first factor \eqref{FakMinvU} is uniformly bounded.
Its subddiagonal blocks add contributions with constants $O(\Lambda_2 L_Y)$, $O(\Lambda_1 L_P)$ to the Lipschitz difference of $\Phi^U$.
Ordering terms appropriately, we obtain
\begin{align}\label{LipBU}
\|\Phi^U(\tilde Z)-\Phi^u(\hat Z)\|\le(L_0^U+L_1^U\varepsilon)\|\tilde Z-\hat Z\|,\\\notag
L_0^U=(\zeta_1^U\Lambda_1+\zeta_2^U\Lambda_2)\Lambda_1\omega,
\ L_1^U=\zeta_3^U(\Lambda_2+\Lambda_3)\omega.
\end{align}
Here, $L^U=L_0^U+L_1^U\varepsilon$ depends on the radius $\varepsilon$ of the neighborhood ${\cal N}_\ZZ$.
\par
The Lipschitz constant for the whole map $\Phi$ is $L=\max\{L^Y,L^P,L_0^U+L_1^U\varepsilon\}$ and we are verifying now that the assumptions  are sufficient to show
\begin{align}\label{LipBGS}
 \|\Phi(\tilde Z)-\Phi(\hat Z)\|\le \frac23\|\tilde Z-\hat Z\|
\end{align}
for $\varepsilon$ sufficiently small.
First, \eqref{Glatth1} corresponds to $L^P\le2/3$ and also implies $L^Y\le2/3$ with an appropriate constant $\zeta_1^P\ge\zeta^Y$.
Then, by \eqref{Glatth2} we have $L_0^U\le1/3$ and $L^U=L_0^U+L_1^U\varepsilon\le2/3$ for $\varepsilon\le1/(3L_1^U)$.
\par
Finally, we choose $\tilde Z\in{\cal N}_\ZZ$ and $\hat Z=0$ where $\|\Phi(0)\|=\|\MM_0^{-1}\tau\|=O(h^{q-1})$.
Since $q\ge2$ by assumption, we may restrict $h_0\le h_1$ such that $\|\Phi(0)\|\le\varepsilon/3$ for $h\le h_0$, and from \eqref{LipBGS} follows
\begin{align*}
 \|\Phi(\tilde Z)\|\le \|\Phi(0)\|+\|\Phi(\tilde Z)-\Phi(0)\|\le \frac13\varepsilon+\frac23\varepsilon=\varepsilon.
\end{align*}
Hence, $\Phi$ maps ${\cal N}_\ZZ$ onto itself and is a contraction proving the existence of a unique fixed point $\check Z=\Phi(\check Z)\in{\cal N}_\ZZ$ and the solution $Z={\bf z}+\check Z$ of the discrete boundary value problem.
Again from \eqref{LipBGS} follows that
\begin{align}\notag
 \|\check Z\|=&3\|\Phi(\check Z)\|-2\|\check Z\|
  \le 3\|\Phi(\check Z)-\Phi(0)\|+3\|\MM_0^{-1}\tau\|-2\|\check Z\|
  \\\label{GlbFeS}
  \le&(3\frac23-2)\|\check Z\|+3\|\MM_0^{-1}\tau\|
  = 3\|\MM_0^{-1}\tau\|.
\end{align}
The assertion now follows from $\|\MM_0^{-1}\tau\|=O(h^{q-1})$, see Lemma~4.1 in \cite{LangSchmitt2022}.
\qed
\par
The global error estimate \eqref{FeSchranke} is rather pessimistic and may be improved by considering the super-convergence conditions \eqref{SupKnv}, \eqref{SupKna}.
However, since the norm in the estimate \eqref{GlbFeS} computes the maximum of all errors $\check Z,\check P,\check U$ only, the lower order $h^q$ may be verified rigorously.
\begin{lemma}
Let all assumptions of Theorem~\ref{TExstGlobal} hold with the full set of order conditions from Table~\ref{TOBed} for $2\le q\le r\le s$, $q<s$,
and let $y^\star,p^\star\in C^{q+1}[0,T^*]$.
Then, the solution $(Y,U,P)$ of \eqref{KKT_state_peer_init}--\eqref{KKT_ctr_peer} with $N_{U^s}=\{0\}$ also satisfies
\begin{align*}
 \|\check Z\|=\max\{\|Y-{\bf y}\|,\|P-{\bf p}\|,\|U-{\bf u}\|\}=O(h^{q}).
\end{align*}
\end{lemma}
{\it Proof:}
The improved error estimates for $Y$ and $P$ follow from the super-convergence effect where the leading error term of $\tau$ is canceled in the product $\MM_0^{-1}\tau$ by the conditions \eqref{SupKnv}, \eqref{SupKna} and a sufficiently fast damping of the remaining modes through assumption \eqref{lbd2} showing that $\|\MM_0^{-1}\tau\|=O(h^q)$, see (49) in \cite{LangSchmitt2022b}.
A proof by a different technique can be found in \cite{SchneiderLangHundsdorfer2018}.
Hence, \eqref{GlbFeS} now yields $\|\check Y\|,\|\check P\|=O(h^q)$.
The errors $\check U_n$ of the control variable in different time intervals are independent and \eqref{ZGlUn} may be solved for $\check U_n$.
There, $\|\tau_n^U\|=O(h^q)$ by \eqref{DeftauU}.
Taking norms we get
\begin{align}\label{FehlerU}
 \|\check U_n\|\le&\,L\omega\big(\|\tau_n^U\|+\|\check Y_n\|+\|\check P_n\|+ \|R_n^U(\check Z)\|\big)
\end{align}
with some constant $L$.
Since $R^U(0)=0$ it follows from \eqref{LipRU} and \eqref{FeSchranke} that $\|R_n^U(\check Z)\|\le\tilde L\|\check Z\|^2=O(h^{2q-2})$.
Finally, \eqref{FehlerU} yields $\|\check U_n\|=O(h^q+h^{2q-2})=O(h^q)$ for $q\ge2$.
\qed
\begin{remark}
Since the lowest order in $\|M_0^{-1}\tau\|$ dominates the error estimate \eqref{GlbFeS}, only order $h^q$ could be proven rigorously.
However, for methods satisfying the conditions \eqref{OBvStrt}, \eqref{OBvStd} and \eqref{SupKnv} with $r>q$, the standard superconvergence argument leads to $\|M_{11}^{-1}\tau^Y\|=O(h^r)$.
Hence, such methods also have global order $h^r$ for pure initial value problems for the state $y(t)$ only.
And in many problems coupling between the unknowns $y,p,u$ seems to be so weak that the improved order $r=q+1$ can also be observed in our numerical experiments following now.
\end{remark}
%
%
\section{Numerical results}\label{SecNum}
We present numerical results for the Peer triplets
\texttt{AP4o43pa}, \texttt{AP4o33pfs} and \texttt{AP4o43p} and compare them with
those obtained for our recently developed triplet \texttt{AP4o43bdf}
from \cite{LangSchmitt2022b,LangSchmitt2023a}.
We note that Matlab codes for \cite{LangSchmitt2023a} are available at
\cite{LangSchmitt2023b}.
As a well-known standard method also the symmetric fourth-order two-stage Gauss method
\cite[Table II.1.1]{HairerWannerLubich2006} is used. The latter one is implemented
along the principles in \cite{Hager2000} using intermediate time points
$t_n+c_ih$ for the control variables. The standard method \texttt{AP4o43bdf}
is based on BDF4 and its well-known stability angle is $\alpha\!=\!73.35^o$. It
satisfies the positivity requirements \eqref{pos-cond} and the additional consistency conditions
\eqref{add_cond} for $q\!=\!2$ which is one
order less than for the new Peer triplets.
Implicit Runge-Kutta methods of Gauss type are
symplectic making them suitable for optimal control
\cite{HairerWannerLubich2006,SanzSerna2016}. However, as all one-step methods they may suffer from order reduction due to their lower stage order $s+1$ compared to the classical order $p\!=\!2s$.
\par
To illustrate the order of convergence, we first consider two unconstrained problems with known analytic
solutions. The first one is a quadratic problem with a mixed term taken from \cite{Hager1976,Hager2000} and
the second one comes from a method-of-lines discretization of a boundary control problem
for the 1D heat equation \cite{LangSchmitt2023a}. Finally, we apply our novel Peer
methods to an optimal control problem for an 1D semilinear reaction-diffusion
model of Schl\"ogl type with cubic nonlinearity, which was intensively studied in
\cite{BuchholzEngelKammannTroeltzsch2013}. We pick the problem of stopping a nucleation process
to show the potential of higher-order methods.
\par
All calculations have been done with Matlab-Version R2019a
on a Latitude 7280 with an i5-7300U Intel processor at 2.7 GHz. We use \texttt{fmincon} with
stop tolerance $10^{-14}$. If not otherwise stated, we apply the \texttt{interior-point} algorithm as
default choice in \texttt{fmincon} and provide the zero control vector as initial guess.
\subsection{A quadratic problem with a mixed term}
The first problem is taken from \cite{Hager2000}. It was originally proposed in
\cite[(P2)]{Hager1976} and includes a mixed term $y_1(t)u(t)$. We consider
\begin{align*}
\mbox{minimize } &\,\frac12\int_{0}^{1}
\left( 1.25\,y_1(t)^2+y_1(t)u(t)+u(t)^2\right)\,dt\\[1mm]
\mbox{subject to } &\,y_1'(t) = 0.5\,y_1(t)+u(t),\quad t\in (0,1],\\
&\, y_1(0) = 1,
\end{align*}
with the optimal solution
\begin{align*}
y_1^\star(t) = \frac{\cosh(1-t)}{\cosh(1)},\quad
u^\star(t) = -\frac{(\tanh(1-t)+0.5)\cosh(1-t)}{\cosh(1)}.
\end{align*}
The optimal costate can be computed from $p_1^\star(t) = -0.5(y_1^\star(t)+2u^\star(t))$.
Introducing a second component $y_2(t)$ and setting $y'_2(t)=1.25\,y_1(t)^2+y_1(t)u(t)+u(t)^2$
with the initial value $y_2(0)=0$, the objective function can be transformed
to the Mayer form $\CC(y(1))=0.5\,y_2(1)$ with the new state vector $y=(y_1,y_2)\T$.
\begin{figure}[t!]
\centering
\includegraphics[width=6.8cm]{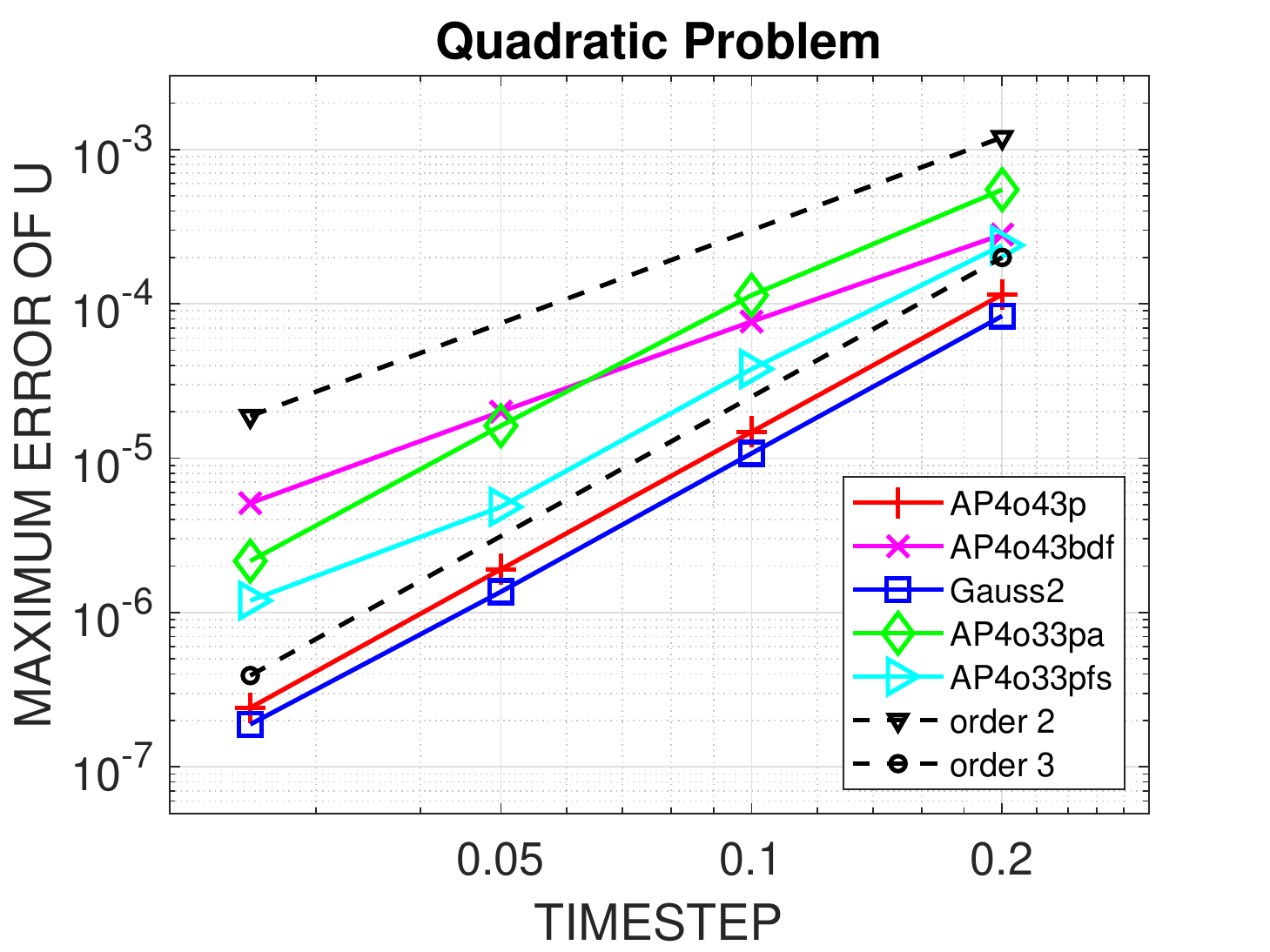}
\includegraphics[width=6.8cm]{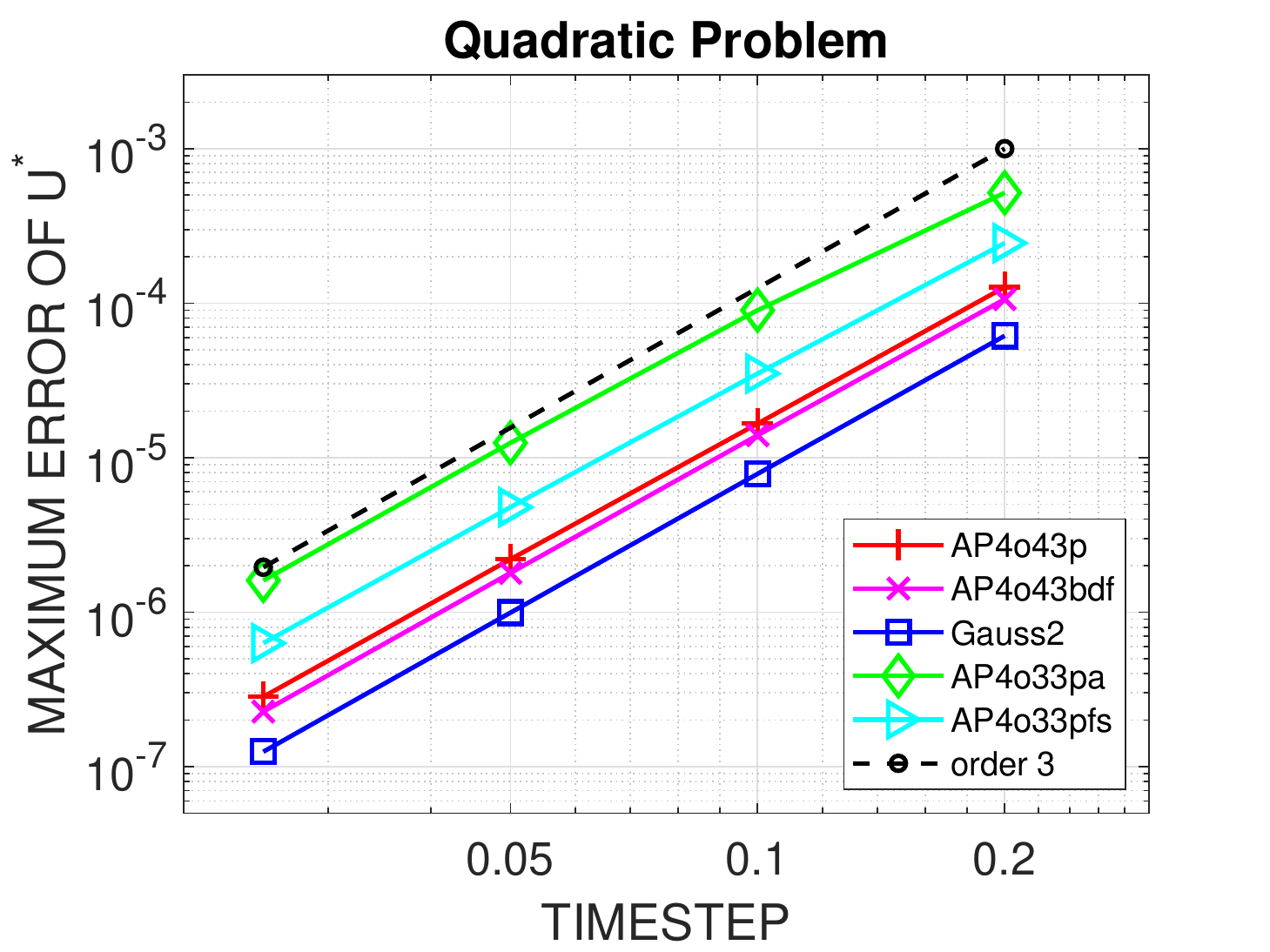}\\[5mm]
\includegraphics[width=6.8cm]{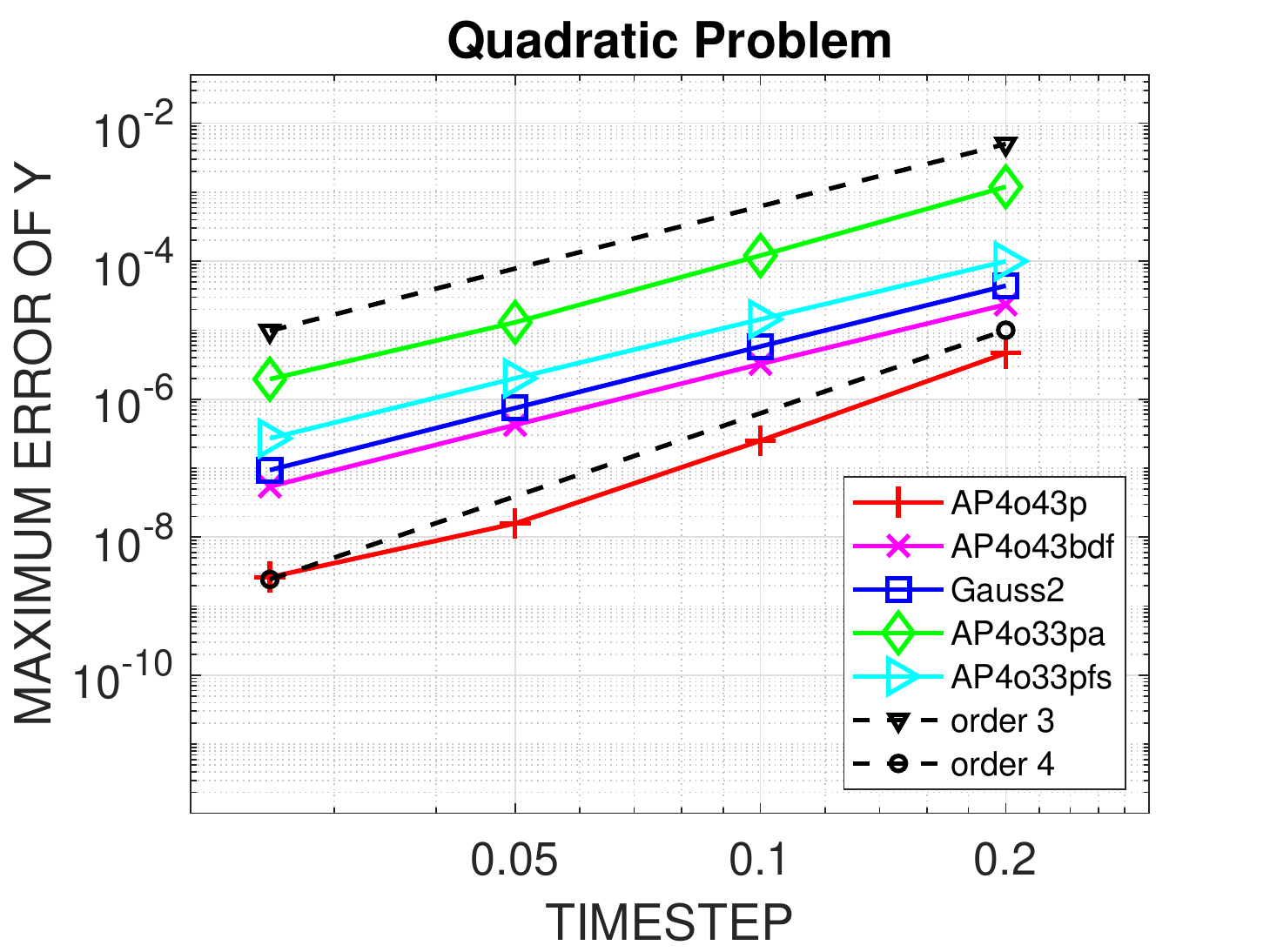}
\includegraphics[width=6.8cm]{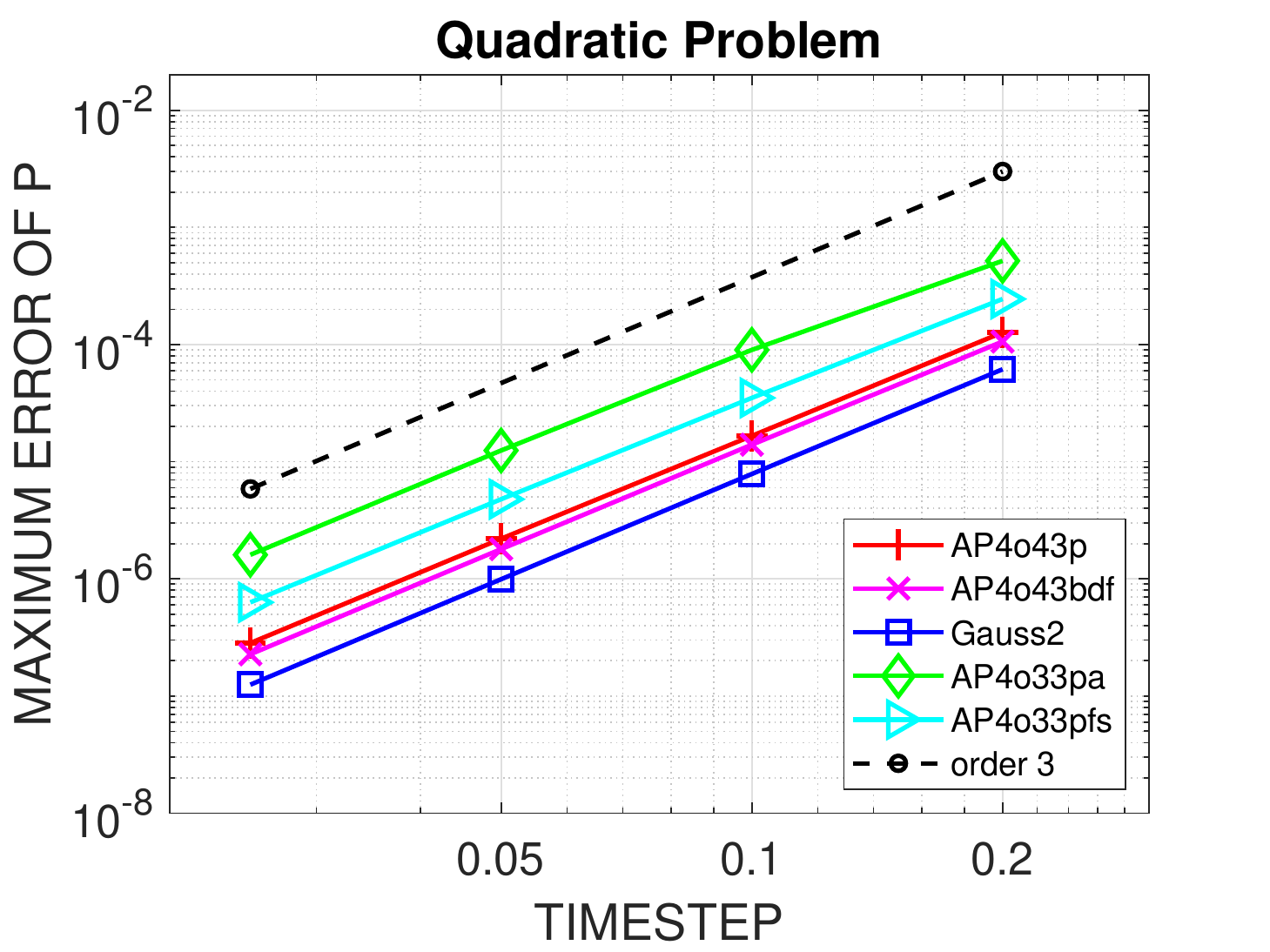}
\parbox{12.5cm}{
\caption{Quadratic Problem. Convergence of the maximal control
errors $\|U_{ni}-u(t_{ni})\|_\infty$ (top left), improved
control errors $\|U^\ddagger_{ni}-u(t_{ni})\|_\infty$ (top right), 
state errors $\|Y_{ni,1}-y_1(t_{ni})\|_\infty$
(bottom left), and adjoint errors $\|P_{ni,1}-p_1(t_{ni})\|_\infty$, (bottom right),
$n=0,\ldots,N,\;i=1,\ldots,s$.
\label{fig:qproblem}}
}
\end{figure}
Numerical results for $N\!+\!1\!=\!5,10,20,40$ are shown in Figure~\ref{fig:qproblem}.
All Peer methods show their theoretical order three for the first component of the adjoint variables, $P_{ni,1}$.
The fourth-order Gauss-2 method drops down to order three due to its lower stage order three.
Order three is also observed for the state variables $Y_{ni,1}$, except for
\texttt{AP4o43p} which achieves its super-convergence order four for the first three runs.
For the Peer methods, the errors of the control vector $U$ as well as the improved
control $U^\ddagger$ obtained from the post-processing in \eqref{ctrpp} decrease with order three as expected.
Since \texttt{AP4o43bdf} satisfies the consistency conditions
\eqref{add_cond} for $u$ with $q\!=\!2$ only, its order in $U$ is two, which is nicely seen.
However, the third-order approximations in the first components of $(Y,P)$ yields order
three for $U^\ddagger$ again. Both methods, \texttt{AP4o33pa} and \texttt{AP4o33pfs}, fall
behind the other ones in terms of accuracy. This is not surprising since their better
stability properties and the dense output feature of the latter one come with larger error constants.
\subsection{Boundary control of an 1D discrete heat equation}
The second problem is taken from \cite{LangSchmitt2023a}. It was especially
designed to provide exact formulas for analytic solutions of an optimal
boundary control problem governed by a one-dimensional discrete heat equation
and an objective function that measures the distance of the final state from
the target and the control costs. Since no spatial discretization errors are
present, numerical orders of time integrators can be observed with high accuracy
without computing reference solutions.
\par
\begin{figure}[t!]
\centering
\includegraphics[width=6.8cm]{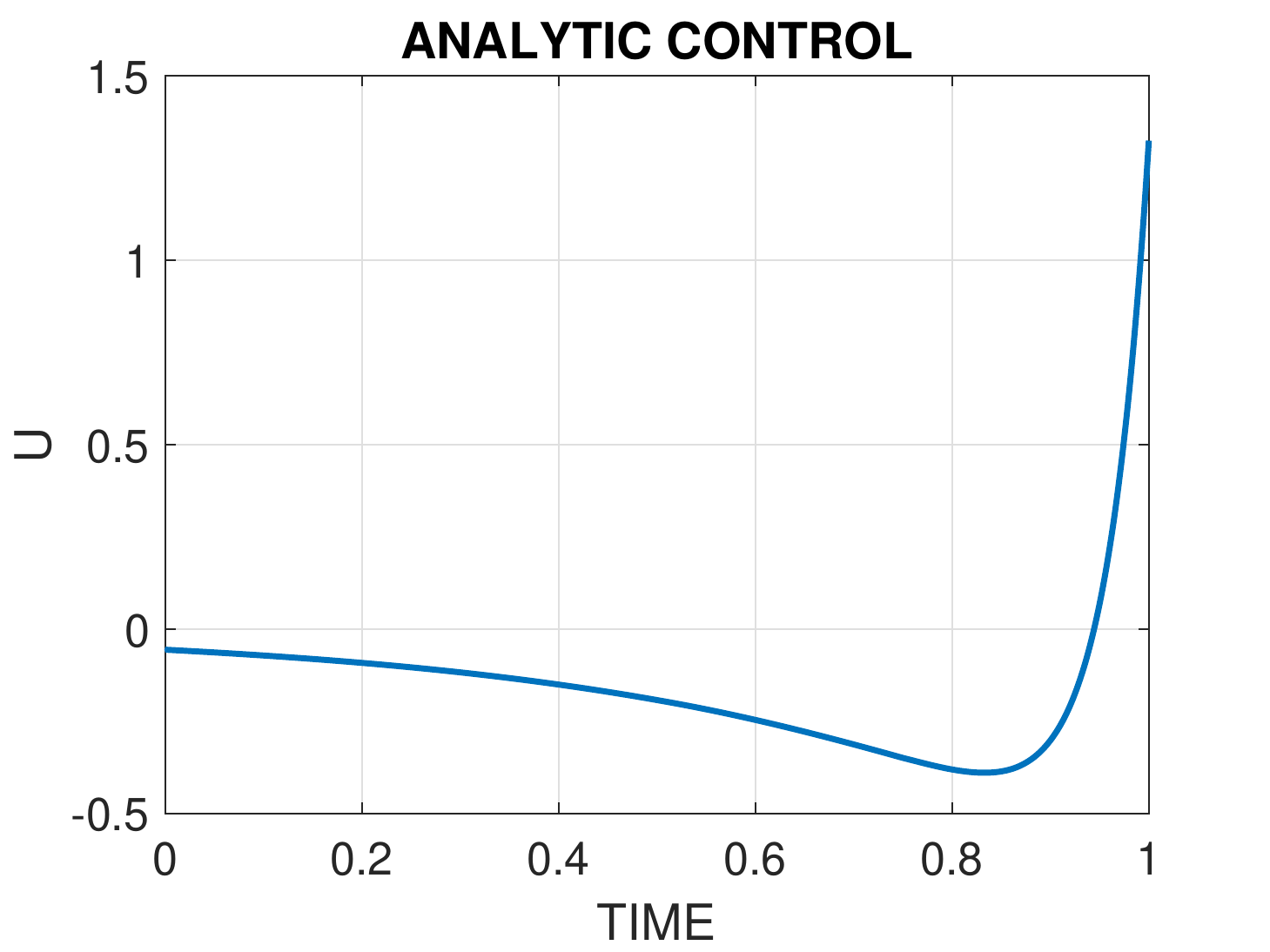}
\includegraphics[width=6.8cm]{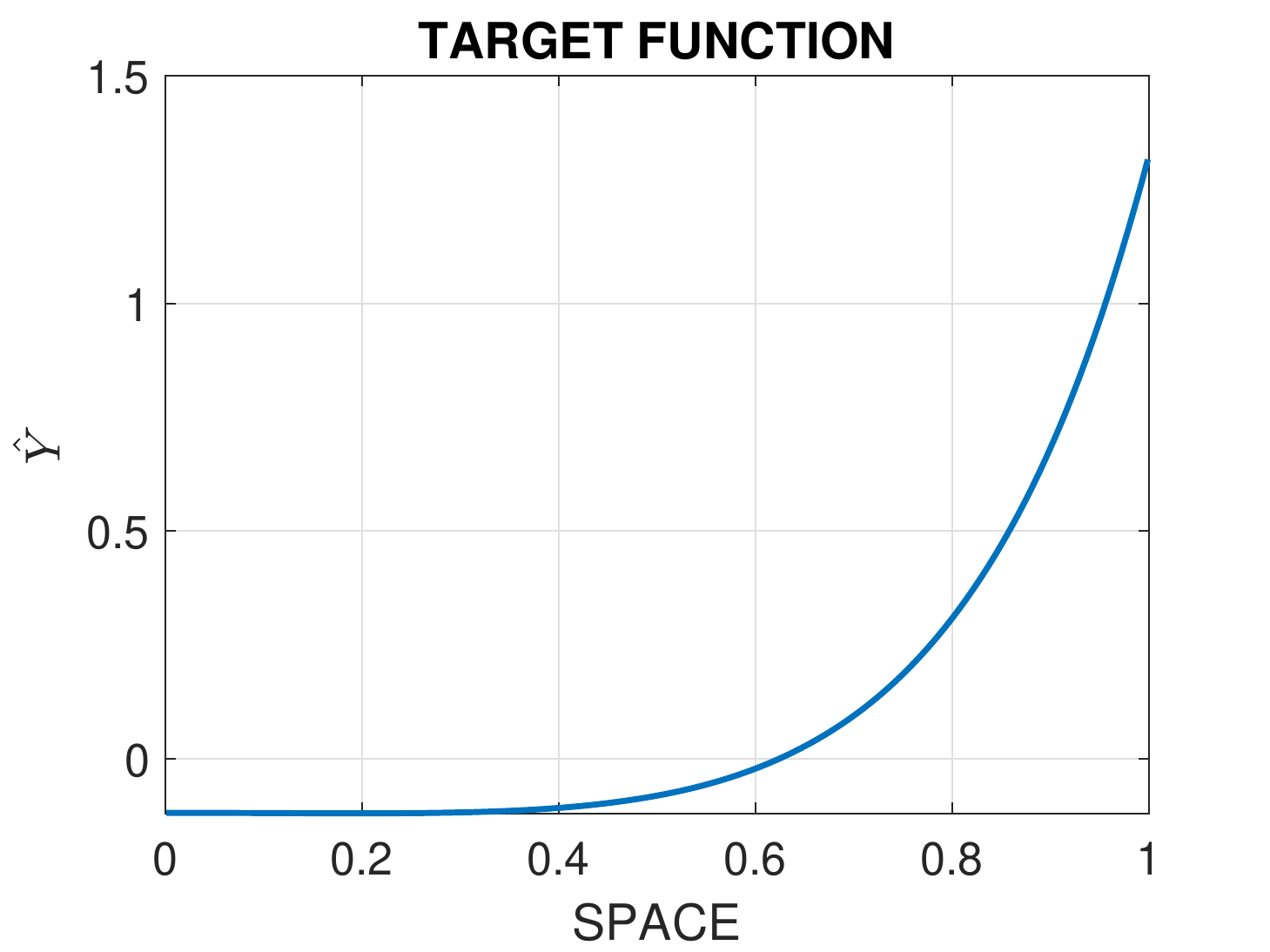}
\parbox{12.5cm}{
\caption{Dirichlet heat problem with $m=500$ spatial points.
Analytic control $u(t)$ (left) and target function $\hat{y}$ (right).
\label{fig:hprobData}}
}
\end{figure}

The optimal control problem reads as follows:
\begin{align*}
\mbox{minimize } &\,\frac12\|y(1)-\hat{y}\|^2_2 +
\frac{1}{2}\int_{0}^{1}u(t)^2\,dt\\[1mm]
\mbox{subject to } &\,y'(t) = Ay(t) +  \gamma e_m u(t),\quad t\in (0,1],\\
&\, y(0) = \ones,
\end{align*}
with
\begin{align*}
 A =&\frac1{(\triangle x)^2}\begin{pmatrix}
   -1&1\\
   1&-2&1\\
   &&\ddots&\ddots&\ddots\\
   &&&1&-2&1\\
   &&&&1&-3
 \end{pmatrix},
\end{align*}
state vector $y(t)\in\R^m$, $\triangle x=1/m$, and
$\gamma=2/(\triangle x)^2$. The components $y_i(t)$
approximate the solution of the continuous 1D heat equation $Y(x,t)$ over the spatial domain
$[0,1]$ at the discrete points $x_i=(i-0.5)\triangle x$, $i=1,\ldots,m$. The corresponding
boundary conditions are $\partial_xY(0,t)=0$ and $Y(1,t)=u(t)$. The matrix $A\in\R^{m\times m}$ 
results from standard central finite differences. Its eigenvalues $\lambda_k$ and
corresponding normalized orthogonal eigenvectors $v^{[k]}$ are given by
\begin{align*}
\lambda_k =&\;-4m^2\sin^2\left( \frac{\omega_k}{2m}\right),\;
\omega_k=\left( k-\frac12 \right)\pi,\\[2mm]
v^{[k]}_i=&\;\nu_k\cos\left( \omega_k\frac{2i-1}{2m}\right),\;
\nu_k = \frac{2}{\sqrt{2m+\sin(2\omega_k)/\sin(\omega_k/m)}},\;\; i,k=1,\ldots,m.
\end{align*}
We follow the test case in \cite{LangSchmitt2023a} and prescribe the sparse control
\begin{align*}
u^\star(t) = -\gamma p_m(t),\quad
p^\star(t) =&\; \delta \left( e^{\lambda_1(T-t)}v^{[1]} +  e^{\lambda_2(T-t)}v^{[2]}\right),
\end{align*}
with $\delta=1/75$, which defines the target vector $\hat{y}$ through
\begin{align*}
\hat{y}(t) =&\,y^\star(T)-\delta\left( v^{[1]} + v^{[2]}\right).
\end{align*}
The coefficients $\eta_k(T)$ of
$y^\star(T)=\sum_{i=1,\ldots,m}\eta_k(T)v^{[k]}$
are given by
\begin{align*}
\eta_k(T) =&\;e^{\lambda_kT}\eta_k(0)-\gamma^2\delta Tv_m^{[k]}\sum_{l=1}^{2}
v_m^{[l]}\varphi_1((\lambda_k+\lambda_l)T)
\end{align*}
where $\eta_k(0)=y^\star(0)\T v^{[k]}$ and $\varphi_1(z):=(e^z-1)/z$. We will compare
the numerical errors for $y(T)$, $p(0)$ and $u(t)$
An approximation $p_h(0)$ for the Peer method is obtained from $p_h(0)=(v\T\otimes I)P_{0}$ with
$v=V_s\mT e_s$, $e_s=(0,\ldots,0,1)\T\in\R^s$.
Note that, compared to \cite{LangSchmitt2023a}, we have changed the sign of the
adjoint variables, i.e., $p\mapsto -p$, to fit into our setting. Introducing an additional
component $y_{m+1}(t)$ and adding the equations $y'_{m+1}(t)=u(t)^2$, $y_{m+1}(0)=0$,
the objective function can be transformed to the Mayer form
\[ \CC(y(1)):=\frac12\,\left( \sum_{i=1}^{m}(y_i(1)-\hat{y}_i)^2+y_{m+1}(1)\right)\]
with the extended
vector $\bar y(1)=(y_1(1),\ldots,y_m(1),y_{m+1}(1))\T$.
We set $m=500$. In Fig.~\ref{fig:hprobData}, the analytic control $u(t)$ and the target function $\hat{y}\in\R^m$ are shown.
\begin{figure}[t!]
\centering
\includegraphics[width=6.8cm]{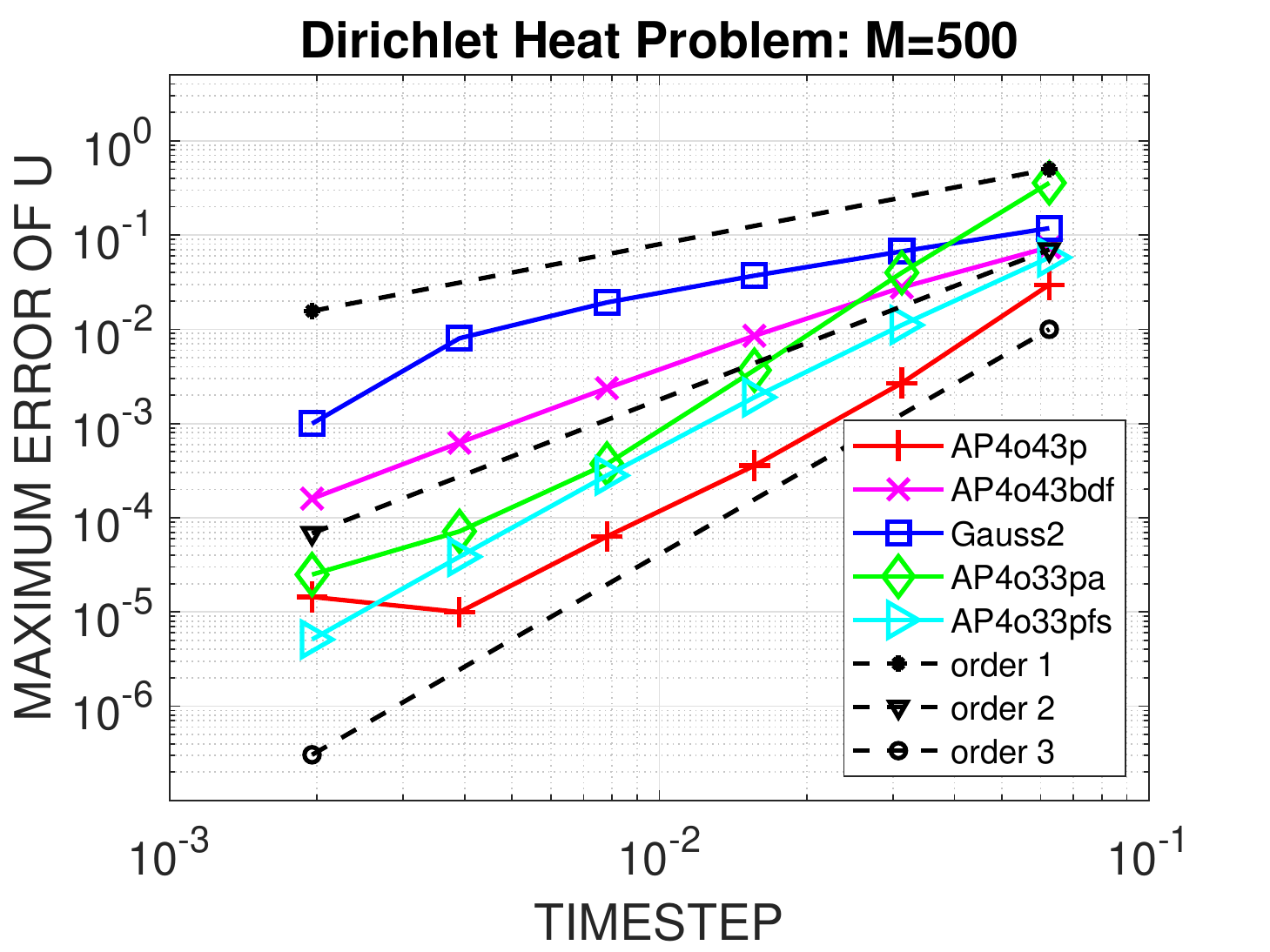}
\includegraphics[width=6.8cm]{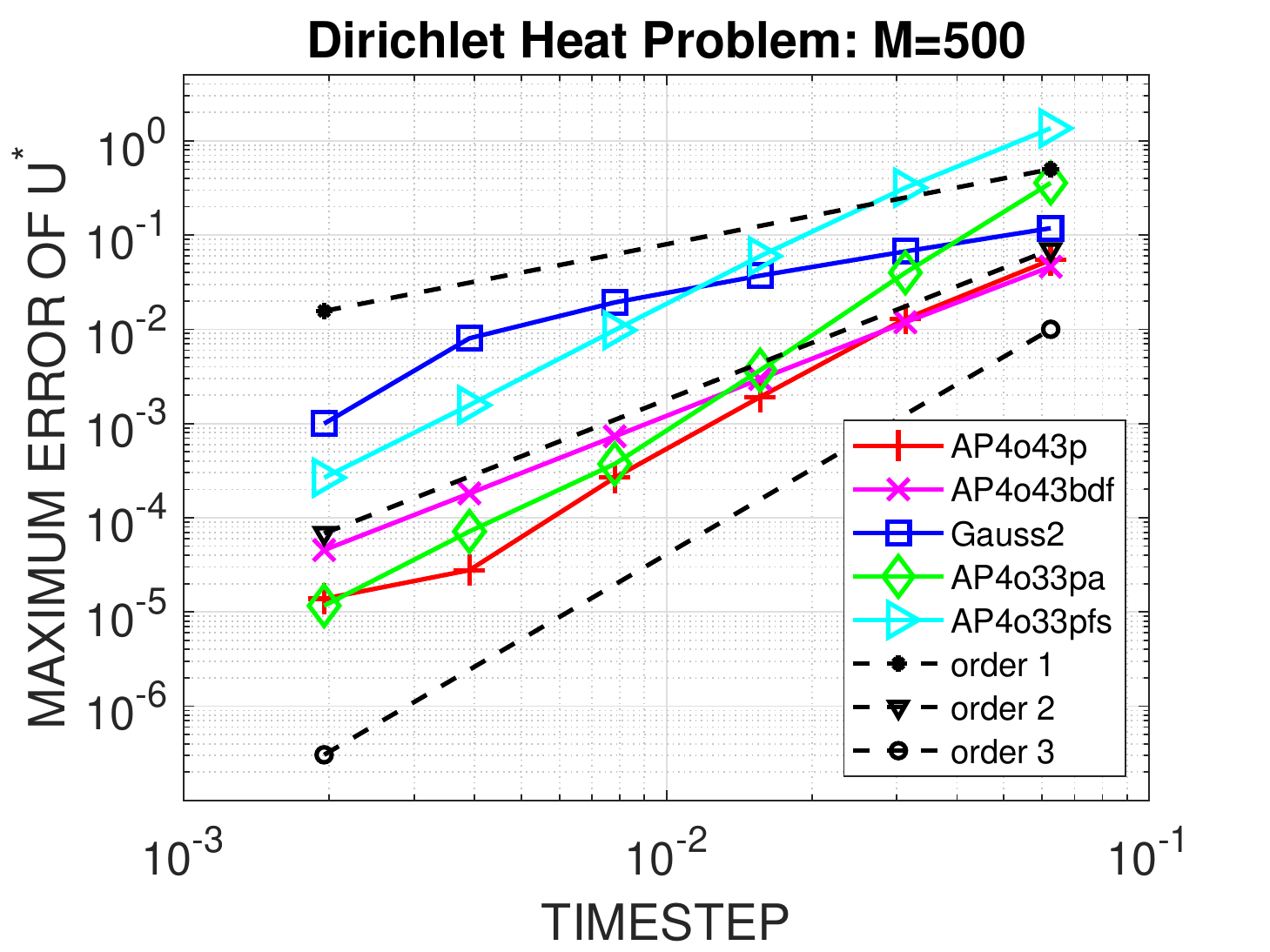}\\[5mm]
\includegraphics[width=6.8cm]{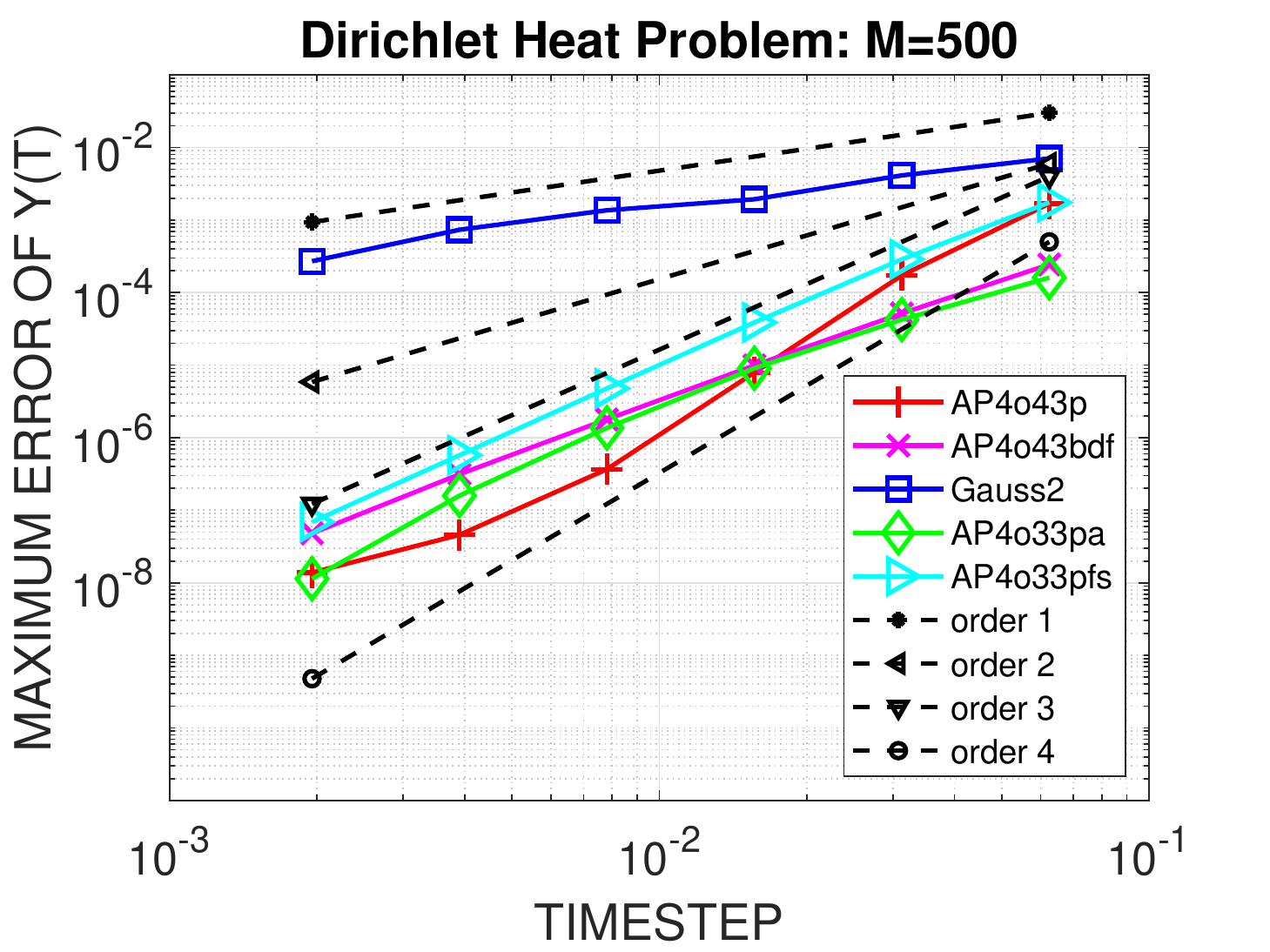}
\includegraphics[width=6.8cm]{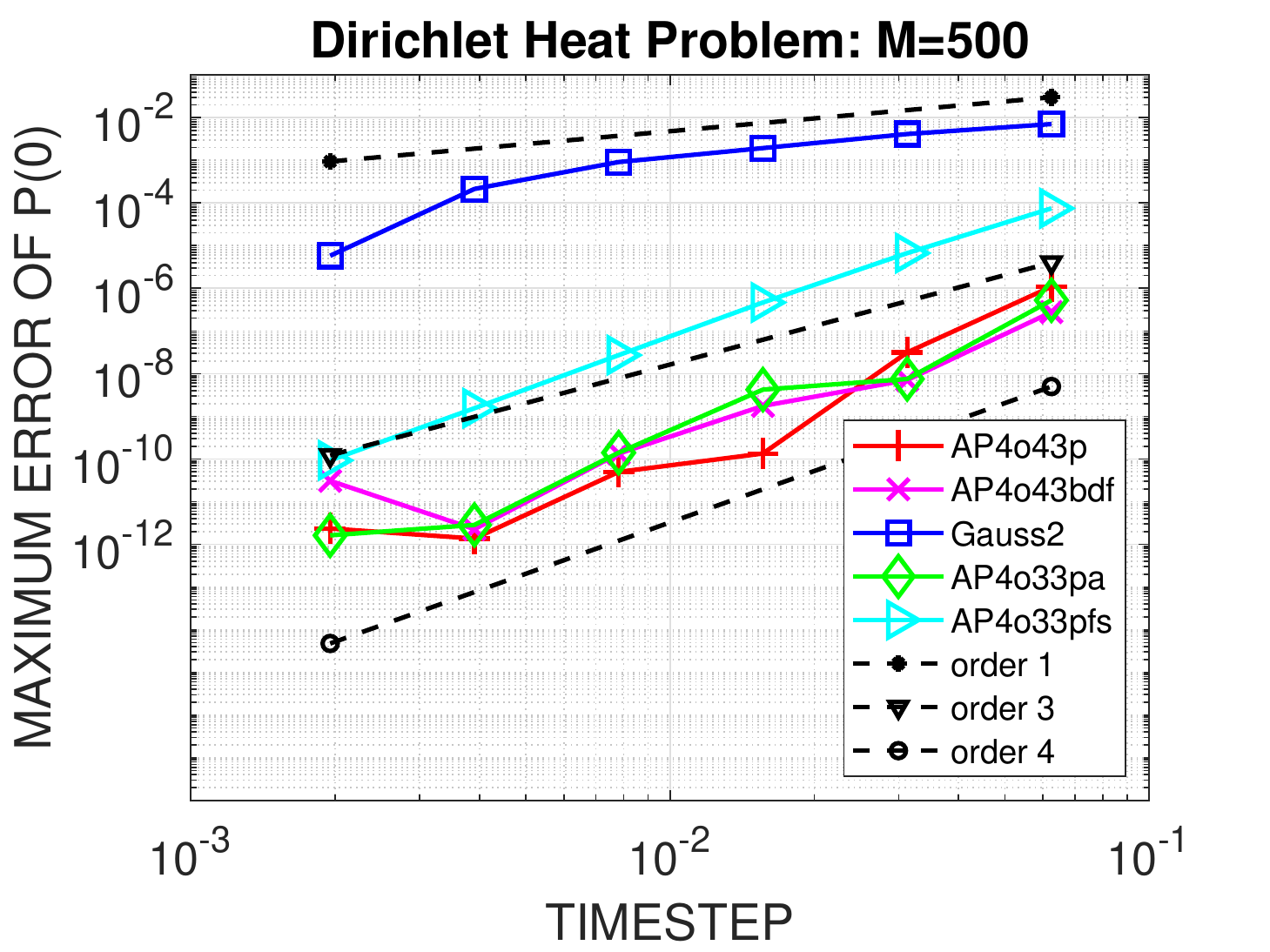}
\parbox{12.5cm}{
\caption{Dirichlet heat problem with $m=500$ spatial points. Convergence of
the maximal control errors $\|U_{ni}-u(t_{ni})\|_\infty$ (top left), improved
control errors $\|U^\ddagger_{ni}-u(t_{ni})\|_\infty$ (top right), state errors
$\|y(T)$$-$$y_h(T)\|_\infty$ (bottom left), and adjoint errors
$\|p(0)$$-$$p_h(0)\|_\infty$ (bottom right).
\label{fig:hproblem}}
}
\end{figure}
We will now discuss the numerical errors obtained from applying
$N+1=2^k$, $k=4,\ldots,9$, time steps. The results are visualized
in Fig.~\ref{fig:hproblem}. As already observed in \cite{LangSchmitt2023a},
the one-step Gauss method of order four suffers from a serious order
reduction to first order in all variables $(y,p,u)$. This phenomenon
is well understood and occurs particularly drastically for time-dependent
Dirichlet boundary conditions \cite{OstermannRoche1992}. This drawback is
shared by all one-step methods due to their insufficient stage order. The
BDF-based \texttt{AP4o43bdf} shows second order convergence in the control,
which is in accordance with the fact that it satisfies \eqref{add_cond}
with $q\!=\!2$ only. This also limits
the accuracy of the state and the adjoint at their endpoints to order
three and two, respectively. Thus, the improvement in the post-processed
control variables $U^\ddagger_{ni}$ from \eqref{ctrpp} is only marginal 
and does not increase the order.
All new Peer methods satisfy \eqref{add_cond} for $q\!=\!3$ and deliver
approximations of the control with order three, except for certain irregularities
in the smallest step. The order of convergence for $y_h(T)$ is three for
\texttt{AP4o33pa} and \texttt{AP4o33pfs}, whereas \texttt{AP4o43p} reaches
fourth-order super-convergence for nearly all time steps. For $p_h(0)$,
\texttt{AP4o33pfs} shows an ideal order three. The other two Peer methods
vary between order three and five and stagnate at the end when errors are
already quite small. The supposed improvement in $U^\ddagger_{ni}$ is not achieved
for the Peer methods. Quite to the contrary, \texttt{AP4o33pfs} loses two
orders of magnitude in accuracy, \texttt{AP4o43p} loses one order. We infer
that for Peer methods, which perform close to their theoretical order,
the approximation quality of $U$ is nearly optimal and post-processing
is not advisable in general. To summarize, the newly constructed Peer methods
significantly improve the approximation of the control with increased order three.
\texttt{AP4o43p} gains from its super-convergence property and performs
remarkably well for this discrete heat problem.
\subsection{Stopping of a nucleation process with distributed control}
\begin{figure}[t]
\centering
\includegraphics[width=6.8cm]{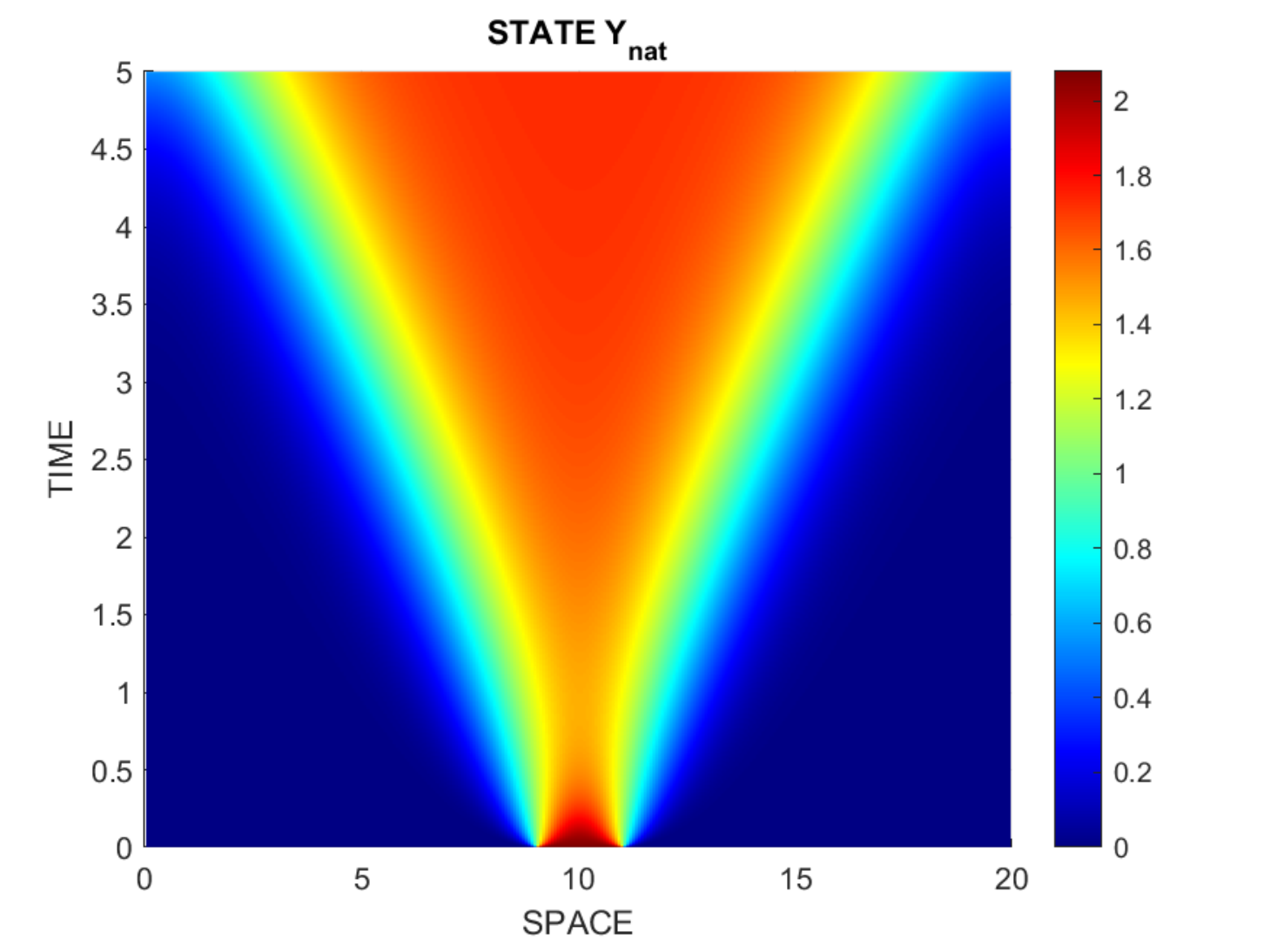}
\hspace{0.1cm}
\includegraphics[width=6.8cm]{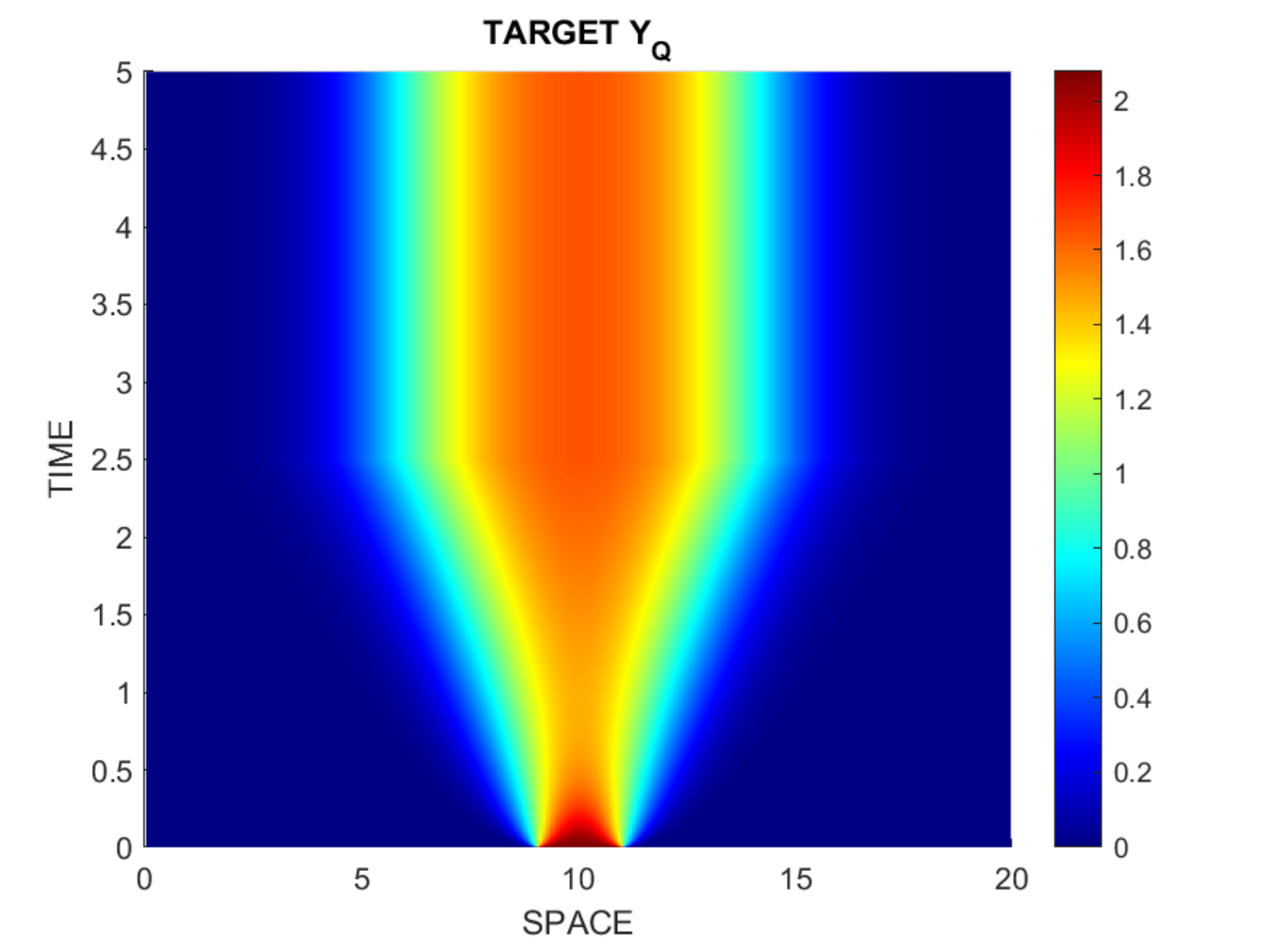}\\[2mm]
\includegraphics[width=6.8cm]{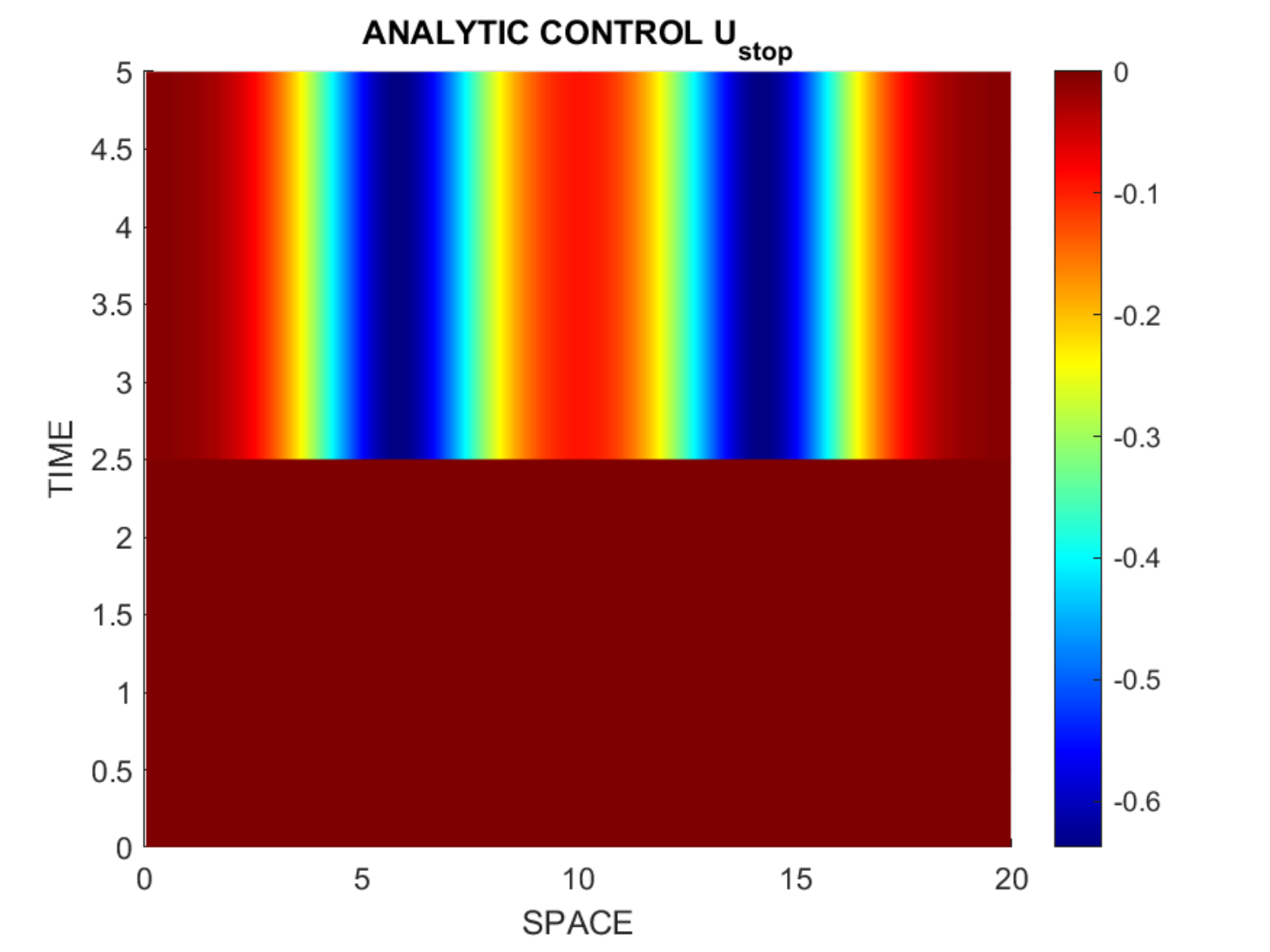}
\parbox{12.5cm}{
\caption{Nucleation process. $Y_{nat}$ for $U=0$ (top, left),
target function $Y_Q$ (top, right) and analytic control $U_{stop}$ (bottom).
\label{fig:nucl_func}}
}
\end{figure}

In our third study, we consider a PDE-constrained optimal control problem
from \cite[Chapter 5.4]{BuchholzEngelKammannTroeltzsch2013} -- stopping of a
nucleation process modelled by a nonlinear reaction-diffusion equation
of Schl\"ogl-type. It reads
\begin{align*}
\mbox{minimize } J:= &\,\frac12 \int_Q (Y(x,t)-Y_Q(x,t))^2\,dxdt +
\frac{\alpha}{2} \int_Q U(x,t)^2\,dxdt\\[2mm]
\mbox{subject to }\hspace*{0.1cm} \partial_tY-\partial_{xx}Y=& \, Y - kY^3 + U(x,t),
\;\; (x,t)\in Q:=(0,L)\times (0,T],\\
Y(x,0) =& \,Y_0(x),\;\; x\in (0,L),\\
\partial_x Y(0,t) =& \,\partial_x Y(L,t) = 0,
\end{align*}
with parameters $\alpha\!=\!10^{-6}$, $L\!=\!20$, $T\!=\!5$, and $k=1/3$. The
initial condition is
\begin{align*}
Y_0(x) =&\, \left\{
\begin{array}{ll}
1.2\,\sqrt{3}, &\,x\in [9,11],\\[1mm]
0, &\,\mbox{else}.
\end{array}
\right.
\end{align*}
The solution $Y_{nat}(x,t)$ for $U(x,t)\equiv 0$ describes a natural nucleation process
with wavelike dispersion to the left and right. The distributed control $U(x,t)$ should
now be chosen in such a way that the dispersion is stopped at $t=2.5$, forcing a
stationary profile $Y_{nat}(x,2.5)$ for the remaining time interval. Therefore, we set
 \begin{align*}
Y_Q(x,t) =&\, \left\{
\begin{array}{ll}
Y_{nat}(x,t), &\, t\in [0,2.5],\\[1mm]
Y_{nat}(x,2.5), &\, t\in (2.5,T],
\end{array}
\right.
\end{align*}
in the objective function $J$. There indeed exists an analytic solution for
such a control,
\begin{align*}
U_{stop}(x,t) =&\, \left\{
\begin{array}{ll}
0, &\,t\le 2.5,\\[1mm]
kY^3_{nat}(x,2.5)-Y_{nat}(x,2.5)-\partial_{xx}Y_{nat}(x,2.5), &\,t>2.5,
\end{array}
\right.
\end{align*}
since $\partial_tY(x,t)$ must vanish for $t\ge 2.5$. The authors of
\cite{BuchholzEngelKammannTroeltzsch2013} proved that the second derivative
$\partial_{xx}Y_{nat}(x,2.5)$ is well defined. The functions $Y_{nat}$,
$Y_Q$, and $U_{stop}$ are plotted in Fig.~\ref{fig:nucl_func}.
We again use standard finite differences on the shifted spatial mesh
$x_i=(i-0.5)\triangle x$, $i=1,\ldots,m$, with $\triangle x=L/m$
to discretize the PDE in space.
The objective function is approximated by a linear spline. This yields,
after transforming to the Mayer form, the
discrete control problem
\begin{align*}
\mbox{minimize } \CC\big(y(T)\big) = y_{m+1}(T) &\\
\mbox{subject to } y'(t) =& \,Ay(t) + G(y(t),u(t))=:f(y(t),u(t)),\quad t\in(0,T],\\
y(0) =& \,
\begin{pmatrix}
(Y_0(x_i))_{i=1}^m\\[1mm]
0
\end{pmatrix},
\end{align*}
with
\begin{align*}
(G(y,u))_i = \left\{
\begin{array}{ll}
-ky_i^3+y_i+u_i,&\,i=1,\ldots,m,\\[1mm]
\displaystyle\frac12 ((y-y_Q)_{i=1}^m)\T M(y-y_Q)_{i=1}^m+\frac{\alpha}{2}u^TMu,&\,
i=m+1,
\end{array}
\right.
\end{align*}
and the matrices
\begin{align*}
A = \frac{1}{(\triangle x)^2}
\begin{pmatrix}
-1 & 1 &&&&\\
1 & -2 & 1 &&&\\
&& \ddots &&&\\
&& 1 & -2 & 1 &\\
&&& 1 & -1 &\\
&&&&& 0
\end{pmatrix},
\quad
M = \frac{\triangle x}{12}
\begin{pmatrix}
10 & 2 &&&\\
2 & 8 & 2 &&\\
&& \ddots &&\\
&& 2 & 8 & 2 \\
&&& 2 & 10 \\
\end{pmatrix}.
\end{align*}
Here, $A\in\R^{m+1,m+1}$, $M\in\R^{m,m}$, $u(t)=(U(x_i,t))_{i=1}^m$,
$y_Q(t)=((Y_Q(x_i,t)_{i=1}^m)$ and $y(t)\approx (Y(x_i,t))_{i=1}^{m+1}$.
The total dimension of the discrete control vectors $(U_{nj})$,
$n=0,\ldots,N$, $j=1,\ldots,s$, is $ms(N+1)$. We set $m\!=\!300$ as in
\cite{BuchholzEngelKammannTroeltzsch2013} and note that $s\!=\!4$ for our
Peer methods. The optimal stopping control $U_{stop}$ is discretized by
\begin{align*}
u_{stop}(t) =&\;\left\{
\begin{array}{ll}
0, &\,t\le 2.5,\\[1mm]
ky_{Q}^3(2.5) - y_{Q}(2.5) - \hat{A}y_{Q}(2.5), &\,t>2.5,
\end{array}
\right.
\end{align*}
where $\hat{A}=(A)_{i,j=1}^m$.
With grid sizes $N\in[24,399]$ in time the excessive demand of memory for the full Hessian of the objective function prohibits its use in Matlab's \texttt{fmincon} subroutine.
However, a closer inspection reveals that
\begin{align*}
\nabla_{U_{nj}U_{nj}}\CC =
h\alpha\sum_{i=1}^{s}\kappa_{ij}^{[n]}(P_{ni})_{m+1}M\in\R^{m,m},
\quad n=0,\ldots,N,\;j=1,\ldots,s,
\end{align*}
are the only entries yielding a sparse tridiagonal Hessian,
see Example~\ref{Ectrprb}.
Furthermore, controls $U_{nj}$ with $\kappa_{ij}^{[n]}=0$
for $i=1,\ldots,s,$ are discarded as noted in Chapter \ref{SGlobErr}.
Hence, we pass the sparse Hessian to
\texttt{fmincon} and switch to the \texttt{trust-region-reflective}
algorithm, which allows a simple way for its allocation.
\par
Let us now present the results for the stopping problem and compare them
to those documented in \cite{BuchholzEngelKammannTroeltzsch2013}. There,
the implicit Euler scheme with $h=1/80$, i.e., $400$ uniform time steps,
has been applied, together with a nonlinear cg method and different
step size rules. Using the optimal control $u_{stop}(t)$ given above in
a forward simulation of the ODE, they found
$\CC=3.4814\cdot 10^{-6}$ as reference value. This nicely compares to our
value $\CC=3.1651\cdot 10^{-6}$ for \texttt{AP4o43p} applied with the
same time steps. In principle, the optimizer should find a solution close
to it when started with $U^{(0)}=u_{stop}$ evaluated at the time points
$t_n+c_ih$.
Computation times and values of the objective function are collected in
Table~\ref{tab:nucl_obj}.
Remarkably, already for $N+1=50$ all Peer methods deliver excellent
approximations in very short time compared to $70-100$ seconds reported in
\cite{BuchholzEngelKammannTroeltzsch2013} for similar calculations.
This is a clear advantage of higher-order methods.
\begin{table}[t!]
\begin{center}
\begin{tabular}{|l|r|r|r|r|}\hline
N+1 & 400 & 200 & 100 & 50 \\ \hline
AP4o43p	& 2.99e-6 &	3.24e-6	& 3.91e-6 &	6.02e-6 \\
CPU time [s] &	176	& 51 & 17 &	7 \\ \hline
AP4o33pa	& 3.76e-6 &	6.49e-6 & 8.12e-6 &	1.53e-5 \\
CPU time [s] &	163	& 56 & 27 &	16 \\ \hline
AP4o33pfs	& 4.17e-6 &	5.82e-6 & 1.10e-5 &	2.62e-5 \\
CPU time [s] &	126	& 72 & 30 &	17 \\ \hline
\end{tabular}
\parbox{12.5cm}{
\caption{Values of the objective function $\CC$ and computing times for $U^{(0)}=u_{stop}$ and $N+1=400,200,100,50$ uniform time steps. The reference
value is $\CC=3.1651\cdot 10^{-6}$.
\label{tab:nucl_obj}}
}
\end{center}
\end{table}
\begin{figure}[t!]
\centering
\includegraphics[width=6.8cm]{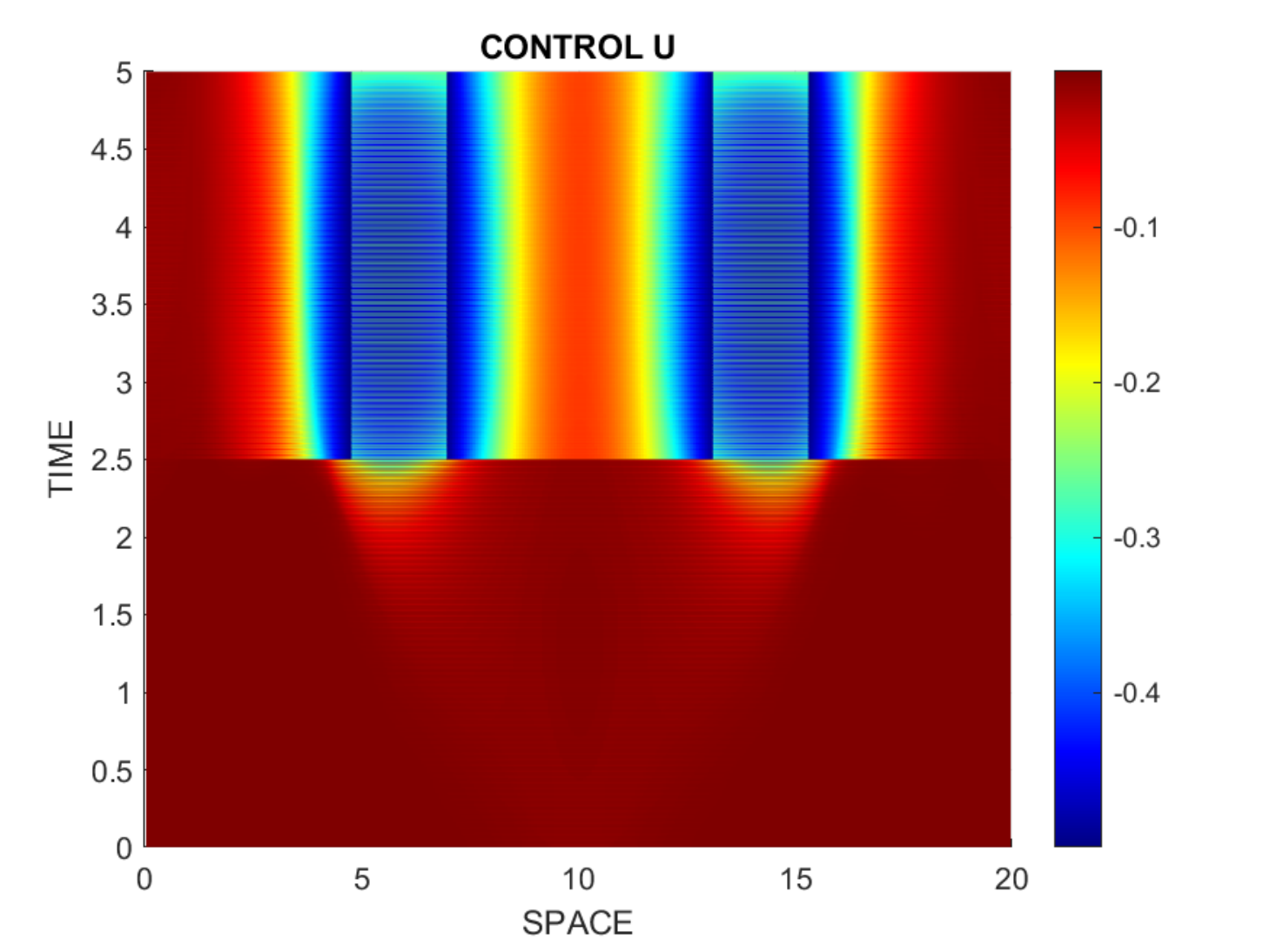}
\hspace{0.1cm}
\includegraphics[width=6.8cm]{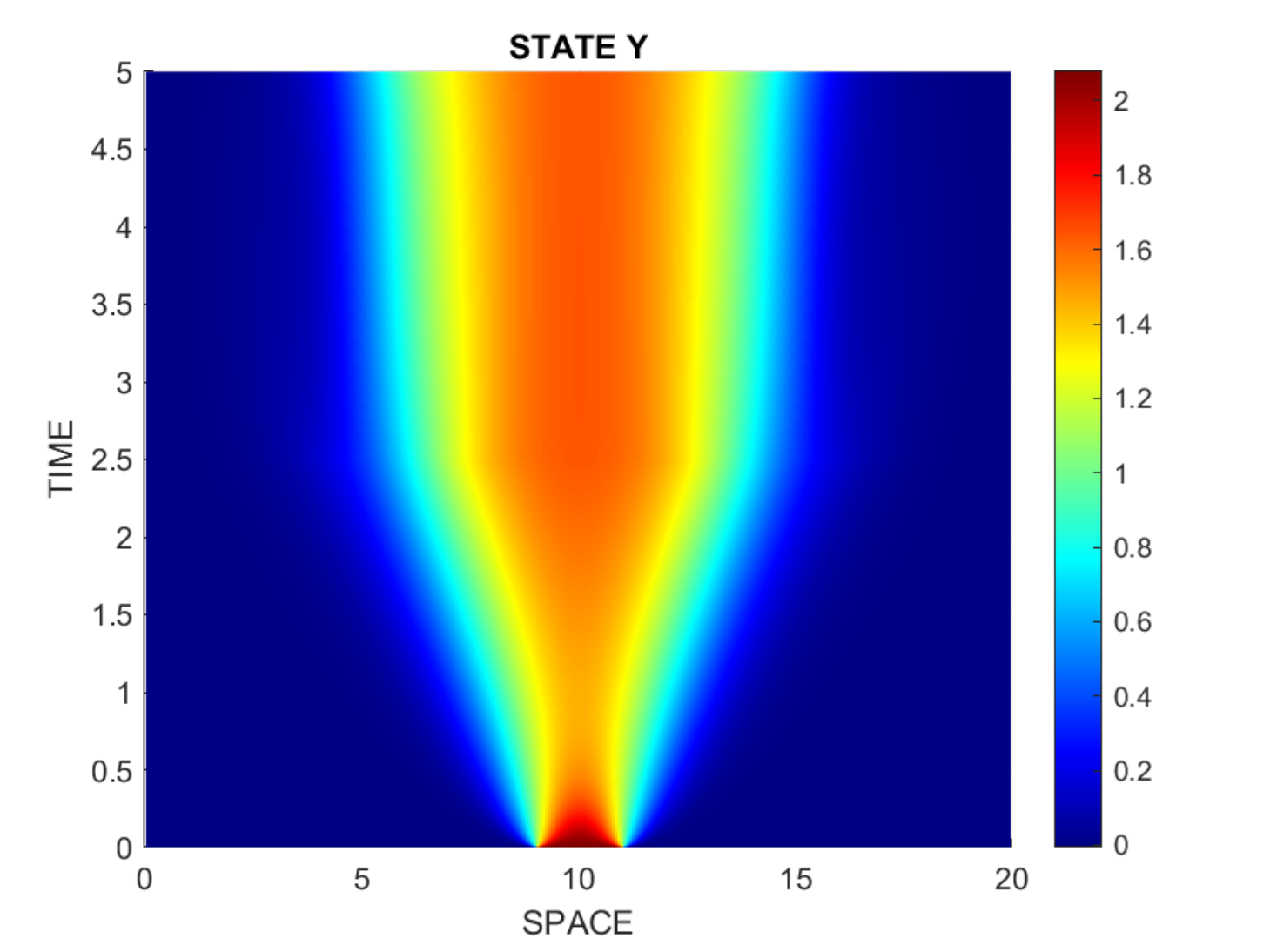}\\[2mm]
\includegraphics[width=6.4cm]{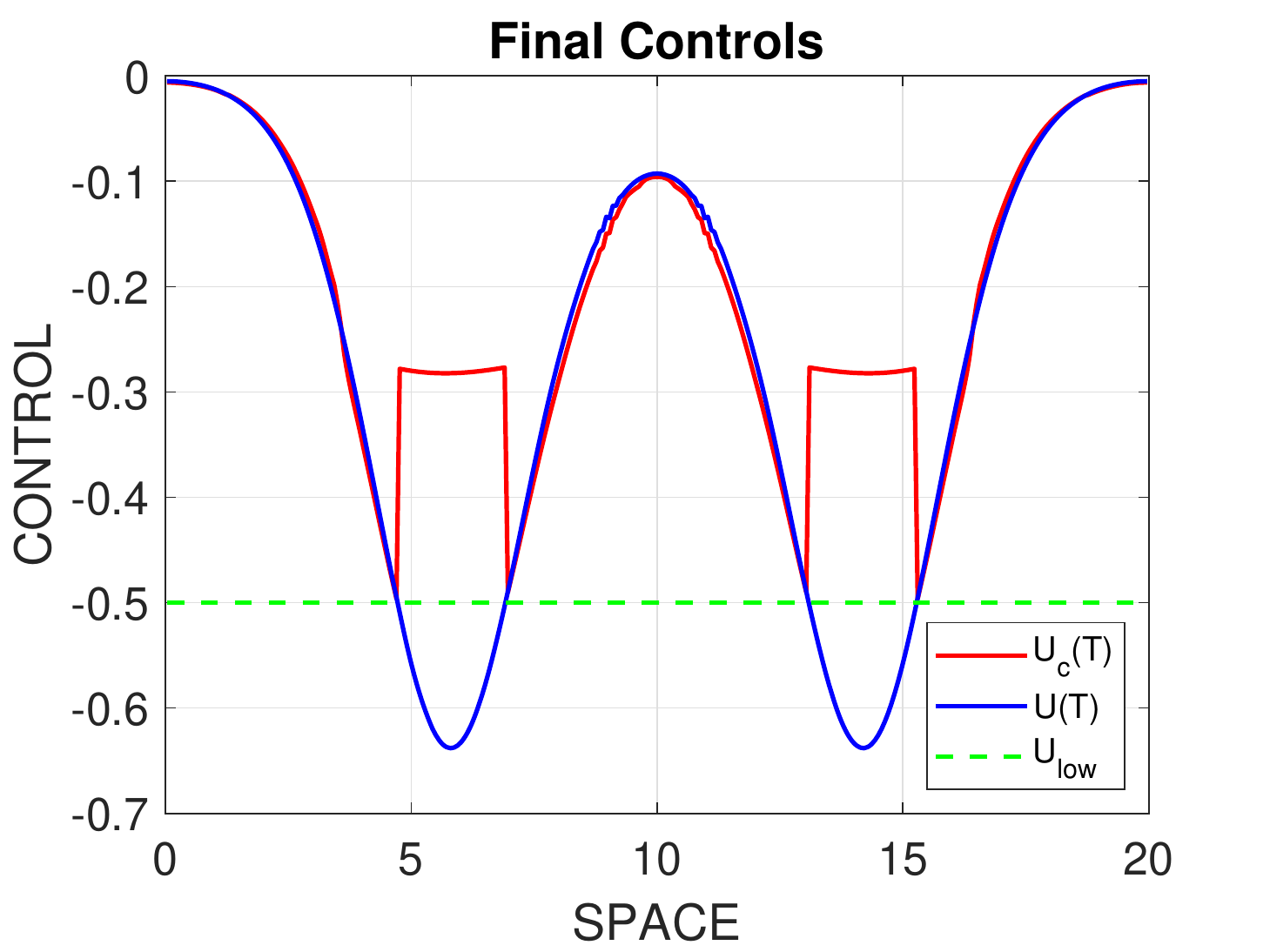}
\hspace{0.1cm}
\includegraphics[width=6.4cm]{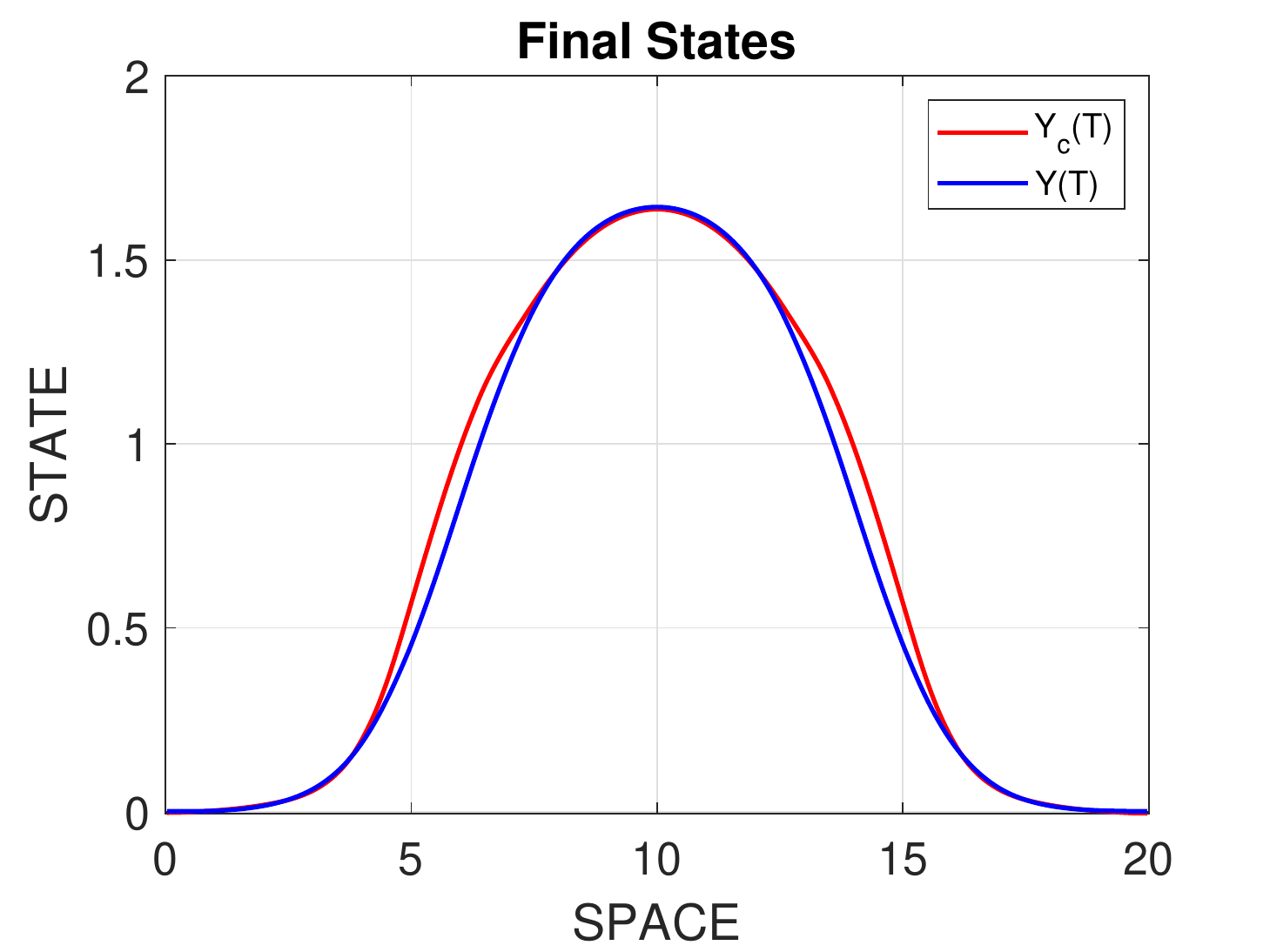}
\parbox{12.5cm}{
\caption{Nucleation process with box constraints for the control. Computed
constrained control $U_c(t)$ (top left), state $Y_h(t)$ (top right),
constraint and unconstraint controls at final time $t=T$
(bottom left), and corresponding final states (bottom right)
approximated by \texttt{AP4o43p} with $50$ uniform time steps.
\label{fig:nucl_cstr}}
}
\end{figure}
Choosing $U^{(0)}=\beta\,u_{stop}$ with $\beta=0.99$,
the authors of \cite{BuchholzEngelKammannTroeltzsch2013}
already discovered slow convergence and tiny deviations from $u_{stop}$ in the
shape of the computed optimal control, which were clearly
visible in their plots. In contrast, all controls computed by the Peer
methods stay close to the overall picture shown in Figure~\ref{fig:nucl_func}
even for $\beta=0.95,0.50$, and for $50$ times steps.
The maximal pointwise control errors range around $6\cdot 10^{-3}$ and $5\cdot 10^{-2}$, respectively.
\par
As a last (speculative) test we impose box constraints of the form
\begin{align*}
u(t) \in U_{ad}:=\{ u(t)\in [L^{\infty}(0,T)]^m
 :\,-0.5\le u_i(t)\le 0,\;i=1,\ldots,m,\;t\in (0,T)\}.
\end{align*}
Now the explicitly given optimal control $U_{stop}$ violates the prescribed
bounds with its minimum value $-0.638$. We apply \texttt{AP4o43p} with $50$ uniform
time steps and set $U^{(0)}=u_{stop}$. \texttt{fmincon} first restricts the control
vector to the admissible set $U_{ad}$ and after $66$ seconds and 142
iteration steps it delivers a solution with $\CC=0.0323$. The stopping process is still
quite satisfactory. Details are plotted in Figure \ref{fig:nucl_cstr}. We get
nearly the same solution for 400 uniform time steps.
Interestingly, the restricted analytic optimal control $\hat{u}_{stop}\in U_{ad}$ only yields $\CC=0.0850$, which is larger than that of the Peer solution by a factor of $2.6$.
\section{Summary}\label{SecSum}
We have upgraded our four-stage implicit Peer triplets constructed in \cite{LangSchmitt2022b}
to meet the additional order conditions and positivity requirements for an efficient use in a gradient-based iterative solution algorithm for ODE constrained optimal control problems. Using
super-convergence for both the state and adjoint variables, 
an almost A-stable method \texttt{AP4o33pa} of the order pair $(3,3)$ with stability angle $\alpha=89.90^o$
could be found. We also considered the class of FSAL methods, where the last stage of the previous time step equals the first stage of the new step, and came up with the
$A(77.53^o)$-stable method \texttt{AP4o33pfs}. A notable theoretical result is that there is no
BDF4-based triplet that improves the second-order control approximation of our recently developed \texttt{AP4o43bdf} \cite{LangSchmitt2022b} to the present setting.
Increasing the order of the forward scheme,
an $A(59.78^o)$-stable method \texttt{AP4o43p} of the order pair $(4,3)$ was constructed. 
All methods show their theoretical orders in the
numerical experiments and clearly outperform the fourth-order symplectic Runge-Kutta-Gauss method in the boundary control problem for an 1D discrete heat equation proposed in
\cite{LangSchmitt2023a} to study order reduction phenomena. The new Peer triplets also perform
remarkably well for a PDE-constrained optimal control problem which models the stopping
of a nucleation process driven by a reaction-diffusion equation of Schl\"ogl-type. In future
work, we will equip our Peer triplets with variable step sizes to further improve their
efficiency.

\vspace{0.5cm}

\par
\noindent\textbf{Funding:}
The first author is supported by the Deutsche Forschungsgemeinschaft
(German Research Foundation) within the collaborative research center
TRR154 {\em ``Mathematical modeling, simulation and optimisation using
the example of gas networks''} (Project-ID 239904186, TRR154/3-2022, TP B01).
\par
\noindent\textbf{Data Availability Statement:} Data will be made available on reasonable request.
\section*{Appendix}
The coefficient matrices which define the Peer triplets \texttt{AP4o33pa},
\texttt{AP4o33pfs} and \texttt{AP4o43p} discussed above are presented here.
We provide exact rational numbers for the node vector $\cc$ and give
numbers with $16$ digits for all matrices.
It is sufficient to only show pairs $(A_n,K_n)$ and
the node vector $\cc=(c_1,\ldots,c_s)\T$ for \texttt{AP4o43p} and some additional data for \texttt{AP4o33pa} and \texttt{AP4o33pfs}, since all other parameters can be easily computed from the following relations:
\begin{align*}
B_n=&\,(A_nV_4-K_nV_4\tilde E_4+R_n)P_4V_4^{-1},\; a=A_0\eins,\; w=A_N\T\eins,
 \;v=V_4\mT e_1,
\end{align*}
with the special matrices
\begin{align*}
V_4=\big(\eins,\cc,\cc^2,\cc^{3}),\quad\PP_4=\Big({j-1\choose i-1}\Big)_{i,j=1}^{4},
\quad\tilde E_4=\big(i\delta_{i+1,j}\big)_{i,j=1}^4\,.
\end{align*}
The matrices $R_0,R_N$ and $R\equiv R_n$ for $1\le n\le N-1$ are slack variables at order 4.
They are provided for \texttt{AP4o33pa}
and \texttt{AP4o33pfs} only since they vanish for the triplet \texttt{AP4o43p} ($R_n\equiv 0$).

{
\small
\everymath{\displaystyle}

\subsection*{\textbf{A1: Coefficients of \texttt{AP4o33pa}}}
\[ \cc\T=\left(\frac{46}{5253},
\frac{29}{51}, \frac{1723}{2193}, \frac{17131}{12189}
\right) \]
\begin{align*}
 A_0=\begin{pmatrix}
 -1.157765450537458& 4.180419822183092&-3.571237514138118& 0.4344668789817266\\[1mm]
  9.320046415868424&-20.43515251977805& 20.53668079758682& -2.660420735071554\\[1mm]
 -9.502446854904932& 18.14294953408145&-17.88837560028214&  2.643706254438956\\[1mm]
  1.573865446847084&-2.968198110625862& 2.201646466132119& 0.1498151692184390 \end{pmatrix}
\end{align*}
\begin{align*}
 K_0=\begin{pmatrix}
  0.1525423728813559& 0.06343283582089552&-0.04424778761061947& 0\\[1mm]
  0.2455414494142291&  0.3479528534959272&  0.2643445483279409& 0\\[1mm]
 -0.2389119757586965&  0.3687279250113433& -0.2354279614690257& 0\\[1mm]
 0.03447092342852595& -0.05320115087852647& 0.03711064142489613&\kappa_{44}^{[0]}
\end{pmatrix}
\end{align*}
with $\kappa_{44}^{[0]}=0.2479535745634692$.
\begin{align*}
 A=\begin{pmatrix}
 0.7073170731707317 & 0 & 0 & 0 \\[1mm]
 -1.458044769359054 & 2.011111111111111 & 0 & 0 \\[1mm]
 0.8963499143698150 &-3.446643123594083 &  2.170212765957447 & 0 \\[1mm]
0.08807733909162651 &0.3555507383436048 &-0.8914986166587666 & 0.5675675675675676
 \end{pmatrix}
\end{align*}
\begin{align*}
 K=\diag\left(0.2240817025504534, 0.2911518627633785, 0.2558139534883721, 0.2289524811977960 \right)
\end{align*}
\begin{align*}
 R=\begin{pmatrix}
 0 & 0 & 0 & -0.2105994034490964\\[1mm]
 0 & 0 & 0 &  0.1876445792137739\\[1mm]
 0 & 0 & 0 & -0.1297946665997080\\[1mm]
 0 & 0 & 0 &  0.1527494908350306
 \end{pmatrix}
\end{align*}
\begin{align*}
 A_N=\begin{pmatrix}
  0.03570841538693515& 0.4969703797836259& 0& 0\\[1mm]
  2.797947998593283&-2.717111089179658&  1.827587054105035&-0.3120359279234260\\[1mm]
 -3.797058467469895& 4.498208855806741& -2.913725127809472& 0.8173416699480771\\[1mm]
  0.4837073832344139&0.1093148794369315&-0.4021296652058669& 0.07527364129327442
 \end{pmatrix}
\end{align*}
\begin{align*}
 K_N=\begin{pmatrix}
  0.2323465386026342& 0.08709000303247828& 0& 0\\[1mm]
  0.0006578497520678987&-0.2800336616814694& 0& 0\\[1mm]
 -0.0006400881985255662& 0.5062443715754399& 0.32694879378132385& 0\\[1mm]
 0.00009235381026342189&-0.07304242875763006& 0&\kappa_{44}^{[N]}
 \end{pmatrix}
\end{align*}
with $\kappa_{44}^{[N]}=0.01004801943170234$.
\begin{align*}
 R_N=\begin{pmatrix}
  0& 0& 0& -0.1751101070505921\\[1mm]
  0& 0& 0&  0.2296022411517165\\[1mm]
  0& 0& 0& -0.5247365005443616\\[1mm]
  0& 0& 0& -0.07622773831802632
 \end{pmatrix}
\end{align*}

\subsection*{\textbf{A2: Coefficients of \texttt{AP4o33pfs}}}
\[ \cc\T=\left(0, \frac{9}{86}, \frac{321}{602}, 1 \right) \]
\begin{align*}
 A_0=\begin{pmatrix}
 1.333333333333333 & 0 & 0 &  0\\[1mm]
-2.789814648187671 & 2.243282202070159 & 0.06686328023669716 &0.01646570267735142\\[1mm]
 4.349477807846901 &-6.391186028966211 & 2.276667661951199   & -0.06058221663260115\\[1mm]
-6.567438826613935 & 9.406667237260441 &-4.671899050533916   &  1.788163545558252
\end{pmatrix}
\end{align*}
\begin{align*}
 K_0=\diag\left(0, 0.2868808051464541, 0.4845433642003949, 0.2814200916147642\right)
\end{align*}
\begin{align*}
 A=\begin{pmatrix}
 0.7857142857142857 & 0 & 0 & 0\\[1mm]
-2.028837530067695  & 2.203900659027200 & 0 & 0\\[1mm]
 4.063000939519495  &-6.340099591541239 & 2.287165301103365 & 0\\[1mm]
-6.494320028787459  & 9.394962342878431 &-4.615533409449387 & 1.744047031603003
 \end{pmatrix}
\end{align*}
\begin{align*}
 K=\diag\left( 0, 0.2754665812532002, 0.4295774647887324, 0.2949559539580673\right)
\end{align*}
\begin{align*}
 R=\begin{pmatrix}
 0& 0& 0&  0\\[1mm]
 0& 0& 0&  0.156340095159149050\\[1mm]
 0& 0& 0& -0.0212049600240154176\\[1mm]
 0& 0& 0& -0.135135135135135135
 \end{pmatrix}
\end{align*}
\begin{align*}
 A_N=\begin{pmatrix}
 1 &  0  & 0 &  0\\[1mm]
-1.037159659693408    &  0.4363577782952090 & 0.6845553714934806 & -0.2064640160522880\\[1mm]
 0.03605110452225963  & -0.5660510638564654 &-0.1074762596776216 & 0.7596425122215622\\[1mm]
 0.001108555171148741 &  0.1296932855612564 &-0.5770791118158589 & 0.4468215038307258
 \end{pmatrix}
\end{align*}
\begin{align*}
 K_N=\begin{pmatrix}
 0.3333333333333333 & 0 & 0 & 0\\[1mm]
-0.3406285072951739 & 0.1264725806602174 & 0 & 0\\[1mm]
 0.1282327493289677 & 0 & 0.5627483658896584 & 0\\[1mm]
-0.03272942952658255 & 0 & 0 & 0.1697266466479663
 \end{pmatrix}
\end{align*}
\begin{align*}
 R_N=\begin{pmatrix}
 0 & 0 & 0 & 0.0463093438915248733\\[1mm]
 0 & 0 & 0 & 0.191797796516481359\\[1mm]
 0 & 0 & 0 & -0.286597642859776972\\[1mm]
 0 & 0 & 0 & 0.1785714285714285754\\[1mm]
 \end{pmatrix}
\end{align*}

\subsection*{\textbf{A3: Coefficients of \texttt{AP4o43p}}}

\[ \cc\T=\left(\frac{4657}{46172}, \frac{43}{97},
\frac{3991}{6596},\frac{21111803999}{23798723875}
\right) \]
\begin{align*}
 A_0=\begin{pmatrix}
  7.666666666666667& -7.952380952380952& 6.428571428571429& -1.0\\[1mm]
 -37.64573385789864& 46.51465022124085& -35.34733224501487& 5.556742966495919\\[1mm]
  38.90401308661976& -51.03310294122830& 39.84674769118604& -5.987622148721481\\[1mm]
  -9.132039686863960& 14.19615134612322&-13.42624214739033& 3.410910572594644 \end{pmatrix}
\end{align*}
\begin{align*}
 K_0=\begin{pmatrix}
  0.2201309814534140& -0.001685331083118719& 0.03214426130560293& 0\\[1mm]
  0.1111845986702137& 0.4311745541022918& -0.1774967804652712& 0\\[1mm]
 -0.1188243074116737& -0.009945644225626329& 0.2279954173163067& 0\\[1mm]
  0.02777498546842700& 0.002324777899894389& -0.04434040826768050&\kappa_{44}^{[0]}
 \end{pmatrix}
\end{align*}
with $\kappa_{44}^{[0]}=0.2883852220354272$.
\begin{align*}
 A=\begin{pmatrix}
  2.080437513028435& 0& 0& 0\\[1mm]
 -6.582767809460944& 2.843481487726957& 0& 0\\[1mm]
  5.640064091163237&-4.381563545251576& 2.010790683327275& 0\\[1mm]
 -1.344827586206897& 3.263399731279439&-4.509045955975008& 1.980031390369082
 \end{pmatrix}
\end{align*}
\begin{align*}
 K=\diag\left( 0.2523093948412364, 0.4504313304404388, 0.0, 0.2972592747183247\right)
\end{align*}
\begin{align*}
 A_N=\begin{pmatrix}
  2.602941176470588& 0.09421300555614037& -1.072906715212599& 0.6\\[1mm]
 -9.770538838886514& 3.643517491998914& 4.765969638829557&-3.172336041397070\\[1mm]
  9.121758438719117&-5.324324324324324&-3.193548387096774& 3.514071174094508\\[1mm]
 -2.137018032260198& 3.217404548657921&-2.956254337680976& 1.067051202531710
  \end{pmatrix}
\end{align*}
\begin{align*}
 K_N=\begin{pmatrix}
  0.2752122060365109& 0& 0.03076923076923077& 0.06493506493506494\\[1mm]
 -0.07088680624623493&\kappa_{22}^{[N]} & -0.1699040256986543& -0.3585636905978095\\[1mm]
  0.07575757575757576& 0& 0.2750926288014159&  0.3832012950339724\\[1mm]
 -0.01770820812361161& 0&-0.04244366487128950& 0.1921737961617600
 \end{pmatrix}
\end{align*}
with $\kappa_{22}^{[N]}=0.3735422712438619$.
}

\bibliographystyle{plain}
\bibliography{bibpeeropt}

\end{document}